\documentclass[11pt,CJK]{article}
\usepackage{amsfonts,amsmath,amssymb,amsthm}
\usepackage{mathrsfs}
\usepackage{latexsym}
\usepackage{graphicx}
\usepackage{fancyhdr,cite}
\usepackage{indentfirst}
\usepackage{color}

\theoremstyle{plain}                       
\newtheorem{lemma}{Lemma}[section]
\newtheorem{thm}[lemma]{Theorem}
\newtheorem{cor}[lemma]{Corollary}
\newtheorem{remark}[lemma]{Remark}
\newtheorem{definition}[lemma]{Definition}

\theoremstyle{remark}
\newtheorem*{pf}{Proof}

\numberwithin{equation}{section}

\def\Xint#1{\mathchoice
  {\XXint\displaystyle\textstyle{#1}}%
  {\XXint\textstyle\scriptstyle{#1}}%
  {\XXint\scriptstyle\scriptscriptstyle{#1}}%
  {\XXint\scriptscriptstyle\scriptscriptstyle{#1}}%
  \!\int}
\def\XXint#1#2#3{{\setbox0=\hbox{$#1{#2#3}{\int}$}
  \vcenter{\hbox{$#2#3$}}\kern-.5\wd0}}

\def\dashint{\Xint-}

\newcommand{\myu}[1]{\tilde{#1}_0}
\newcommand{\myl}[2]{\mathcal{#1}_{#2}}
\newcommand{\mys}[1]{#1_\varepsilon}

\begin{document}

\allowdisplaybreaks
\pagestyle{myheadings}\markboth{$~$ \hfill {\rm Q. Xu,} \hfill $~$} {$~$ \hfill {\rm  } \hfill$~$}
\author{Qiang Xu
\thanks{Email: xuqiang09@lzu.edu.cn.}
\thanks{This work was supported by the Chinese Scholar Council (File No. 201306180043).}
\\
School of Mathematics and Statistics, Lanzhou University, \\
Lanzhou, Gansu 730000, PR China.}

\title{\textbf{Uniform Regularity Estimates in Homogenization Theory of Elliptic Systems with Lower Order Terms on the Neumann Boundary Problem} }
\maketitle
\begin{abstract}
In this paper, we mainly employed the idea of the previous paper \cite{QXS} to
study the sharp uniform $W^{1,p}$ estimates with $1<p\leq \infty$ for
more general elliptic systems with the Neumann boundary condition on a bounded $C^{1,\eta}$ domain, arising in homogenization theory.
Based on the skills developed by Z. Shen in \cite{SZW12} and by T. Suslina in \cite{TS2,TS},
we also established the $L^2$ convergence rates on a bounded $C^{1,1}$ domain and a Lipschitz domain, respectively.
Here we found a ``rough'' version of the first order correctors (see $\eqref{f:4.13}$),
which can unify the proof in \cite{SZW12} and \cite{TS}.
It allows us to skip the corresponding convergence results on $\mathbb{R}^d$ that are the preconditions in \cite{TS2,TS}.
Our results can be regarded as an extension of \cite{SZW4} developed by C. Kenig, F. Lin, Z. Shen,
as well as of \cite{TS} investigated by T. Suslina.
\end{abstract}

\section{Introduction and main results}
M. Avellaneda and F. Lin developed the compactness methods in \cite{MAFHL,MAFHL2}
to study uniform regularity estimates for Dirichlet problems in homogenization theory in the end of 1980s.
For the Neumann boundary value problem, it is not until \cite{SZW4} established by C. Kenig, F. Lin and Z. Shen in 2013 that there was
no significant progress on this topic.
Recently, a new method has been introduced in \cite{SZ,SZW12} by S. Armstrong and Z. Shen to arrive at the sharp regularity estimates,
uniformly down to the microscopic scale, without smoothness assumptions,
for Dirichlet and Neumann problems in periodic or non-periodic settings. Meanwhile T. Suslina \cite{TS2,TS} derived the sharp $O(\varepsilon)$
convergence rate in $L^2(\Omega)$ for elliptic systems with either Dirichlet or Neumann boundary conditions in a $C^{1,1}$ domain.

Inspired by these papers, we originally investigated some uniform regularity estimates
for the elliptic operator with rapidly oscillating potentials that is
\begin{equation*}
 \mathfrak{L}_\varepsilon(u_\varepsilon) = -\Delta u_\varepsilon + \frac{1}{\varepsilon}\mathcal{W}(x/\varepsilon)u_\varepsilon + \lambda u_\varepsilon = F
  \quad \text{in} ~~\Omega,
\end{equation*}
where $\mathcal{W}$ is referred to as the rapidly oscillating potential term (see \cite[pp.93]{ABJLGP}).
As we have shown in \cite{QXS}, the operator $\mathfrak{L}_\varepsilon$ is only a special case of $\myl{L}{\varepsilon}$,
and therefore indicates that our results are not very trivial as it seems to be.

Returning to this paper,
neither the well-known compactness methods nor the new developed technique is rigidly used.
Instead we try to  make full use of the previous work in \cite{SZW4}.
On account of these results, we mainly establish the uniform $W^{1,p}$ estimates $(1<p\leq \infty)$,
as well as the $L^2$ convergence rates
for more general elliptic systems with the Nuemann boundary condition in homogenization theory.
More precisely, we consider the following operators depending on parameter $\varepsilon > 0$,
\begin{eqnarray*}
\mathcal{L}_{\varepsilon} =
-\text{div}\big[A(x/\varepsilon)\nabla +V(x/\varepsilon)\big] + B(x/\varepsilon)\nabla +c(x/\varepsilon) +\lambda I
\end{eqnarray*}
where $\lambda\geq 0$ is a constant, and $I=(e^{\alpha\beta})$ is an identity matrix.

Let $d\geq 3$, $m\geq 1$, and $1 \leq i,j \leq d$ and $1\leq \alpha,\beta\leq m$.
Suppose that $A = (a_{ij}^{\alpha\beta})$, $V=(V_i^{\alpha\beta})$, $B=(B_i^{\alpha\beta})$, $c=(c^{\alpha\beta})$ are real measurable functions,
satisfying the following conditions:
\begin{itemize}
\item the uniform ellipticity condition
\begin{equation}\label{a:1}
 \mu |\xi|^2 \leq a_{ij}^{\alpha\beta}(y)\xi_i^\alpha\xi_j^\beta\leq \mu^{-1} |\xi|^2,
 \quad \text{for}~y\in\mathbb{R}^d,~\text{and}~\xi=(\xi_i^\alpha)\in \mathbb{R}^{md},~\text{where}~ \mu>0;
\end{equation}
 (The summation convention for repeated indices is used throughout.)
\item the periodicity condition
\begin{equation}\label{a:2}
A(y+z) = A(y),~~ V(y+z) = V(y),~~ B(y+z) = B(y),~~ c(y+z) = c(y),
~~\text{for}~y\in \mathbb{R}^d ~\text{and}~ z\in \mathbb{Z}^d;
\end{equation}
\item the boundedness condition
\begin{equation}\label{a:3}
 \max\big\{\|V\|_{L^{\infty}(\mathbb{R}^d)},~\|B\|_{L^{\infty}(\mathbb{R}^d)},~\|c\|_{L^{\infty}(\mathbb{R}^d)}\big\}
 \leq \kappa_1,\qquad \text{where}~\kappa_1>0;
\end{equation}
\item the regularity condition
\begin{equation}\label{a:4}
 \max\big\{ \|A\|_{C^{0,\tau}(\mathbb{R}^d)},~ \|V\|_{C^{0,\tau}(\mathbb{R}^d)}\big\} \leq \kappa_2,
 \qquad \text{where}~\tau\in(0,1)~\text{and}~\kappa_2 > 0.
\end{equation}
\end{itemize}
Set $\kappa = \max\{\kappa_1,\kappa_2\}$, and we say $A\in \Lambda(\mu,\tau,\kappa)$ if $A = A(y)$ satisfies the conditions
$\eqref{a:1}$, $\eqref{a:2}$ and$\eqref{a:4}$.
Throughout this paper, we always assume that $\Omega$ is a bounded $C^{1,\eta}$ domain with $\eta\in [\tau,1)$,
and $L_\varepsilon = -\text{div}[A(x/\varepsilon)\nabla]$ is the elliptic operator from \cite{SZW4}, unless otherwise stated.

We focus on the Neumann boundary value problem:
\begin{equation}\label{pde:1.1}
 \left\{\begin{aligned}
  \mathcal{L}_\varepsilon(u_\varepsilon) & = \text{div}(f) + F  & \qquad \text{in}~~~~\Omega,\\
  \mathcal{B}_\varepsilon(u_\varepsilon)& = g - n\cdot f &\qquad \text{on} ~~\partial\Omega,  
 \end{aligned}\right.
\end{equation}
where
\begin{equation*}
 \big[\mathcal{B}_\varepsilon(u_\varepsilon)\big]^\alpha
 = n_i(x)V_i^{\alpha\beta}(x/\varepsilon)u_\varepsilon^\beta
 + n_i(x)a_{ij}^{\alpha\beta}(x/\varepsilon)\frac{\partial u_\varepsilon^\beta}{\partial x_j}
\end{equation*}
denotes the conormal derivative of $u_\varepsilon$ with respect to $\mathcal{L}_\varepsilon$, and $n = (n_1,\cdots,n_d)$ is the outward unit normal vector to $\partial\Omega$.

The following are the main results of this paper.

\begin{thm}[$W^{1,p}$ estimates]\label{thm:1.1}   
Let $1<p<\infty$. Suppose that $A\in\emph{VMO}(\mathbb{R}^d)$ satisfies $\eqref{a:1}$, $\eqref{a:2}$, and other coefficients satisfy $\eqref{a:3}$.
Let $f\in L^p(\Omega;\mathbb{R}^{md})$, $F\in L^q(\Omega;\mathbb{R}^m)$ and $g\in B^{-1/p,p}(\partial\Omega;\mathbb{R}^m)$,
where $q=\frac{pd}{d+p}$ if $p\geq 2$ and $q=p$ if $1<p<2$.
Assume that $u_\varepsilon\in W^{1,p}(\Omega;\mathbb{R}^m)\cap H^1(\Omega;\mathbb{R}^m)$ is the weak solution to $\eqref{pde:1.1}$,
Then we have the uniform estimate
\begin{equation}\label{pri:1.1}
\|u_\varepsilon\|_{W^{1,p}(\Omega)} \leq C\big\{\|f\|_{L^{p}(\Omega)} + \|F\|_{L^q(\Omega)} + \|g\|_{B^{-1/p,p}(\partial\Omega)}\big\},
\end{equation}
where $C$ depends only on $\mu,\omega,\kappa,\lambda,m,d,p$ and $\Omega$.
\end{thm}

Here we refer the reader to \cite[pp.2283]{SZW10} for the class of $\mathrm{VMO}(\mathbb{R}^d)$,
and $B^{\sigma,p}(\partial\Omega;\mathbb{R}^m)$ denotes the $L^p$ Besov space of order $\sigma$ (see \cite{RAA}).
Compared to the proof of \cite[Theorem 1.1]{SZW4},
(where the solution belongs to the space of the equivalent classes in $H^1(\Omega;\mathbb{R}^m)$
with respect to the relation $u\sim v \Leftrightarrow u-v\in \mathbb{R}^{m}$),
we can not employ Poincar\'e inequality freely any longer. In other words, here we lack the equivalence between
$\|u_\varepsilon\|_{W^{1,p}(\Omega)}$ and $\|\nabla u_\varepsilon\|_{L^{p}(\Omega)}$.
Thus we have to first estimate $\|u_\varepsilon\|_{W^{1,p}(\Omega)}$ for $p\geq2$ by interpolation inequalities,
and then we can infer the quantity of $\|\nabla u_\varepsilon\|_{L^p}$ and $\|u_\varepsilon\|_{L^p}$ for $1<p<2$
by duality, respectively.
We remark that there are no periodicity or regularity assumptions on the coefficients of the lower order terms, and the estimate $\eqref{pri:1.1}$
still holds when $\Omega$ is a bounded $C^1$ domain (see \cite{SZW12}). In recent years,
the uniform $W^{1,p}$ estimates for different type of operators
in homogenization theory have been studied extensively
(see \cite{SZ,MAFHL,LI,GZ,JGZSLS,SGZWS,SZW4,SZWSL,QXS} and their references).


\begin{thm}[Lipschitz estimates]\label{thm:1.2} 
 Suppose that $A\in\Lambda(\mu,\tau,\kappa)$, $V$ satisfies $\eqref{a:2}$ and $\eqref{a:4}$, and other coefficients
satisfy $\eqref{a:3}$.
Let $p>d$ and $0<\sigma<1$. Then, for any $F\in L^p(\Omega;\mathbb{R}^m)$, $f\in C^{0,\sigma}(\Omega;\mathbb{R}^{md})$ and
$g\in C^{0,\sigma}(\partial\Omega;\mathbb{R}^m)$, the weak solution $u_\varepsilon$ to $\eqref{pde:1.1}$
satisfies the uniform estimate
\begin{equation}\label{pri:1.2}
 \|\nabla u_\varepsilon\|_{L^\infty(\Omega)} \leq C \big\{\|f\|_{C^{0,\sigma}(\Omega)}+\|F\|_{L^p(\Omega)}
 +\|g\|_{C^{0,\sigma}(\partial\Omega)}\big\},
\end{equation}
where $C$ depends only on $\mu,\tau,\kappa,\lambda,p,d,m,\eta,$ and $\Omega$.
\end{thm}

We point out that $\eqref{pri:1.2}$ can not be improved even with $C^\infty$ data and domain.
In virtue of the compactness methods, the estimate $\eqref{pri:1.2}$ was established in \cite[Theorem 1.2]{SZW4}
for $L_\varepsilon$ with the Neumann boundary condition under the additional symmetry condition
$A^* = A$, (that is $a_{ij}^{\alpha\beta} = a_{ji}^{\beta\alpha}$),
while this symmetry condition was removed recently in \cite{SZ}.
By means of \cite[Theorem 1.2]{SZW4},
the proof of this theorem can be completed by the method analogous to that used in \cite{QXS}.
We first construct the Neumann boundary corrector associated with $V$ via
\begin{equation}\label{pde:1.2}
\left\{\begin{aligned}
 L_\varepsilon(\Psi_{\varepsilon,0}) & = \text{div}(V_\varepsilon) & \qquad & \text{in} \quad~\Omega, \\
 \frac{\partial\Psi_{\varepsilon,0}}{\partial \nu_\varepsilon} & = n\cdot\big(\widehat{V} - V_\varepsilon\big)
 &\qquad & \text{on}\quad \partial\Omega,
\end{aligned}\right.
\end{equation}
where $\partial/\partial\nu_\varepsilon = n \cdot A(x/\varepsilon)\nabla$, and $\widehat{V}$ is defined in \eqref{f:6.1}. Then
$\Psi_{\varepsilon,0}$ yields the transformation: $u_\varepsilon = \Psi_{\varepsilon,0}v_\varepsilon$ such that $u_\varepsilon$ and $v_\varepsilon$
solve the following equations
\begin{equation*}
(\text{D}_1)\left\{\begin{aligned}
\myl{L}{\varepsilon}(u_\varepsilon) &= \text{div}(f) + F &\quad&\text{in}~~ \Omega,\\
\myl{B}{\varepsilon}(u_\varepsilon) &= g -n\cdot f                &\quad&\text{on}~ \partial\Omega,
\end{aligned}\right.
\qquad \text{and}\qquad
(\text{D}_2)
\left\{\begin{aligned}
L_\varepsilon(v_\varepsilon) &= \text{div}(\tilde{f}) + \tilde{F} &\quad&\text{in}~~\Omega,\\
\frac{\partial v_\varepsilon}{\partial\nu_\varepsilon}&= \tilde{g} -n\cdot\tilde{f}    &\quad&\text{on}~\partial\Omega,
\end{aligned}\right.
\end{equation*}
respectively. Then the result of \cite[Theorem 1.2]{SZW4} can be directly applied to $(\text{D}_2)$.
Due to the same reason as explained in \cite{QXS}, we need to derive
\begin{equation*}
\|\Psi_{\varepsilon,0}-I\|_{L^\infty(\Omega)} \leq C\varepsilon \ln(r_0/\varepsilon+2),
\qquad r_0 = \text{diam}(\Omega),
\end{equation*}
which follows from the decay estimates of Neumann matrixes defined in \cite[pp.916]{SZW4} (see Theorem $\ref{thm:3.1}$),
as well as
\begin{equation*}
\|\nabla \Psi_{\varepsilon,0}\|_{C^{0,\sigma_1}(\Omega)} = O(\varepsilon^{-\sigma_1})
\qquad \text{and}\qquad
\|\nabla u_\varepsilon\|_{C^{0,\sigma_1}(\Omega)} = O(\varepsilon^{-\sigma_2}) \quad\text{as}~~ \varepsilon \to 0,
\end{equation*}
which are the main conclusions of Corollary $\ref{cor:3.1}$ and Lemma $\ref{lemma:5.7}$,
where $0<\sigma_1<\sigma_2<1$ are independent of $\varepsilon$.
The above two estimates together with $\eqref{pri:1.1}$ guarantee
that the right hand side of $(\text{D}_2)$ can be uniformly bounded by the given data in Theorem $\ref{thm:1.2}$.

Note that the right hand side of $(\text{D}_2)$ which involves $\text{div}(\tilde{f})$
is, as a matter of fact, more general than that in \cite[Theorem 1.2]{SZW4}. We find a simple way inspired by \cite{SZWSL} to
derive the uniform Lipschitz estimate for the weak solution $u_\varepsilon$ to $L_\varepsilon(u_\varepsilon) = \text{div}(f)$ in
$\Omega$ and $\partial u_\varepsilon/\partial \nu_\varepsilon = -n\cdot f $ on $\partial\Omega$. The key ingredient is to construct
the auxiliary functions $\{v_{\varepsilon,k}\}_{k=1}^d$, which satisfy
 \begin{equation*}
  L_\varepsilon(v_{\varepsilon,k}) = 0 \quad\text{in}~~\Omega,
  \qquad \partial v_{\varepsilon,k}/\partial \nu_\varepsilon = -n_k I \quad\text{on}~~\partial\Omega
 \end{equation*}
 for $k=1,\cdots,d$, where $n_k$ is the $k$'th component of the outward unit normal to $\partial\Omega$. This is of help
 to eliminating the bad term $\int_{\partial\Omega}\nabla_x N_\varepsilon(x,z)n_k(z)dS(z)$ in $\eqref{f:3.3}$
 (see the proof of Lemma $\ref{lemma:3.3}$), where $N_\varepsilon$ denotes the Neumann matrixes associated with $L_\varepsilon$.
 Then the desired result: $\|\nabla u_{\varepsilon}\|_{L^\infty(\Omega)}\leq C\|f\|_{C^{0,\sigma}(\Omega)}$
 follows from \cite[Theorem 1.2]{SZW4} at once.
 We remark that this argument also works for common elliptic operators,
 provided that there are some decay estimates of corresponding Neumann matrixes previously. At the end of the paragraph, we mention that both
 the compactness method in \cite{MAFHL,MAFHL2,SZW4}  and the new argument developed in \cite{SZ,SZW12}
 may be still valid for the Neumann problem $\eqref{pde:1.1}$, however it will be quite complicated as compared with our method.
 For more references on this topic, see \cite{SZ,MAFHL,MAFHL2,GZ,SGZWS,SZW4,CC,SZW12,SZWSL,QXS}.

\begin{thm}[$L^2$ convergence rates]\label{thm:1.3}
Suppose that the coefficients of $\mathcal{L}_\varepsilon$ satisfies $\eqref{a:1}$, $\eqref{a:2}$ and $\eqref{a:3}$.
Assume that $u_\varepsilon,u_0$  are the weak solutions to
\begin{equation}\label{pde:1.5}
(\emph{H}_\varepsilon)~\left\{\begin{aligned}
\mathcal{L}_\varepsilon(u_\varepsilon) & = F  &\qquad& \emph{in}~~\Omega,\\
\mathcal{B}_\varepsilon(u_\varepsilon) & = g  &\qquad& \emph{on}~\partial\Omega,
\end{aligned}\right.
\qquad\qquad
(\emph{H}_0)~\left\{\begin{aligned}
\mathcal{L}_0(u_0) & = F  &\qquad& \emph{in}~~\Omega,\\
\mathcal{B}_0(u_0) & = g  &\qquad& \emph{on}~\partial\Omega,
\end{aligned}\right.
\end{equation}
respectively. We obtain the following two results:
\begin{itemize}
\item[\emph{(i)}] if $\Omega$ is a bounded $C^{1,1}$ domain, and
$F\in L^2(\Omega;\mathbb{R}^m)$ and $g\in B^{1/2,2}(\partial\Omega;\mathbb{R}^m)$ are given,
then
\begin{equation}\label{pri:1.3}
\|u_\varepsilon - u_0\|_{L^2(\Omega)} \leq C\varepsilon\big\{\|F\|_{L^2(\Omega)}+\|g\|_{B^{1/2,2}(\partial\Omega)}\big\};
\end{equation}
\item[\emph{(ii)}] let $\Omega$ be a bounded Lipschitz domain. Suppose that
$F\in L^\frac{2d}{d+1}(\Omega;\mathbb{R}^m)$ and $g\in L^{2}(\partial\Omega;\mathbb{R}^m)$, and $A$ additionally satisfies $A^*=A$, then
we have
\begin{equation}\label{pri:1.4}
\|u_\varepsilon - u_0\|_{L^2(\Omega)} \leq C\varepsilon^{\frac{1}{2}}\big\{\|F\|_{L^{\frac{2d}{d+1}}(\Omega)}+\|g\|_{L^{2}(\partial\Omega)}\big\},
\end{equation}
\end{itemize}
Moreover, if $u_0\in H^2(\Omega)$, then $\|u_\varepsilon - u_0\|_{L^{\frac{2d}{d-1}}(\Omega)}\leq C\varepsilon\|u_0\|_{H^2(\Omega)}$,
where $C$ depends only on $\mu,\kappa,m,d$ and $\Omega$.

%
%
\end{thm}

We note that $B^{1/2,2}(\partial\Omega;\mathbb{R}^m) = H^{1/2}(\partial\Omega;\mathbb{R}^m)$.
Here $\myl{L}{0}$ and $\myl{B}{0}$ from the homogenized equation $(\text{H}_0)$ are defined in $\eqref{pde:6.3}$.
As mentioned before, we find a new type of the first order corrector $\varepsilon\chi_{k}(x/\varepsilon)\varphi_k$ with $k=0,1,\ldots,d$,
which together leads to
\begin{equation}\label{f:4.13}
w_\varepsilon^\beta = u_\varepsilon^\beta - u_0^\beta
-\varepsilon\sum_{k=0}^d\chi_k^{\beta\gamma}(x/\varepsilon)\varphi_k^\gamma,
\end{equation}
where
$\varphi=(\varphi_k^\gamma)\in H^2(\mathbb{R}^d;\mathbb{R}^{md})$, and $\{\chi_{k}\}_{k=0}^d$
are correctors associated with $\myl{L}{\varepsilon}$ (see \eqref{pde:6.1} and \eqref{pde:6.2}).
Then we acquire the ``first but rough'' estimate of
$\|w_\varepsilon\|_{H^1(\Omega)}$ in Lemma $\ref{lemma:4.4}$ by energy methods
(where we borrow the idea from T. Suslina in \cite[pp.3475-3477]{TS}),
and the next step is therefore reduced to choose suitable $\varphi_k$ to obtain the convergence rates under the different conditions.
In author's point of view, Lemma $\ref{lemma:4.4}$ plays a role as the bifurcation in the whole proof of Theorem $\ref{thm:1.3}$.

Precisely speaking, when $\Omega$ is a bounded $C^{1,1}$ domain, following the arguments developed by T. Suslina \cite{TS2,TS},
it is not hard to derive $\|w_\varepsilon\|_{H^1(\Omega)}=O(\varepsilon^{\frac{1}{2}})$
by choosing $\varphi_0 = \bar{S}_\varepsilon(\tilde{u}_0)$
and $\varphi_k = \bar{S}_\varepsilon(\nabla_k\tilde{u}_0)$ (see Remark $\ref{rm:5.1}$).
But we fail to obtain
$\|w_\varepsilon\|_{L^2(\Omega)}=O(\varepsilon)$ due to the bad estimate $\eqref{pri:4.9}$ for the boundary terms.
So, we turn to construct $w_\varepsilon$ by choosing $\varphi_0 = S_\varepsilon(\psi_{4\varepsilon}u_0)$ and
$\varphi_k = S_\varepsilon(\psi_{4\varepsilon}\nabla_ku_0)$ to avoid handling the boundary things (see Lemma $\ref{lemma:4.5}$),
then by duality argument, we can arrive at the estimate $\eqref{pri:1.3}$ (see Lemma $\ref{lemma:6.5}$,
we mention that our proof actually relies on the new development in \cite{SZW12}).

For a bounded Lipschitz domain, due to the skills improved by Z. Shen in \cite{SZW12},
we set $\varphi_0 = S_\varepsilon^2(\psi_{2\varepsilon}u_0)$
and $\varphi_k = S_\varepsilon^2(\psi_{2\varepsilon}\nabla_k u_0)$ (see Lemma $\ref{lemma:4.2}$), and similarly
obtain $\|w_\varepsilon\|_{H^1(\Omega)} = O(\varepsilon^{\frac{1}{2}})$, which leads to the estimate
$\eqref{pri:1.4}$.
The progress is that we additionally employ the radial maximal function
coupled with non-tangential maximal function (see $\eqref{def:2}$ and $\eqref{def:3}$)
to analysis the boundary behavior of the solution $u_0$ to ($\text{H}_0$).
We remark that the thinking is originally arose by C. Kenig, F. Lin and Z. Shen in \cite{SZW2}.
If $u\in H^2(\Omega;\mathbb{R}^m)$, using duality method again, we can derive the sharp estimate
$\|w_\varepsilon\|_{L^{\frac{2d}{d-1}}(\Omega)} = O(\varepsilon)$.
Although $\eqref{pri:1.4}$ is not sharp, it opens a new door to reach the Rellich estimate (see Remark $\ref{re:5.2}$).
We mention that under some conditions, it is possible to remove the assumption of $u_0\in H^2(\Omega;\mathbb{R}^m)$ by
using a subtle technique,
and we will show that in another paper.

Here $S_\varepsilon,\bar{S}_\varepsilon$ and $\psi_{2\varepsilon}$ are defined in $\eqref{f:6.5}$, $\eqref{def:6.3}$ and $\eqref{f:4.20}$, respectively.
$\tilde{u}_0$ is an extension function of $u_\varepsilon$ on $\mathbb{R}^d$.
Thus to varying degrees,
the two results with different preconditions between \cite{TS} and \cite{SZW2}
may be reduced by Lemma $\ref{lemma:4.4}$ to figure out some proper first order correctors.
Note that both $\eqref{pri:1.3}$ and $\eqref{pri:1.4}$ do not require the smoothness assumptions on the coefficients
of $\myl{L}{\varepsilon}$,  where the estimate $\eqref{pri:1.3}$ is sharp.

We also remark that the counterpart of Theorem $\ref{thm:1.3}$ in \cite{SZ,SZW12}
performs a crucial role as the start point for various sorts of uniform regularity estimates,
which marks the new way to regularity theory of homogenization problems.
We refer the reader to \cite{SZ,ABJLGP,BMSHSTA,BMSHSTA2,CDJ,GG1,GG2,SGZWS,VSO,SZW2,SZW1,ODVB,SZW12,TS2,TS,SAT2,QXS} and references therein
for more results. The assumption of $d\geq 3$ is not essential but easy to organize the paper,
since we usually have different way to handle the corresponding problem in the case of $d=2$.

This paper is organized as follows. We introduce some definitions, symbols and remarks in Section 2.
The proof of Theorem $\ref{thm:1.1}$ is shown in Section 3, and the uniform Lipschitz estimate is studied in Section 4.
Finally, we discuss the $L^2$ convergence rates in Section 5.

\section{Preliminaries}

Define the correctors $\chi_k = (\chi_{k}^{\alpha\beta})$, $0\leq k\leq d$, associated with $\mathcal{L}_\varepsilon$ as follows:
\begin{equation}\label{pde:6.1}
\left\{ \begin{aligned}
 &L_1(\chi_k) = \text{div}(V)  \quad \text{in}~ \mathbb{R}^d, \\
 &\chi_k\in H^1_{per}(Y;\mathbb{R}^{m^2})~~\text{and}~\int_Y\chi_k dy = 0
\end{aligned}
\right.
\end{equation}
for $k=0$, and
\begin{equation}\label{pde:6.2}
 \left\{ \begin{aligned}
  &L_1(\chi_k^\beta + P_k^\beta) = 0 \quad \text{in}~ \mathbb{R}^d, \\
  &\chi_k^\beta \in H^1_{per}(Y;\mathbb{R}^m)~~\text{and}~\int_Y\chi_k^\beta dy = 0
 \end{aligned}
 \right.
\end{equation}
for $1\leq k\leq d$, where $Y = (0,1]^d \cong \mathbb{R}^d/\mathbb{Z}^d$, and $H^1_{per}(Y;\mathbb{R}^m)$ denotes the closure
of $C^\infty_{per}(Y;\mathbb{R}^m)$ in $H^1(Y;\mathbb{R}^m)$.
Note that $C^\infty_{per}(Y;\mathbb{R}^m)$ is the subset of $C^\infty(Y;\mathbb{R}^m)$, which collects all $Y$-periodic vector-valued functions
(see \cite[pp.56]{ACPD}). By asymptotic expansion arguments, we obtain the homogenized operator
\begin{equation*}
 \mathcal{L}_0 = -\text{div}(\widehat{A}\nabla+ \widehat{V}) + \widehat{B}\nabla + \widehat{c} + \lambda I,
\end{equation*}
where $\widehat{A} = (\hat{a}_{ij}^{\alpha\beta})$, $\widehat{V}=(\hat{V}_i^{\alpha\beta})$,
$\widehat{B} = (\hat{B}_i^{\alpha\beta})$ and $\widehat{c}= (\hat{c}^{\alpha\beta})$ are given by
\begin{equation}\label{f:6.1}
\begin{aligned}
\hat{a}_{ij}^{\alpha\beta} = \int_Y \big[a_{ij}^{\alpha\beta} + a_{ik}^{\alpha\gamma}\frac{\partial\chi_j^{\gamma\beta}}{\partial y_k}\big] dy, \qquad
\hat{V}_i^{\alpha\beta} = \int_Y \big[V_i^{\alpha\beta} + a_{ij}^{\alpha\gamma}\frac{\partial\chi_0^{\gamma\beta}}{\partial y_j}\big] dy, \\
\hat{B}_i^{\alpha\beta} = \int_Y \big[B_i^{\alpha\beta} + B_j^{\alpha\gamma}\frac{\partial\chi_i^{\gamma\beta}}{\partial y_j}\big] dy, \qquad
\hat{c}^{\alpha\beta} = \int_Y \big[c^{\alpha\beta} + B_i^{\alpha\gamma}\frac{\partial\chi_0^{\gamma\beta}}{\partial y_i}\big] dy.
\end{aligned}
\end{equation}

\begin{definition}\label{def:1}
\emph{Let $f=(f_i^\alpha)\in L^2(\Omega;\mathbb{R}^{md})$,
$F\in L^{\frac{2d}{d+2}}(\Omega;\mathbb{R}^m)$ and $g\in B^{-1/2,2}(\partial\Omega;\mathbb{R}^m)$.
We say $u_\varepsilon\in H^1(\Omega;\mathbb{R}^m)$ is a weak solution to $\eqref{pde:1.1}$,
if $u_\varepsilon$ satisfies
\begin{eqnarray}\label{pde:1.4}
 \mathrm{B}_\varepsilon[u_\varepsilon,\phi] = -\int_\Omega f^\alpha\cdot\nabla\phi^\alpha dx + \int_\Omega F^\alpha\phi^{\alpha} dx
 + <g,\phi>_{B^{-1/2,2}(\partial\Omega)\times B^{1/2,2}(\partial\Omega)}
\end{eqnarray}
for any $\phi\in H^1(\Omega;\mathbb{R}^m)$, where
\begin{equation}\label{pde:1.6}
 \mathrm{B}_\varepsilon[u_\varepsilon,\phi]
 = \int_\Omega \Big\{a_{ij,\varepsilon}^{\alpha\beta}\frac{\partial u_\varepsilon^\beta}{\partial x_j}
 + V_{i,\varepsilon}^{\alpha\beta}u_\varepsilon^\beta\Big\}\frac{\partial \phi^\alpha}{\partial x_i} dx
 + \int_\Omega \Big\{B_{i,\varepsilon}^{\alpha\beta}\frac{\partial u_\varepsilon^\beta}{\partial x_i}
  + c_\varepsilon^{\alpha\beta} u_\varepsilon^\beta + \lambda u_\varepsilon^\alpha \Big\}\phi^\alpha dx
\end{equation}
is the bilinear form associated with $\mathcal{L}_\varepsilon$.}
\end{definition}

\begin{remark}
\emph{Choose $\phi^\alpha = 1$ in $\eqref{pde:1.4}$, and then we have the compatibility condition
\begin{equation}\label{C:1}
 \int_\Omega \Big(B_i^{\alpha\beta}(x/\varepsilon)\frac{\partial u_\varepsilon^\beta}{\partial x_i}
 + c^{\alpha\beta}(x/\varepsilon) u_\varepsilon^\beta \Big)dx + \lambda\int_\Omega u_\varepsilon^\alpha dx
 =\int_\Omega F^\alpha dx +<g^\alpha,1>_{B^{-1/2,2}(\partial\Omega)\times B^{1/2,2}(\partial\Omega)}
\end{equation}
for $\alpha = 1,\ldots,m$, which implies the counterpart of $\eqref{C:1}$ in \cite{SZW4} since $B=0, c=0$ and $\lambda=0$ there.}
\end{remark}

\begin{remark}
We similarly define the bilinear form $\mathrm{B}_0$ associated with $\mathcal{L}_0$ as
\begin{equation}
 \mathrm{B}_0[u_0,\phi]
 = \int_\Omega \Big\{\hat{a}_{ij}^{\alpha\beta}\frac{\partial u_0^\beta}{\partial x_j}
 + \hat{V}_{i}^{\alpha\beta}u_0^\beta\Big\}\frac{\partial \phi^\alpha}{\partial x_i} dx
 + \int_\Omega \Big\{\hat{B}_{i}^{\alpha\beta}\frac{\partial u_0^\beta}{\partial x_i}
  + \hat{c}^{\alpha\beta} u_0^\beta + \lambda u_0^\alpha \Big\}\phi^\alpha dx
\end{equation}
for any $u_0,\phi\in H^1(\Omega;\mathbb{R}^m)$.
\end{remark}

We remark that not only the solutions $u_\varepsilon$ approaches to $u_0$ weakly in $H^1(\Omega;\mathbb{R}^m)$,
but also the flows converge, i.e.
$A_\varepsilon\nabla u_\varepsilon + V_\varepsilon u_\varepsilon \rightharpoonup \widehat{A}\nabla u_0 + \widehat{V}u_0$, and
$B_\varepsilon\nabla u_\varepsilon + (c_\varepsilon+\lambda I)u_\varepsilon \rightharpoonup
\widehat{B}\nabla u_0 + (\widehat{c}+\lambda I)u_0$
weakly in $L^2(\Omega;\mathbb{R}^m)$ as $\varepsilon\to 0$ (see \cite[pp.31]{VSO}).
In other words, for any $\phi\in H^1(\Omega;\mathbb{R}^m)$
we have $\mathrm{B}_\varepsilon[u_\varepsilon,\phi]\to \mathrm{B}_0[u_0,\phi]$ as $\varepsilon\to 0$.
Thus by Definition $\ref{def:1}$, we have
\begin{equation}\label{pde:1.7}
 \mathrm{B}_0[u_0,\phi] = -\int_\Omega f^\alpha\cdot\nabla\phi^\alpha dx + \int_\Omega F^\alpha\phi^{\alpha} dx
 + <g,\phi>_{B^{-1/2,2}(\partial\Omega)\times B^{1/2,2}(\partial\Omega)}.
\end{equation}
Let $\phi\in H^1_0(\Omega;\mathbb{R}^m)$ be an arbitrary test function in $\eqref{pde:1.7}$, and integrating by parts we can derive
$\mathcal{L}_0(u) = \text{div}(f) + F$ in $\Omega$.
Then we return to $\eqref{pde:1.7}$ with any $\phi\in H^1(\Omega;\mathbb{R}^m)$,
and obtain $n\cdot\big(\widehat{A}\nabla u_0 + \widehat{V}u_0\big) = g$ on $\partial\Omega$. By writing
$\myl{B}{0}(u_0) = n\cdot\big(\widehat{A}\nabla u_0 + \widehat{V} u_0\big)$,
we can express the homogenized problem related to the Neumann problem $\eqref{pde:1.1}$ as following:
\begin{equation}\label{pde:6.3}
\left\{\begin{aligned}
\myl{L}{0}(u_0) & = \text{div}(f) + F &\qquad&\text{in}~~~\Omega,\\
\myl{B}{0}(u_0) & = g &\qquad&\text{on}~~\partial\Omega.
\end{aligned}\right.
\end{equation}
We refer the reader to \cite[pp.103]{ABJLGP} or \cite[pp.31]{VSO} for more comments.
We also point out that $\myl{L}{0}$ is still an elliptic operator and if $A^* = A$, then we may have
$\mu |\xi|^2 \leq \hat{a}_{ij}^{\alpha\beta}\xi_i^\alpha\xi_j^\beta\leq \mu^{-1}|\xi|^2$ for any $\xi=(\xi_i^\alpha)\in\mathbb{R}^{md}$
(see \cite[pp.23]{ABJLGP}).

\begin{definition}
\emph{Define the adjoint operator $\mathcal{L}_\varepsilon^*$ as
\begin{equation*}
 \mathcal{L}_\varepsilon^*
 = -\text{div}\Big\{A^*(x/\varepsilon)\nabla + B^{*}(x/\varepsilon)\Big\}
 + V^{*}(x/\varepsilon)\nabla
 + c^{*}(x/\varepsilon) + \lambda I,
\end{equation*}
while the corresponding boundary operator becomes
\begin{equation*}
 \mathcal{B}_\varepsilon^* = n\cdot \big[A^{*}(x/\varepsilon)\nabla+B^*(x/\varepsilon)\big].
\end{equation*}
Furthermore, the related bilinear form is given by
\begin{equation*}
 \mathrm{B}_\varepsilon^*[v_\varepsilon,\phi]
 = \int_\Omega \Big\{a_{ij,\varepsilon}^{\alpha\beta}\frac{\partial v_\varepsilon^\alpha}{\partial x_i}
 + B_{j,\varepsilon}^{\alpha\beta} v_\varepsilon^\alpha\Big\} \frac{\partial \phi^\beta}{\partial x_j} dx
 + \int_\Omega \Big\{V_{i,\varepsilon}^{\alpha\beta}\frac{\partial v_\varepsilon^\beta}{\partial x_i}
  + c_\varepsilon^{\alpha\beta} v_\varepsilon^\alpha + \lambda v_\varepsilon^\beta \Big\} \phi^\beta dx
\end{equation*}
for any $v_\varepsilon,\phi\in H^1(\Omega;\mathbb{R}^m)$.}
\end{definition}

\begin{remark}
\emph{ If $u_\varepsilon,v_\varepsilon\in H^1(\Omega;\mathbb{R}^m)$ are two weak solutions to
\begin{equation*}
\left\{\begin{aligned}
 \mathcal{L}_\varepsilon(u_\varepsilon) & = \text{div}(f)+F &\quad \text{in}~~~\Omega,\\
 \mathcal{B}_\varepsilon(u_\varepsilon) & = g-n\cdot f      &\quad \text{on}~\partial\Omega,
\end{aligned}\right.
\qquad\quad
\left\{\begin{aligned}
\mathcal{L}_\varepsilon^*(v_\varepsilon) &= \text{div}(h) + H &\quad \text{in}~~~\Omega,\\
\mathcal{B}_\varepsilon^*(v_\varepsilon) &= b -n\cdot h &\quad \text{on}~\partial\Omega,
\end{aligned}\right.
\end{equation*}
respectively, then we have the second Green's formula
\begin{equation}\label{G:1}
<\mathcal{L}_\varepsilon(u_\varepsilon),v_\varepsilon> - <u_\varepsilon,\mathcal{L}_\varepsilon^*(v_\varepsilon)>
= -<\mathcal{B}_\varepsilon(u_\varepsilon),v_\varepsilon> + <u_\varepsilon,\mathcal{B}_\varepsilon^*(v_\varepsilon)>,
\end{equation}
by noting that $\mathrm{B}_\varepsilon[u_\varepsilon,v_\varepsilon] = \mathrm{B}_\varepsilon^*[v_\varepsilon,u_\varepsilon]$.
Moreover, if $f_i,F,H\in L^2(\Omega;\mathbb{R}^m)$ with $i=1,\cdots,d$, and $h,b,g$ vanish, then we have
\begin{equation}\label{f:2.10}
\int_\Omega u_\varepsilon H dx
= -\int_\Omega f\cdot\nabla v_\varepsilon dx + \int_\Omega Fv_\varepsilon dx.
\end{equation} }
\end{remark}

\begin{remark}
\emph{To handle the convergence rates, we define some auxiliary functions via
\begin{equation}\label{def:6.1}
b_{ik}^{\alpha\gamma}(y) = \hat{a}_{ik}^{\alpha\gamma} - a_{ik}^{\alpha\gamma}(y)
- a_{ij}^{\alpha\beta}(y)\frac{\partial}{\partial y_j}\big\{\chi_{k}^{\beta\gamma}\big\}, \qquad
b_{i0}^{\alpha\gamma}(y) =
\hat{V}_i^{\alpha\gamma} - V_{i}^{\alpha\gamma}(y)-a_{ij}^{\alpha\beta}(y)\frac{\partial}{\partial y_j}\big\{\chi_0^{\beta\gamma}\big\},
\end{equation}
and
\begin{equation}\label{def:6.2}
\left.\begin{aligned}
\Delta\vartheta_{i}^{\alpha\gamma}
 & = W_i^{\alpha\gamma} := \hat{B}_i^{\alpha\gamma} - B_{i}^{\alpha\gamma}(y)
- B_{j}^{\alpha\beta}(y)\frac{\partial}{\partial y_j}\big\{\chi_i^{\beta\gamma}\big\} & \quad \text{in} ~~\mathbb{R}^d,
&\qquad \int_Y \vartheta_i^{\alpha\beta}(y) dy = 0 ,\\
\Delta \vartheta_0^{\alpha\gamma}
&= W_0^{\alpha\gamma} := \hat{c}^{\alpha\gamma} - c^{\alpha\gamma}(y)
- B_{i}^{\alpha\beta}(y)\frac{\partial}{\partial y_i}\big\{\chi_0^{\beta\gamma}\big\} & \quad \text{in}~~\mathbb{R}^d ,
&\qquad \int_Y \vartheta_0^{\alpha\beta}(y) dy = 0 .
\end{aligned}\right.
\end{equation}
We mention that the existence of $\vartheta_{k}$ is given by \cite[Theorem 4.28]{ACPD}
on account of $\dashint_{Y}\vartheta_k^{\alpha\gamma}(y)dy = 0$ for $k=0,1,\ldots,d$. Furthermore
it is not hard to see that $\vartheta_k^{\alpha\gamma}$ is periodic and belongs to $H^2_{loc}(\mathbb{R}^d)$.}
\end{remark}

\begin{lemma}[Cacciopolli's inequality]
Suppose that $A$ satisfies $\eqref{a:1}$ and $\eqref{a:2}$,
and $f$ is a $Y$-periodic function in $L^2(Y;\mathbb{R}^{md})$.
Let $\chi$ be a weak solution of $\emph{div}(A\nabla\chi + f) = 0$ in $\mathbb{R}^d$.
Then, for any $q\in\mathbb{R}^m$ and $B\subset 2B$ with $0<r\leq 1$, we have
\begin{equation}\label{pri:2.13}
\Big(\dashint_{B} |\nabla\chi|^2 dy\Big)^{1/2}
\leq \frac{C}{r}\Big(\dashint_{2B}|\chi-q|^2 dy\Big)^{1/2} + C\Big(\dashint_{2B}|f|^2 dy\Big)^{1/2},
\end{equation}
where $C$ depends only on $\mu,m,d$.
\end{lemma}
\begin{pf}
The proof is standard, and we provide a proof for the sake of completeness. Let $v =(\chi-q)\phi^2$ be a test function with any
$q\in\mathbb{R}^{m}$, where
$\phi = 1$ in $B$, $\phi=0$ outside $2B$ and $|\nabla\phi|\leq C/r$. Then we have
\begin{equation}
\int_{\mathbb{R}^d}\phi^2A\nabla\chi\nabla\chi dy
+ \int_{\mathbb{R}^d}\phi A\nabla\chi\nabla\phi(\chi-q) dy
=\int_{\mathbb{R}^d} \phi^2f\cdot\nabla\chi dy
+ 2\int_{\mathbb{R}^d}\phi f\cdot\nabla\phi(\chi-q) dy.
\end{equation}
By using $\eqref{a:1}$ and H\"older's inequality coupled with Young's inequality, we derive
\begin{equation*}
\frac{\mu}{4}\int_{\mathbb{R}^d}\phi^2|\nabla\chi|^2 dy
\leq C\bigg\{\int_{\mathbb{R}^d}|\nabla\phi|^2|\chi-q|^2 dy +
\int_{\mathbb{R}^d} \phi^2|f|^2 dy
\bigg\},
\end{equation*}
and then have
\begin{equation*}
\int_B |\nabla\chi|^2 dy
\leq C\bigg\{\frac{1}{r^2}\int_{2B}|\chi-q|^2dy + \int_{2B}|f|^2 dy\bigg\}.
\end{equation*}
This implies the estimate $\eqref{pri:2.13}$, and we complete the proof.
\qed
\end{pf}

\begin{lemma}\label{lemma:6.4}
There exist $E_{jik}^{\alpha\gamma}\in H^1_{per}(Y)$ with $k = 0,1,\ldots,d$, such that
\begin{equation}\label{pri:6.6}
 b_{ik}^{\alpha\gamma} = \frac{\partial}{\partial y_j}\big\{E_{jik}^{\alpha\gamma}\big\}
 \qquad\text{and}\qquad
 E_{jik}^{\alpha\gamma} = - E_{ijk}^{\alpha\gamma},
\end{equation}
where $1\leq i,j\leq d$ and $1\leq\alpha,\gamma\leq m$. Moreover if $\chi_{k}$ is H\"older continuous,
then $E_{jik}^{\alpha\gamma}\in L^\infty(Y)$.
\end{lemma}
\begin{pf}
We provide a proof here for the sake of completeness. Based on the observation that
\begin{equation}\label{f:6.2}
 \dashint_{Y} b_{ik}^{\alpha\gamma}(y) dy = 0
 \qquad \text{and}\qquad
 \frac{\partial}{\partial y_i}\big\{b_{ik}^{\alpha\gamma}\big\} = 0,
\end{equation}
(The first equality of $\eqref{f:6.2}$ follows from $\eqref{f:6.1}$, and second one follows from $\eqref{pde:6.1}$ and $\eqref{pde:6.2}$.)
we can construct auxiliary periodic function $\theta_{ik}^{\alpha\gamma}\in H^2_{loc}(\mathbb{R}^d)$ satisfying
\begin{equation}\label{pde:6.4}
 \Delta \theta_{ik}^{\alpha\gamma} = b_{ik}^{\alpha\gamma} \quad \text{in}~~\mathbb{R}^d,
 \qquad \text{and}\qquad
 \dashint_{Y} \theta_{ik}^{\alpha\gamma}(y) dy = 0.
\end{equation}
Note that the equality $\eqref{f:6.2}$ guarantees the existence of $\theta_{ik}^{\alpha\gamma}$ in $H_{per}^{1}(Y)$
(see\cite[theorem 4.28]{ACPD}).
Let $E_{jik}^{\alpha\gamma} = \frac{\partial}{\partial y_j}\big\{\theta_{ik}^{\alpha\gamma}\big\}
- \frac{\partial}{\partial y_i}\big\{\theta_{jk}^{\alpha\gamma}\big\}$, obviously $E_{jik}^{\alpha\gamma} = -E_{ijk}^{\alpha\gamma}$.
Moreover we have
\begin{equation*}
 \frac{\partial}{\partial y_j}\big\{E_{jik}^{\alpha\gamma}\big\} = \Delta\theta_{ik}^{\alpha\gamma}
 - \frac{\partial^2}{\partial y_j\partial y_i}\big\{\theta_{jk}^{\alpha\gamma}\big\} =  b_{ik}^{\alpha\gamma},
\end{equation*}
by noting the fact that
$\frac{\partial}{\partial y_i}\big\{b_{ik}^{\alpha\gamma}\big\} = 0$ in $\mathbb{R}^d$ implies
$\Delta\frac{\partial}{\partial y_i}\big\{\theta_{ik}^{\alpha\gamma}\big\}=0$. From Liouvill's theorem, it follows that
$\frac{\partial}{\partial y_i}\big\{\theta_{ik}^{\alpha\gamma}\big\} = C$, and therefore we have
$\frac{\partial^2}{\partial y_j\partial y_i}\big\{\theta_{jk}^{\alpha\gamma}\big\} = 0$.

We now turn to prove $E_{jik}^{\alpha\gamma}\in L^\infty(Y)$ based on $\chi_{k}\in C^{0,\sigma}(\mathbb{R}^d)$.
Due to Cacciopolli's inequality $\eqref{pri:2.13}$,
\begin{equation*}
 \int_{B(x,r)}|\nabla\chi_k|^2 dy \leq \frac{C}{r^2} \int_{B(x,2r)}|\chi_k - \chi_k(x)|^2dy + Cr^d\leq Cr^{2\sigma+d-2}
\end{equation*}
for any $r\in(0,1)$ and $x\in Y$, and this leads to
\begin{equation}\label{f:6.3}
\dashint_{B(x,r)} |b_{ik}^{\alpha\gamma}|^2 dy \leq C\big\{1+r^{2\sigma-2}\big\}\leq Cr^{2\sigma-2}.
\end{equation}
Without loss of generality, we may assume $Y$ is centered at 0. Let $\varphi\in C_0^\infty(2Y)$ be a cut-off function satisfying $\varphi = 1$ in
$B(0,2/3)$, $\varphi = 0$ outside $B(0,4/5)$, and $|\nabla\varphi| \leq C$. We thus localize the equation $\eqref{pde:6.4}$ as
$\Delta\big(\varphi\theta_{ik}^{\alpha\gamma}\big) = \varphi b_{ik}^{\alpha\gamma} + 2\nabla\varphi\nabla\theta_{ik}^{\alpha\gamma}
 + \Delta\varphi\theta_{ik}^{\alpha\gamma}$ in $\mathbb{R}^d$. In view of the Newton potential, we acquire
\begin{equation*}
\begin{aligned}
\|\nabla\theta_{ik}^{\alpha\gamma}\|_{L^\infty(Y)}
&\leq C\Big\{1+\int_{2Y}\frac{|b_{ik}^{\alpha\gamma}|dy}{|x-y|^{d-1}}\Big\}
\leq C\Big\{1+\sum_{k=0}^{\infty}\int_{2^{-k}\leq|x-y|\leq 2^{-k+1}}\frac{|b_{ik}^{\alpha\gamma}|dy}{|x-y|^{d-1}}\Big\}\\
&\leq C\Big\{1+\sum_{k=0}^\infty2^{k(d-1)}\cdot2^{(1-k)d}\Big(\dashint_{|x-y|\leq 2^{-k+1}}|b_{ik}^{\alpha\gamma}|^2dy\Big)^{1/2}\Big\}\\
&\leq C\{1+\sum_{k=0}^\infty 2^{-\sigma k}\} < \infty,
\end{aligned}
\end{equation*}
where the estimate $|\nabla G(x,y)|\leq C|x-y|^{1-d}$ for Newton potential $G$ is employed in the first inequality,
and then we use H\"older's inequality in the third inequality and the estimate $\eqref{f:6.3}$ in last one.
This completes the proof.
\qed

\end{pf}

\begin{definition}\label{def:6.4}
\emph{Fix $\zeta\in C_0^\infty(B(0,1/2))$ such that $\zeta\geq 0$ and $\int_{\mathbb{R}^d}\zeta = 1$. Define
\begin{equation}\label{f:6.5}
S_\varepsilon(f)(x) = f*\zeta_\varepsilon(x) = \int_{\mathbb{R}^d} f(x-y)\zeta_\varepsilon(y)dy,
\end{equation}
where $\zeta_\varepsilon(x) = \varepsilon^{-d}\zeta(x/\varepsilon)$.}
\end{definition}

\begin{lemma}\label{lemma:6.1}
Let $f\in L^p(\mathbb{R}^d)$ for some $1\leq p<\infty$. Then for any $h\in L^p_{per}(\mathbb{R}^d)$,
\begin{equation}\label{pri:6.1}
\|h(\cdot/\varepsilon)S_\varepsilon(f)\|_{L^p(\mathbb{R}^d)}\leq C \|h\|_{L^p(Y)}\|f\|_{L^p(\mathbb{R}^d)},
\end{equation}
where $C$ depends only on $d$.
\end{lemma}

\begin{lemma}\label{lemma:6.2}
Let $f\in W^{1,p}(\mathbb{R}^d)$ for some $1<p<\infty$. Then we have
\begin{equation}\label{pri:6.2}
\|S_\varepsilon(f) - f\|_{L^p(\mathbb{R}^d)} \leq C\varepsilon\|\nabla f\|_{L^p(\mathbb{R}^d)},
\end{equation}
and furthermore obtain
\begin{equation}\label{pri:6.3}
\|S_\varepsilon(f)\|_{L^2(\mathbb{R}^d)}\leq C\varepsilon^{-1/2}\|f\|_{L^q(\mathbb{R}^d)}
\qquad \text{and}\qquad
\|S_\varepsilon(f)-f\|_{L^2(\mathbb{R}^d)}\leq C\varepsilon^{1/2}\|\nabla f\|_{L^q(\mathbb{R}^d)},
\end{equation}
where $q = \frac{2d}{d+1}$, and $C$ depends only on $d$.
\end{lemma}

\begin{remark}\label{re:4.3}
\emph{The proof of Lemma $\ref{lemma:6.1}$ and $\ref{lemma:6.2}$ can be found in \cite[pp.8]{SZW12}. If we define $\bar{S}_\varepsilon$ to be
\begin{equation}\label{def:6.3}
\bar{S}_\varepsilon(f)(x) = \dashint_{Y}f(x-\varepsilon z)dz
\end{equation}
for any $f\in L^2(\mathbb{R}^d)$
(which is referred to as the Steklov smoothing
operator originally applied to homogenization problem  by V.V. Zhikov in \cite{ZVVPSE} and further developed in \cite{TS2,TS}),
then the estimates $\eqref{pri:6.1}$ and $\eqref{pri:6.2}$ are still true for this kind of operators (see \cite[pp.3459]{TS}).
By definition, it is clear to see that the operator $S_\varepsilon$ defined in \cite{SZW12}
plays the same role as the Steklov smoothing operator $\bar{S}_\varepsilon$ in \cite{TS},
but seems more refined in view of the extra property $\eqref{pri:6.3}$ it satisfied. }
\end{remark}

\begin{remark}
\emph{For simplicity of presentation, if $f$ is a periodic function, we will denote $f(x/\varepsilon)$ by $f_{\varepsilon}(x)$.
For example, we usually write
$A_\varepsilon(x) = A(x/\varepsilon)$ and $\chi_{k,\varepsilon}(x) = \chi_k(x/\varepsilon)$,
and their components follow the same simplified way as well.}
\end{remark}

\begin{lemma}\label{lemma:6.3}
Let $\Omega$ be a bounded $C^1$ domain, and $\bar{S}_\varepsilon$ be given in $\eqref{def:6.3}$.
Suppose that $u\in H^{1}(\mathbb{R}^d;\mathbb{R}^m)$, and $f\in L_{per}^2(Y)$. Then we have
\begin{equation}\label{pri:4.4}
 \|f_\varepsilon\bar{S}_\varepsilon(u)\|_{L^2(\Omega\setminus\Sigma_\varepsilon)}
 \leq C\varepsilon^{\frac{1}{2}}\|f\|_{L^2(Y)}\|u\|_{H^{1}(\mathbb{R}^d)},
\end{equation}
and if $f\in H^1_{per}(Y)$, then
\begin{equation}\label{pri:4.9}
\|f_\varepsilon \bar{S}_\varepsilon(u)\|_{L^2(\partial\Omega)}\leq C\varepsilon^{-\frac{1}{2}}\|f\|_{H^1(Y)}\|u\|_{H^1(\mathbb{R}^d)}
\end{equation}
where $C$ depends only on $m,d$ and $\Omega$.
\end{lemma}

\begin{pf}
The proof of $\eqref{pri:4.4}$  will not be reproduced here.
We refer the reader to \cite[Lemma 3.3]{TS} and its references.
We now prove the estimate $\eqref{pri:4.9}$.
Let $\varrho \in C_0^1(\mathbb{R}^d;\mathbb{R}^d)$ be a
vector field such that $\big<\varrho,n\big> \geq c >0$ on $\partial\Omega$, and then
\begin{equation*}
\begin{aligned}
\|f_\varepsilon \bar{S}_\varepsilon(u)\|_{L^2(\partial\Omega)}^2
&\leq \frac{1}{c}\int_{\partial\Omega}<\varrho, n>|f_\varepsilon \bar{S}_\varepsilon(u)|^2 dS
= \frac{1}{c}\int_\Omega \text{div}\big(\varrho|f_\varepsilon \bar{S}_\varepsilon(u)|^2\big) dx\\
& = \frac{1}{c}\Big\{\int_\Omega \text{div}\varrho |f_\varepsilon \bar{S}_\varepsilon(u)|^2 dx
+ \int_\Omega <\varrho,\nabla(f_\varepsilon \bar{S}_\varepsilon(u))>f_\varepsilon \bar{S}_\varepsilon(u)dx\Big\} \\
&\leq C\big\{\|f_\varepsilon \bar{S}_\varepsilon(u)\|_{L^2(\mathbb{R}^d)}^2
+ \varepsilon^{-1}\|(\nabla f)_\varepsilon \bar{S}_\varepsilon(u)\|_{L^2(\mathbb{R}^d)}
\|f_\varepsilon \bar{S}_\varepsilon(u)\|_{L^2(\mathbb{R}^d)}
+ \|f_\varepsilon \bar{S}_\varepsilon(\nabla u)\|_{L^2(\mathbb{R}^d)}
\|f_\varepsilon \bar{S}_\varepsilon(u)\|_{L^2(\mathbb{R}^d)}
\big\},
\end{aligned}
\end{equation*}
where we use Cauchy's inequality in the last inequality and $(\nabla f)_\varepsilon(x) = \nabla f(x/\varepsilon)$.
Since the estimates $\eqref{pri:6.1}$ and $\eqref{pri:6.2}$ are still valid for $\bar{S}_\varepsilon$ (see Remark $\ref{re:4.3}$),
we obtain
\begin{equation*}
\begin{aligned}
\|f_\varepsilon \bar{S}_\varepsilon(u)\|_{L^2(\partial\Omega)}^2 &\leq C\big\{\|f\|_{L^2(Y)}^2 \|u\|_{L^2(\mathbb{R}^d)}^2
+ \varepsilon^{-1}\|\nabla f\|_{L^2(Y)}\|f\|_{L^2(Y)}\|u\|_{L^2(\mathbb{R}^d)}^2
+ \|f\|_{L^2(Y)}^2\|\nabla u\|_{L^2(\mathbb{R}^d)}\|u\|_{L^2(\mathbb{R}^d)}\big\}\\
&\leq C\varepsilon^{-1}\|f\|_{H^1(Y)}^2\|u\|_{H^1(\mathbb{R}^d)}^2.
\end{aligned}
\end{equation*}
This implies the estimate $\eqref{pri:4.9}$. The proof is completed.
\qed

\end{pf}

\begin{remark}\label{re:4.1}
\emph{Throughout the paper, let $B(P,r)$ denote the open ball centered at $P$ of radius r,
and the symbol $r_0$ only represents the diameter of $\Omega$. Since $\partial\Omega\in C^{1,\eta}$, there exists $R$ such that
for each point $P\in\partial\Omega$ there is a new coordinate system in $\mathbb{R}^d$ obtained
from the standard Euclidean coordinate system translation and rotation so that
$P=(0,0)$ and
\begin{equation*}
 B(P,R)\cap \Omega
 = B(P,R)\cap\big\{(x^\prime,x_d)\in\mathbb{R}^d:x^\prime\in\mathbb{R}^{d-1}
 ~\text{and}~ x_d>\phi(x^\prime)\big\},
\end{equation*}
where $\phi\in C^{1,\eta}(\mathbb{R}^{d-1})$ is a boundary function with $\phi(0) = 0$ and
$\|\nabla\phi\|_{C^{0,\eta}(\mathbb{R}^d)}\leq M_0$. Note that the pair of $(\eta,M_0)$ indicates the boundary character of $\Omega$.
To describe boundary estimates, we need more notation: let
\begin{equation*}
\begin{aligned}
D(P,r) & = B(P,r)\cap\Omega, \qquad\quad\Delta(P,r) & = B(P,r)\cap\Omega.
\end{aligned}
\end{equation*}
and $kD(P,r) = D(P,kr)$ and $k\Delta(P,r) = \Delta(P,kr)$ with $k>0$. We usually omit the center point and the radius of
$B(P,r),D(P,r),\Delta(P,r)$ without confusion.
In the paper, saying a constant $C$ depends on $\Omega$ means this constant involves both $(\eta,M_0)$ and $|\Omega|$,
where $|\Omega|$ denotes the volume of $\Omega$.}
\end{remark}

In the following, we introduce the Schauder estimates for ``classical''
Neumann problem. For this purpose, we set $L = -\text{div}(A\nabla)$ and
$\mathcal{L} = -\text{div}(A\nabla+V)+B\nabla +c +\lambda I$, where the coefficients $A,V,B,c$ do not dependent on $\varepsilon$.
Besides, $\partial u\setminus\partial\nu = n\cdot A\nabla u$ denotes the conormal derivative associated with $L$,
and $\mathcal{B}(u) = n\cdot\big(A\nabla u + Vu\big)$ represents the same thing for $\mathcal{L}$.
For simplicity of presentation,
we define
\begin{equation}\label{f:2.12}
 \mathcal{R}(r;D;\Delta;F;f;g)
 =  r\Big(\dashint_{D}|F|^pdx\Big)^{\frac{1}{p}}
+\|f\|_{L^\infty(D)}
+ r^\sigma[f]_{C^{0,\sigma}(D)}
 + \|g\|_{L^\infty(2\Delta)}
 + r^\sigma[g]_{C^{0,\sigma}(\Delta)},
\end{equation}
where $f\in C^{0,\sigma}(\Omega;\mathbb{R}^{md})$ with $\sigma\in(0,1)$, $F\in L^p(\Omega;\mathbb{R}^m)$ with $1\leq p<\infty$,
and $g\in C^{0,\sigma}(\Omega;\mathbb{R}^{m})$. Note that $D = D(P,r)$ and $\Delta = \Delta(P,r)$ for $P\in\partial\Omega$.

\begin{lemma}\label{lemma:2.4}
Let $\Omega$ be a bounded $C^{1,\tau}$ domain.
Suppose $A$ satisfies $\eqref{a:1}$ and $\eqref{a:4}$, and $f\in C^{0,\sigma}(\Omega;\mathbb{R}^{md})$ with $\sigma\in(0,\tau]$ and
$F\in L^p(\Omega;\mathbb{R}^d)$ with $p>d$. Let $u$ be a weak solution to
$L(u) = \emph{div}(f) + F$ in $D(Q,4r)$ and $\partial u/\partial\nu = g - n\cdot f$ on $\Delta(Q,4r)$, where $Q\in\partial\Omega$.
Then we have the following boundary estimates:
\begin{itemize}
\item[\emph{(i)}] the Schauder estimate
\begin{equation}\label{pri:2.5}
[\nabla u]_{C^{0,\sigma}(D)}
\leq Cr^{-\sigma}\bigg\{\Big(\dashint_{2D}|\nabla u|^2 dx\Big)^{\frac{1}{2}}
+\mathcal{R}(r;2D;2\Delta;F;f;g)\bigg\};
\end{equation}
\item[\emph{(ii)}] the Lipschitz estimate
\begin{equation}\label{pri:2.6}
\|\nabla u\|_{L^{\infty}(D)}
\leq C\bigg\{\Big(\dashint_{2D}|\nabla u|^2 dx\Big)^{\frac{1}{2}}
+\mathcal{R}(r;2D;2\Delta;F;f;g)\bigg\},
\end{equation}
\end{itemize}
where $\mathcal{R}(r;2D;2\Delta;F;f;g)$ is defined in $\eqref{f:2.12}$,
and $C$ depends on $\mu,\tau,\kappa,m,d,\sigma,p$
and the character of $\Omega$.
\end{lemma}

\begin{pf}
The estimate $\eqref{pri:2.5}$ can be found in \cite[Theorem 5.52]{GML}, and we thus omit the proof.
For (ii), we can straightforward derive $\eqref{pri:2.6}$ from $\eqref{pri:2.5}$ as follows. For any $x,y\in D$, we have
\begin{equation*}
 |\nabla u(x)| \leq |\nabla u(x) - \nabla u(y)| + |\nabla u(y)|
 \leq r^\sigma[\nabla u]_{C^{0,\sigma}(D)} + |\nabla u(y)|.
\end{equation*}
Integrating both sides with respect to $y$ on $D$ and divided by $|D|$, we arrive at
\begin{equation*}
 |\nabla u(x)| \leq r^\sigma[\nabla u]_{C^{0,\sigma}(D)}|D|^{-1} + \Big(\dashint_{D}|\nabla u| dy\Big)^{1/2}
 \leq C\bigg\{ \Big(\dashint_{2D}|\nabla u|^2 dy\Big)^{1/2} + \mathcal{R}(r;2D;2\Delta;F;f;g)  \bigg\},
\end{equation*}
where we use H\"older's inequality in the last inequality. This implies the desired estimate.
\qed
\end{pf}

\begin{remark}
\emph{Let $u$ be the weak solution to $L(u) = \text{div}(f)+F$ in $\Omega$
and $\partial u/\partial\nu = g-n\cdot f$ on $\partial\Omega$, and the corresponding assumptions are given in Lemma $\ref{lemma:2.4}$.
Based on the estimates $\eqref{pri:2.5}$, $\eqref{pri:2.6}$ and \cite[Theorem 5.14]{MGLM} (the corresponding interior estimate),
it is not hard to derive the global estimates:
\begin{equation}\label{pri:2.9}
\max\big\{\|\nabla u\|_{L^\infty(\Omega)},[\nabla u]_{C^{0,\sigma}(\Omega)}\big\}
\leq C\big\{\|f\|_{C^{0,\sigma}(\Omega)} + \|F\|_{L^p(\Omega)} + \|g\|_{C^{0,\sigma}(\partial\Omega)}\big\},
\end{equation}
where $C$ depends on $\mu,\tau,m,d,p,\sigma$ and $\Omega$. In fact, under the assumptions of Lemma $\eqref{lemma:2.4}$.
It is not hard to derive $W^{1,p}$ estimates for $2\leq p<\infty$ (by using the methods developed in \cite{GZ1}) that
\begin{equation}\label{pri:2.10}
 \|\nabla u\|_{L^p(\Omega)}
 \leq C\big\{\|f\|_{L^p(\Omega)}+\|F\|_{L^{\frac{pd}{p+d}}(\Omega)}+\|g\|_{B^{-1/p,p}(\partial\Omega)}\big\}.
\end{equation}
We finally mention that the $W^{1,p}$ estimates actually hold for $1<p<2$ by the duality argument.}
\end{remark}

\begin{lemma}\label{lemma:2.5}
Let $\Omega$ be a bounded $C^{1,\tau}$ domain.
Suppose that the coefficients of $\mathcal{L}$ satisfy $\eqref{a:1}$, $\eqref{a:3}$ and $\eqref{a:4}$,
and $f\in C^{0,\sigma}(\Omega;\mathbb{R}^{md})$ with $\sigma\in(0,\tau]$ and
$F\in L^p(\Omega;\mathbb{R}^d)$ with $p>d$. Let $u$ be a weak solution to
$\mathcal{L}(u) = \emph{div}(f) + F$ in $D(Q,4r)$ and $\mathcal{B}(u) = g - n\cdot f$ on $\Delta(Q,4r)$, where $Q\in\partial\Omega$.
Then we have the following boundary estimates:
\begin{itemize}
\item[\emph{(i)}] the Lipschitz estimate
\begin{equation}\label{pri:2.8}
\|\nabla u\|_{L^{\infty}(D)}
\leq C\bigg\{\frac{1}{r}\Big(\dashint_{2D}|u|^2 dx\Big)^{\frac{1}{2}}
+\mathcal{R}(r;2D;2\Delta;F;f;g)\bigg\},
\end{equation}
\item[\emph{(ii)}] the Schauder estimate
\begin{equation}\label{pri:2.7}
[\nabla u]_{C^{0,\sigma}(D)}
\leq Cr^{-\sigma}\bigg\{\frac{1}{r}\Big(\dashint_{2D}|u|^2 dx\Big)^{\frac{1}{2}}
+\mathcal{R}(r;2D;2\Delta;F;f;g)\bigg\};
\end{equation}
\end{itemize}
where $ \mathcal{R}(r;2D;2\Delta;F;f;g)$ is defined in $\eqref{f:2.12}$,
and $C$ depends on $\mu,\tau,\kappa,m,d,\sigma,p$ and the character of $\Omega$.
\end{lemma}

\begin{pf}
There are several methods to derive the above estimates. We plan to
show the corresponding global estimates and then employ localization argument (introduced in \cite[Remark 2.11]{QXS})
to establish the desired estimates. We rewrite $\mathcal{L}(u) = \text{div}(f)+F$ as
$L(u) = \text{div}(f+Vu)+F-B\nabla u -(c+\lambda I)u$ in $\Omega$, and $\mathcal{B}(u) = g-n\cdot f$
as $\partial u/\partial\nu = g-n\cdot (f+Vu)$ on $\partial\Omega$. Then it follows from $\eqref{pri:2.9}$ that
\begin{equation*}
\begin{aligned}
\|\nabla u\|_{L^\infty(\Omega)}
&\leq C\big\{\|f\|_{C^{0,\sigma}(\Omega)} + \|F\|_{L^p(\Omega)} + \|g\|_{C^{0,\sigma}(\partial\Omega)}
+ \|u\|_{C^{0,\sigma}(\Omega)} + \|u\|_{W^{1,p}(\Omega)}\big\}\\
& \leq \big\{\|f\|_{C^{0,\sigma}(\Omega)} + \|F\|_{L^p(\Omega)} + \|g\|_{C^{0,\sigma}(\partial\Omega)}
+2\|\nabla u\|_{L^\infty(\Omega)}^{\sigma}\|u\|_{L^\infty(\Omega)}^{1-\sigma} + \|u\|_{W^{1,p}(\Omega)} \big\} \\
& \leq C\big\{\|f\|_{C^{0,\sigma}(\Omega)} + \|F\|_{L^p(\Omega)} + \|g\|_{C^{0,\sigma}(\partial\Omega)}
+ \|u\|_{W^{1,p}(\Omega)} + \|u\|_{L^\infty(\Omega)}\big\}
+ \frac{1}{2}\|\nabla u\|_{L^\infty(\Omega)},
\end{aligned}
\end{equation*}
where we use the interpolation inequality in the second inequality, and Young's inequality in the last one. This implies
\begin{equation}\label{pri:2.11}
\begin{aligned}
\|\nabla u\|_{L^\infty(\Omega)}
&\leq C\big\{\|f\|_{C^{0,\sigma}(\Omega)} + \|F\|_{L^p(\Omega)} + \|g\|_{C^{0,\sigma}(\partial\Omega)}
+ \|u\|_{W^{1,p}(\Omega)}\big\}\\
&\leq C\big\{\|f\|_{C^{0,\sigma}(\Omega)} + \|F\|_{L^p(\Omega)} + \|g\|_{C^{0,\sigma}(\partial\Omega)}\big\}
\end{aligned}
\end{equation}
where we use the Sobolev embedding theorem $\|u\|_{L^\infty(\Omega)}\leq C\|u\|_{W^{1,p}(\Omega)}$ (for $p>d$) in the first inequality,
and $W^{1,p}$ estimates (which follows from $\eqref{pri:2.10}$ by using the same idea as in the proof of Theorem 1.1,
and we thus omit details here) in the last one. Moreover, combining $\eqref{pri:2.9}$ and $\eqref{pri:2.11}$, we have
\begin{equation}\label{pri:2.12}
\begin{aligned}
\big[\nabla u\big]_{C^{0,\sigma}(\Omega)}
&\leq C\big\{\|f\|_{C^{0,\sigma}(\Omega)} + \|F\|_{L^p(\Omega)} + \|g\|_{C^{0,\sigma}(\partial\Omega)}
+ \|u\|_{W^{1,\infty}(\Omega)}\big\} \\
&\leq C\big\{\|f\|_{C^{0,\sigma}(\Omega)} + \|F\|_{L^p(\Omega)} + \|g\|_{C^{0,\sigma}(\partial\Omega)} \big\}.
\end{aligned}
\end{equation}
Finally, by the localization argument (shown in \cite[Remark 2.11]{QXS}),
we can derive $\eqref{pri:2.8}$, $\eqref{pri:2.7}$ from $\eqref{pri:2.11}$ and $\eqref{pri:2.12}$,  respectively.
The proof is complete.
\qed
\end{pf}


\begin{remark}\label{re:2.2}
\emph{The following notations will be used frequently:
\begin{itemize}
\item distance function $\delta(x)=\text{dist}(x,\partial\Omega)$, where $x\in\Omega$;
\item average of function $\bar{f}_E = \dashint_{E} f(x) dx = \frac{1}{|E|}\int_E f(x) dx $,
where $E$ is a subset of $\mathbb{R}^d$, and the subscript of $\bar{f}_E$ is usually omitted;
\item boundary layer $\Omega\setminus\Sigma_r$, where $\Sigma_{r}=\{x\in\Omega:\text{dist}(x,\partial\Omega)>r\}$ with $r>0$;
\item cut-off function $\psi_r$ (associated with $\Sigma_r$), satisfying
\begin{equation}\label{f:4.20}
\psi_{r} =1 \quad\text{in}\quad \Sigma_{2r},
\qquad\psi_{r} = 0 \quad\text{outside}\quad \Sigma_{r}, \qquad\text{and}\quad
|\nabla\psi_{r}|\leq C/r;
\end{equation}
\item  level set $S_r = \big\{x\in\Omega:\text{dist}(x,\partial\Omega) = r\big\}$;
\item internal diameter $r_{00} = \min\{\text{dist}(x,y):x,y\in\partial\Omega\}$, and layer constant $c_0 = r_{00}/10$.
\end{itemize}}
\end{remark}

\begin{remark}\label{re:2.1}
\emph{For $0\leq r< c_0$,
we may assume that there exist homeomorphisms $\Lambda_r: \partial\Omega\to \partial\Sigma_r = S_r$ such that $\Lambda_0(Q) = Q$,
$|\Lambda_r(Q) - \Lambda_t(P)| \sim |r-t| + |Q-P|$ and
$|\Lambda_r(Q) - \Lambda_t(Q)|\leq C\text{dist}(\Lambda_r(Q),S_t)$ for any
$r>s$ and $P,Q\in\partial\Omega$ (which are bi-Lipschitz maps, see \cite[pp.1014]{SZW2}). Especially, we may have
$\max_{r\in[0,c_0]}\{\|\nabla\Lambda_r\|_{L^\infty(\partial\Omega)},\|\nabla(\Lambda_r^{-1})\|_{L^\infty(\partial\Omega)}\} \leq C(\eta,M_0)$. For a function $h$,
we define the radial maximal function
$\mathcal{M}(h)$ on $\partial\Omega$ as
\begin{equation}\label{def:2}
 \mathcal{M}(h)(Q) = \sup\big\{|h(\Lambda_r(Q))|: 0\leq r\leq  c_0\big\} \quad\qquad \forall ~Q\in\partial\Omega.
\end{equation}
We mention that the radial maximal function will play an important role in the study of convergence rates for Lipschitz domains
(we refer the reader to\cite{SZW2} for the original thinking,
and we also refer the reader to \cite[Theorem 5.1]{CTM} for the existence of such bi-Lipschitz maps).}
\end{remark}

\begin{definition}
The non-tangential maximal function of $u$ is defined by
\begin{equation}\label{def:3}
(u)^*(Q) = \sup\big\{ |u(x)|:x\in \Gamma_{N_0}(Q)\big\} \qquad\quad \forall~ Q\in\partial\Omega,
\end{equation}
where $\Gamma_{N_0}(Q) = \{x\in\Omega:|x-Q|\leq N_0\delta(x)\}$ is the cone with vertex $Q$ and aperture $N_0$,
and $N_0>1$ is sufficiently large.
\end{definition}

\begin{remark}
\emph{Let $h\in L^p(\Omega)$ with $1\leq p<\infty$. For any $r\in(0,c_0)$ ($c_0$ and $\Lambda_r$ are given in Remark $\ref{re:2.2}$),
we can show the estimate of $\|h\|_{L^p(\Omega\setminus\Sigma_r)}$. By $\eqref{def:2}$, we note that
$h(\Lambda_r(x))\leq \mathcal{M}(h)(x)$ a.e. $x\in\partial\Omega$ for all $r\in (0,c_0)$.
Then
\begin{equation}\label{pri:6.5}
\begin{aligned}
 \int_{\Omega\setminus\Sigma_r} |h|^p dx
 &= \int_0^r\int_{S_t=\Lambda_{t}(\partial\Omega)} |h(y)|^p dS_t(y)dt \\
 & = \int_0^r\int_{\partial\Omega} |h(\Lambda_t(z))|^p |\nabla\Lambda_{t}|dS(z)dt \leq Cr\int_{\partial\Omega} |\mathcal{M}(h)|^p dS
 \leq  Cr\int_{\partial\Omega} |(h)^{*}|^p dS,
\end{aligned}
\end{equation}
where $C$ depends only on $p$ and the boundary character. We note that the first equality is based on so-called co-area formula $\eqref{f:2.11}$,
and we use the change of variable in the second one. Besides,
the first inequality follows from Remark $\ref{re:2.1}$. In the last one,
it is not hard to see $\mathcal{M}(h)(Q)\leq (h)^*(Q)$ by comparing the Definition $\eqref{def:3}$ with $\eqref{def:2}$.}

\emph{We now explain the co-area formula used here.
Let $Z(0;r)=\{x\in\Omega:0<\delta(x)\leq r\}$,
then $Z(0;r) = \Omega\setminus\Sigma_r$.
Here we point out $|\nabla\delta(x)| =1$ a.e. $x\in\Omega$ without the proof (see \cite[pp.142]{LCE1}).
In view of co-area formula (see \cite[Theorem 3.13]{LCE1}), we have
\begin{equation}\label{f:2.11}
 \int_{\Omega\setminus\Sigma_r} |h|^p dx = \int_{Z(0;r)} |h|^p dx
 = \int_{0}^{r}\int_{\{x\in\Omega:\delta(x)=t\}}\frac{|h|^p}{|\nabla\delta|}d\mathcal{H}^{d-1}dt
 =\int_{0}^{r}\int_{S_t}|h|^p dS_tdt,
\end{equation}
where $S_r=\{x\in\Omega:\delta(x) = t\}$, $d\mathcal{H}^{d-1}$ is the ($d-1$)-dimensional Hausdorff measure,
and $dS_t=d\mathcal{H}^{d-1}(S_t)$ denotes the surface measure of $S_t$.}
\end{remark}

\begin{lemma}\label{lemma:2.6}
Let $\Omega$ be a Lipschitz domain, and $\mathcal{M}$ associated with $c_0$ is defined in Remark $\ref{re:2.1}$.
Then for any $h\in H^1(\Omega)$, we have the following estimate
\begin{equation}\label{pri:2.14}
 \|\mathcal{M}(h)\|_{L^2(\partial\Omega)}
 \leq C\|h\|_{H^1(\Omega\setminus\Sigma_{c_0})},
\end{equation}
where $C$ depends only on $d,c_0$ and the character of $\Omega$.
\end{lemma}

\begin{pf}
The proof is based on the fundamental theorem of calculus and definition of the radial maximal function,
and our proof follows the idea from \cite[Proposition 8.4]{SZW2} (they actually proved a weighted inequality).
For any $0<s\leq t<c_0$ and $P\in\partial\Omega$, we first have
\begin{equation*}
h(\Lambda_{t}(P)) - h(\Lambda_{s}(P)) = \int_s^t \frac{d}{d t}\{h(\Lambda_r(P))\} dr
\end{equation*}
and then
\begin{equation}
|h(\Lambda_{t}(P))| \leq C \int_0^{c_0}|\nabla h(\Lambda_r(P))| dr + |h(\Lambda_{s}(P))|.
\end{equation}
Integrating both sides of the above inequality with respect to $s$, and then divided by $c_0$, we arrive at
\begin{equation*}
|h(\Lambda_{t}(P))| \leq C \int_0^{c_0} \big(|\nabla h(\Lambda_r(P))| + |h(\Lambda_{r}(P))|\big) dr
\leq C \Big(\int_0^{c_0}\big(|\nabla h(\Lambda_r(P))|^2 + |h(\Lambda_{r}(P))|^2\big) dr\Big)^{1/2}.
\end{equation*}
This implies
\begin{equation*}
|\mathcal{M}(h)(Q)|^2 \leq C \int_0^{c_0}\big(|\nabla h(\Lambda_r(P))|^2 + |h(\Lambda_{r}(P))|^2\big) dr.
\end{equation*}
Hence we obtain
\begin{equation*}
\int_{\partial\Omega} |\mathcal{M}(h)(Q)|^2 dS
\leq C\int_0^{c_0}\int_{\partial\Omega} \big(|\nabla h(\Lambda_r(P))|^2 + |h(\Lambda_{r}(P))|^2\big) dSdr
\leq C\int_{\Omega\setminus\Sigma_{c_0}} (|\nabla u|^2 + |u|^2) dx,
\end{equation*}
and this gives the desired estimate $\eqref{pri:2.14}$.
\qed

\end{pf}

Let $N_\varepsilon(x,z)=[N_\varepsilon^{\alpha\gamma}(x,z)]$ denote the Neumann matrix associated with $L_\varepsilon$
in $\Omega$ with pole at $z$, which solves the following Neumann boundary value problem:
\begin{equation*}
\left\{\begin{aligned}
 L_\varepsilon(N_\varepsilon(\cdot,z)) &= \mathbf{\delta}_z(x)I &\quad&\text{in}~~\Omega,\\
 \frac{\partial}{\partial\nu_\varepsilon}\big\{N_\varepsilon(\cdot,z)\big\}
 & = -|\partial\Omega|^{-1}I &\quad&\text{on}~\partial\Omega,
\end{aligned}\right.
\end{equation*}
where $\delta_z(x)$ is the Dirac delta function with pole at $z$.

\begin{remark}
\emph{If $A\in\Lambda(\mu,\tau,\kappa)$,
then we have the decay estimates of the Neumann matrix as follows:
\begin{equation}\label{pri:6.4}
\left\{\begin{aligned}
|N_\varepsilon(x,z)| &\leq C|x-z|^{2-d}, \\
|\nabla_xN_\varepsilon(x,z)| + |\nabla_zN_\varepsilon(x,z)| &\leq C|x-z|^{1-d},\\
|\nabla_x\nabla_zN_\varepsilon(x,z)|&\leq C|x-z|^{-d}
\end{aligned}\right.
\end{equation}
for any $x,z\in \Omega$ and $x\not=z$, where $C$ depends on $\mu,\tau,\kappa,m,d$ and $\Omega$.
We mention that the symmetry condition $A^* = A$ is not necessary any longer, due to the new method developed in \cite{SZ,SZW12}.
We also refer the reader to \cite{SZW4} for the proof of $\eqref{pri:6.4}$.}
\end{remark}

\section{$W^{1,p}$ Estimates}

\begin{lemma}[$H^1$ estimates]\label{lemma:2.1}
Let $\Omega$ be a Lipschitz domain.
Suppose that the coefficients of $\myl{L}{\varepsilon}$ satisfy $\eqref{a:1}$ and $\eqref{a:3}$.
Then for any $f\in L^{2}(\Omega;\mathbb{R}^{md})$, $F\in L^{\frac{2d}{d+2}}(\Omega;\mathbb{R}^m)$ and
$g\in B^{-1/2,2}(\partial\Omega;\mathbb{R}^m)$, there exists a unique weak solution
$u_\varepsilon\in H^1(\Omega;\mathbb{R}^m)$ to $\eqref{pde:1.1}$, whenever $\lambda\geq\lambda_0(\mu,\kappa,m,d)$, and $\lambda_0$ is sufficiently large.
Moreover, we have the uniform estimate
\begin{equation}\label{pri:2.4}
\|u_\varepsilon\|_{H^1(\Omega)} \leq C\big\{\|f\|_{L^{2}(\Omega)}+\|F\|_{L^{\frac{2d}{d+2}}(\Omega)}
+\|g\|_{B^{-1/2,2}(\partial\Omega)}\big\},
\end{equation}
where $C$ depends only on $\mu,m,d$ and $\Omega$.
\end{lemma}
\begin{pf}
 We write out the proof for the sake of completeness.
 First, we need to verify the boundedness and coercivity of $\mathrm{B}_\varepsilon[\cdot,\cdot]$. The easy one is the boundedness:
 \begin{equation*}
 \big|\mathrm{B}_\varepsilon[u_\varepsilon,\phi]\big| \leq C(\mu,\kappa,m,d)\|u_\varepsilon\|_{H^1(\Omega)}\|\phi\|_{H^1(\Omega)}
 \end{equation*}
 for any $u_\varepsilon,\phi\in H^1(\Omega;\mathbb{R}^m)$. Then set $\phi = u_\varepsilon$ in $\eqref{pde:1.6}$, and a routine computation
 gives rise to the coercivity:
 \begin{equation}\label{f:2.3}
  \mathrm{B}_\varepsilon[u_\varepsilon,u_\varepsilon] \geq \frac{\mu}{2}\int_\Omega |\nabla u_\varepsilon|^2 dx
  + \big[\lambda - C(\mu,\kappa,m,d)\big]\int_\Omega |u_\varepsilon|^2 dx
  \geq c \|u_\varepsilon\|_{H^1(\Omega)}^2,
 \end{equation}
 where $c = \min\{\mu/2,\lambda-C(\mu,\kappa,m,d)\}$. So we can choose a sufficiently large number $\lambda_0 =\lambda_0(\mu,\kappa,m,d)$
 such that $c = \mu/2$ whenever $\lambda>\lambda_0$.

 The next thing is to prove $\mathfrak{F}\in H^{-1}_0(\Omega;\mathbb{R}^m)=\big[H^1(\Omega;\mathbb{R}^m)\big]^*$,
 where $\mathfrak{F}$ denotes the right-hand side of $\eqref{pde:1.4}$ in the sense of
 \begin{equation*}
  <\mathfrak{F},\phi>
  = -\int_\Omega f^\alpha\cdot\nabla\phi^\alpha dx + \int_\Omega F^\alpha\phi^{\alpha} dx
 + <g,\phi>_{B^{-1/2,2}(\partial\Omega)\times B^{1/2,2}(\partial\Omega)}
 \end{equation*}
for any $\phi\in H^1(\Omega;\mathbb{R}^m)$. It is apparent to see that
 \begin{equation*}
 \begin{aligned}
 |<\mathfrak{F},\phi>|
 &\leq C\big\{\|f\|_{L^{2}(\Omega)} + \|F\|_{L^{\frac{2d}{d+2}}(\Omega)}
 +\|g\|_{B^{-1/2,2}(\partial\Omega)}\big\}\|\phi\|_{H^1(\Omega)},
\end{aligned}
\end{equation*}
and this leads to
\begin{equation}\label{f:2.4}
\|\mathfrak{F}\|_{H^{-1}_0(\Omega)}
\leq C\big\{\|f\|_{L^{2}(\Omega)} + \|F\|_{L^{\frac{2d}{d+2}}(\Omega)}
 +\|g\|_{B^{-1/2,2}(\partial\Omega)}\big\},
\end{equation}
where $C$ depends on $m,d,\Omega$.

Finally, due to the Lax-Milgram theorem, there exists a unique weak solution $u_\varepsilon\in H^1(\Omega;\mathbb{R}^m)$ to $\eqref{pde:1.1}$,
such that $ \mathrm{B}_\varepsilon[u_\varepsilon,\phi] = <\mathfrak{F},\phi>$ holds for all $\phi\in H^1(\Omega;\mathbb{R}^m)$.
and the estimate $\eqref{pri:2.4}$ follows from $\eqref{f:2.3}$ and $\eqref{f:2.4}$. The proof is done.
\qed
\end{pf}

\begin{remark}
\emph{The adjoint operator $\mathcal{L}_\varepsilon^*$ has the same results as that in Lemma $\ref{lemma:2.1}$
and the corresponding compatibility condition $\eqref{C:1}$.}
\end{remark}

\begin{lemma}\label{lemma:2.2}
Let $2\leq p<\infty$. Suppose $A\in \emph{VMO}(\mathbb{R}^d)$ satisfies $\eqref{a:1}$ and $\eqref{a:2}$.
Let $f\in L^{p}(\Omega;\mathbb{R}^{md})$, $F\in L^q(\Omega;\mathbb{R}^m)$
and $g\in B^{-1/p,p}(\partial\Omega;\mathbb{R}^m)$, where $q=\frac{pd}{p+d}$.
Then if $F$ and $g$ satisfy the
compatibility condition
$\int_\Omega F^\alpha dx + <g^\alpha,1> = 0$, the weak solution $u_\varepsilon\in W^{1,p}(\Omega;\mathbb{R}^m)$ to
$L_\varepsilon(u_\varepsilon) = \emph{div}(f) + F$ in $\Omega$ and $\partial u_\varepsilon/\partial\nu_\varepsilon = g - n\cdot f$
on $\partial\Omega$ satisfies the uniform estimate
\begin{equation}\label{pri:2.2}
\|\nabla u_\varepsilon\|_{L^p(\Omega)} \leq C\big\{\|f\|_{L^{p}(\Omega)}+\|F\|_{L^q(\Omega)} + \|g\|_{B^{-1/p,p}(\partial\Omega)}\big\},
\end{equation}
where $C$ depends only on $\mu,\omega,m,d,p$ and $\Omega$.
\end{lemma}

\begin{pf}
In fact, the estimate $\eqref{pri:2.2}$ holds for $1<p<\infty$. The proof can be found in \cite[Theorem 1.1]{SZW4}
and thus is not reproduced here.
\qed
\end{pf}

\begin{remark}
\emph{Before approaching the proof of Theorem $\ref{thm:1.1}$, we first introduce the following interpolation inequalities.
For any $u\in W^{1,p}(\Omega;\mathbb{R}^m)$ with $2< p<\infty$,
and for any $\delta>0$, there exists a constant $C_\delta$ depending on $\delta,p,m,d$ and $\Omega$, such that
\begin{equation}\label{f:2.8}
 \|u\|_{L^p(\Omega)} \leq \delta\|\nabla u\|_{L^p(\Omega)} + C_\delta\|u\|_{L^2(\Omega)}
 \quad \text{and} \quad\|\nabla u\|_{L^q(\Omega)} \leq \delta\|\nabla u\|_{L^p(\Omega)} + C_\delta\|\nabla u\|_{L^2(\Omega)},
\end{equation}
where $q=\frac{pd}{d+p}$.
Outline of proof: (i) The first inequality of $\eqref{f:2.8}$ follows from a contradiction argument. (ii)
The second one comes from H\"older's inequality and Young's inequality with $\delta$ in the case of $p>\frac{2d}{d-2}$, while it is
trivial for $2< p\leq \frac{2d}{d-2}$ provided $C_\delta \geq 1$.
The details are left to the reader (or see \cite{RAA}).}
\end{remark}

\begin{flushleft}
\textbf{Proof of Theorem \ref{thm:1.1}}\textbf{.}\quad
The case of $p=2$ has been investigated in Lemma $\ref{lemma:2.1}$.
We first consider the case of $2< p<\infty$, and rewrite the $\eqref{pde:1.1}$ as
\begin{equation*}
\left\{\begin{aligned}
 L_\varepsilon(u_\varepsilon) &= \text{div}(\tilde{f}) + \tilde{F} &\quad& \text{in}~\Omega \\
 \frac{\partial u_\varepsilon}{\partial \nu_\varepsilon} &= g-n\cdot\tilde{f} &\quad& \text{on}~\partial\Omega,
\end{aligned}\right.
\end{equation*}
where
\begin{equation*}
\begin{aligned}
&\tilde{f} = f^\alpha+V_\varepsilon u_\varepsilon, \\
& \tilde{F} = F - B_\varepsilon\nabla u_\varepsilon -(c_\varepsilon+\lambda I)u_\varepsilon.
\end{aligned}
\end{equation*}\end{flushleft}
On account of $\eqref{C:1}$, we have $ \int_\Omega\tilde{F}^\alpha dx + <g^\alpha,1> = 0$. It follows from
$\eqref{pri:2.2}$ and $\eqref{f:2.8}$ that
\begin{equation*}
\begin{aligned}
\|\nabla u_\varepsilon\|_{L^p(\Omega)}
&\leq C\big\{\|\tilde{f}\|_{L^{p}(\Omega)} + \|\tilde{F}\|_{L^q(\Omega)} + \|g\|_{B^{-1/p,p}(\partial\Omega)}\big\}\\
&\leq C\big\{\|f\|_{L^{p}(\Omega)} +\|F\|_{L^q(\Omega)} + \|g\|_{B^{-1/p,p}(\partial\Omega)}
 +\|u_\varepsilon\|_{L^p(\Omega)} + \|u_\varepsilon\|_{W^{1,q}(\Omega)} \big\} \\
&\leq C\big\{\|f\|_{L^{p}(\Omega)} +\|F\|_{L^q(\Omega)} + \|g\|_{B^{-1/p,p}(\partial\Omega)}
+ \|u_\varepsilon\|_{H^1(\Omega)}\big\} + C\delta\|\nabla u_\varepsilon\|_{L^p(\Omega)},
\end{aligned}
\end{equation*}
and this leads to
\begin{equation}
\begin{aligned}
 \|u_\varepsilon\|_{W^{1,p}(\Omega)}
 &\leq 2\|\nabla u_\varepsilon\|_{L^p(\Omega)} + C_1\|u_\varepsilon\|_{L^2(\Omega)} \\
 &\leq C\big\{\|f\|_{L^{p}(\Omega)} +\|F\|_{L^q(\Omega)} + \|g\|_{B^{-1/p,p}(\partial\Omega)}\big\},
\end{aligned}
\end{equation}
where we choose $\delta$ such that $C\delta\leq1/2$,
and note that $B^{-1/p,p}(\partial\Omega;\mathbb{R}^m)\subset B^{-1/2,2}(\partial\Omega;\mathbb{R}^m)$ holds for
$2\leq p<\infty$. We also use the estimate $\eqref{pri:2.4}$ in the computation above.

In the case of $1<p<2$, we employ the duality argument.
For any $h\in C^{1}_0(\Omega;\mathbb{R}^{md})$,
let $v_\varepsilon\in W^{1,p^\prime}(\Omega;\mathbb{R}^m)\cap H^1(\Omega;\mathbb{R}^m)$ ($p^\prime=\frac{p}{p-1}$) be the weak solution of
$\mathcal{L}_\varepsilon^*(v_\varepsilon) = \text{div}(h)$ in $\Omega$ and $\mathcal{B}_\varepsilon^*(v_\varepsilon) = 0$
on $\partial\Omega$. In view of $\eqref{G:1}$, we have
\begin{equation*}
\begin{aligned}
-\int_\Omega \nabla u_\varepsilon\cdot h &= <u_\varepsilon,\mathcal{L}_\varepsilon^*(v_\varepsilon)>
 = <\mathcal{L}_\varepsilon(u_\varepsilon),v_\varepsilon> + <\mathcal{B}_\varepsilon(u_\varepsilon),v_\varepsilon>\\
 &= -\int_\Omega f\cdot\nabla v_\varepsilon dx + \int_\Omega F v_\varepsilon dx + <g,v_\varepsilon> =: <\mathfrak{F},v_\varepsilon>.
\end{aligned}
\end{equation*}
This coupled with
\begin{equation*}
\begin{aligned}
|<\mathfrak{F},v_\varepsilon>| &\leq \|f\|_{L^{p}(\Omega)}\|\nabla v_\varepsilon\|_{L^{p^\prime}(\Omega)}
+ \|F\|_{L^p(\Omega)}\|v_\varepsilon\|_{L^{p^\prime}(\Omega)} + C\|g\|_{B^{-1/p,p}(\partial\Omega)}\|v_\varepsilon\|_{W^{1,p^\prime}(\Omega)}\\
& \leq C\big\{\|f\|_{L^{p}(\Omega)}+\|F\|_{L^p(\Omega)}+\|g\|_{B^{-1/p,p}(\partial\Omega)}\big\}\|h\|_{L^{p^\prime}(\Omega)}
\end{aligned}
\end{equation*}
indicates
\begin{equation}\label{f:2.5}
 \|\nabla u_\varepsilon\|_{L^{p}(\Omega)} \leq C\big\{\|f\|_{L^{p}(\Omega)}+\|F\|_{L^p(\Omega)}+\|g\|_{B^{-1/p,p}(\partial\Omega)}\big\}.
\end{equation}

We proceed with the same idea to derive
\begin{equation}\label{f:2.9}
 \|u_\varepsilon\|_{L^{p}(\Omega)} \leq C\big\{\|f\|_{L^{p}(\Omega)}+\|F\|_{L^p(\Omega)}+\|g\|_{B^{-1/p,p}(\partial\Omega)}\big\},
\end{equation}
and this together with $\eqref{f:2.5}$ gives $\eqref{pri:1.1}$ in the case of $1<p<2$.

To see $\eqref{f:2.9}$, for any $H\in L^{p^\prime}(\Omega;\mathbb{R}^m)$,
let $w_\varepsilon\in W^{1,p^\prime}(\Omega;\mathbb{R}^m)\cap H^1(\Omega;\mathbb{R}^m)$ be the solution of
$\mathcal{L}_\varepsilon^*(w_\varepsilon) = H$ in $\Omega$ and $\mathcal{B}_\varepsilon^*(w_\varepsilon)=0$. Then we have
$\int_\Omega u_\varepsilon H dx = <\mathfrak{F},w_\varepsilon>$ and
\begin{equation*}
|<\mathfrak{F},w_\varepsilon>|\leq C\{\|f\|_{L^{p}(\Omega)}+\|F\|_{L^p(\Omega)}+\|g\|_{B^{-1/p,p}(\partial\Omega)}\}\|H\|_{L^{p^\prime}(\Omega)}.
\end{equation*}
The estimate $\eqref{f:2.9}$ straightforward follows, and we complete the proof.
\qed

\begin{cor}\label{cor:2.1}
 Suppose that the coefficients of $\mathcal{L}_\varepsilon$ satisfy the same conditions as in Theorem $\ref{thm:1.1}$.
 Let $p>d$ and $\sigma = 1-d/p$. Assume $u_\varepsilon\in W^{1,p}(\Omega;\mathbb{R}^m)\cap H^1(\Omega;\mathbb{R}^m)$ is the weak solution to
 $\eqref{pde:1.1}$, where $f\in L^p(\Omega;\mathbb{R}^{md})$, $F\in L^q(\Omega;\mathbb{R}^m)$ with $q=\frac{pd}{d+p}$,
 and $g\in L^\infty(\partial\Omega;\mathbb{R}^m)$. Then we have the uniform H\"older estimate
 \begin{equation}\label{pri:2.3}
  \|u_\varepsilon\|_{C^{0,\sigma}(\Omega)} \leq C\big\{\|f\|_{L^p(\Omega)}+\|F\|_{L^q(\Omega)}+\|g\|_{L^\infty(\partial\Omega)}\big\},
 \end{equation}
 where $C$ depends on $\mu,\omega,\kappa,\lambda,m,d,p$ and $\Omega$.
\end{cor}

\begin{pf}
The conclusion is immediately derived from the Sobolev embedding theorem, and we omit the proof here (see \cite[Corollary 3.8]{QXS}).
\qed
\end{pf}

\section{Lipschitz Estimates}

\begin{lemma}\label{lemma:3.3}
Suppose that $A\in\Lambda(\mu,\tau,\kappa)$. Let $f=(f_i^\alpha)\in C^{0,\sigma}(\Omega;\mathbb{R}^{md})$, $F\in L^p(\Omega;\mathbb{R}^m)$
and $g\in C^{0,\sigma}(\partial\Omega;\mathbb{R}^m)$, where $p>d$ and $\sigma\in (0,1)$.
Then if $F$ and $g$ satisfy the compatibility condition $\int_\Omega F^\alpha dx + \int_{\partial\Omega} g^\alpha dS(x) = 0$,
the weak solution to $L_\varepsilon(u_\varepsilon) = \emph{div}(f)+ F$ in $\Omega$
and $\partial u_\varepsilon/\partial\nu_\varepsilon = g-n\cdot f$ on $\partial\Omega$ satisfies the uniform estimate
\begin{equation}\label{pri:3.7}
\|\nabla u_\varepsilon\|_{L^\infty(\Omega)} \leq C\big\{\|f\|_{C^{0,\sigma}(\Omega)}+\|F\|_{L^p(\Omega)}
+ \|g\|_{C^{0,\sigma}(\partial\Omega)}\big\},
\end{equation}
where $C$ depends on $\mu,\tau,\kappa,m,d,p,\sigma$ and $\Omega$.
\end{lemma}

\begin{pf}
 Since the estimate $\eqref{pri:3.7}$ has been proved in \cite[Theorem 1.2]{SZW4} in the case of $f= 0$,
 it reduces to prove the estimate $\eqref{pri:3.7}$ for the weak solution to the Neumman problem:
 \begin{equation}\label{f:3.4}
  L_\varepsilon(u_\varepsilon) = \text{div}(f) \quad\text{in}~~\Omega,
  \qquad \partial u_\varepsilon/\partial \nu_\varepsilon = -n\cdot f \quad\text{on}~~\partial\Omega.
 \end{equation}
 Owing to the fact that $u_\varepsilon$ in $\eqref{f:3.4}$ can be formulated by the Neumann matrix as follows
 \begin{equation*}
 \begin{aligned}
  u_\varepsilon(x) & = \int_{\Omega} N_\varepsilon(x,y)\text{div}(f)(y) dy - \int_{\partial\Omega} N_\varepsilon(x,y)n(y)\cdot f(y) dS \\
  & = -\int_\Omega \nabla_y N_\varepsilon(x,y)\cdot f(y) dy
 \end{aligned}
 \end{equation*}
 for any $x\in\Omega$, (for convenience, we omit the upper index.) and we have
 \begin{equation}\label{f:3.3}
 \begin{aligned}
 \nabla u_\varepsilon(x)
 &= -\int_\Omega \nabla_x\nabla_y N_\varepsilon(x,y)[f(y)-f(x)]dy - f(x)\int_{\Omega}\nabla_x \nabla_y N_\varepsilon(x,y)dy\\
 & = -\int_\Omega \nabla_x\nabla_y N_\varepsilon(x,y)[f(y)-f(x)]dy -f_k(x)\int_{\partial\Omega}\nabla_x N_\varepsilon(x,y)n_k(y)dS(y),
 \end{aligned}
 \end{equation}
 where the summation convention is used to the subscript $k$ from $1$ to $d$.
 To deal with the second term in the last equality of $\eqref{f:3.3}$, we construct the following elliptic systems
 \begin{equation*}
  \left\{\begin{aligned}
  L_\varepsilon(v_{\varepsilon,k}) &= 0 &\qquad \text{in} ~~~\Omega,\\
  \frac{~\partial v_{\varepsilon,k}}{\partial\nu_\varepsilon} & = n_kI &\qquad \text{on}~\partial\Omega
  \end{aligned}\right.
 \end{equation*}
 for $k=1,\cdots,d$, where $n_k$ is the $k$'th component of the unit outward normal vector to $\partial\Omega$.
 In view of \cite[Theorem 1.2]{SZW4}, we arrive at
 \begin{equation}\label{f:3.5}
  \|\nabla v_{\varepsilon,k}\|_{L^\infty(\Omega)}\leq C\|n\|_{C^{0,\eta}(\partial\Omega)} \leq C,
 \end{equation}
where we note that $\|n\|_{C^{0,\eta}(\partial\Omega)}$ is bounded by a constant depending on $\eta,M_0$.
We also point out that the integration of $n_k$ on $\partial\Omega$ vanishes, which guarantees the existence of $v_{\varepsilon,k}$.
In addition, $v_{\varepsilon,k}$ is a matrix-valued function.

 Proceeding as in the proof above, the transposed matrix of $v_{\varepsilon,k}$ which is denoted by $v_{\varepsilon,k}^*$ can be formulated by
 \begin{equation*}
  v_{\varepsilon,k}^*(x) = \int_{\partial\Omega} N_\varepsilon(x,y)n_k(y)dS(y).
 \end{equation*}
 Inserting the above expression into the formula $\eqref{f:3.3}$ gives
 \begin{equation*}
  \nabla u_\varepsilon(x)
   = -\int_\Omega \nabla_x\nabla_y N_\varepsilon(x,y)[f(y)-f(x)]dy - f_k(x)\nabla v_{\varepsilon,k}^*(x).
 \end{equation*}
 Then on account of $\eqref{pri:6.4}$ and $\eqref{f:3.5}$, we derive
 \begin{equation}\label{f:3.6}
  \|\nabla u_\varepsilon\|_{L^\infty(\Omega)}\leq C[f]_{C^{0,\sigma}(\Omega)} + \|f\|_{L^\infty(\Omega)}\|\nabla v^*_{\varepsilon}\|_{L^\infty(\Omega)}
  \leq C\|f\|_{C^{0,\sigma}(\Omega)}
 \end{equation}
 by noting the fact of $\|\nabla v_\varepsilon^*\|_{L^\infty(\Omega)}=\|\nabla v_\varepsilon\|_{L^\infty(\Omega)}$, where $C$ depends only on
 $\mu,\tau,\kappa,m,d,\sigma$ and $\Omega$.

 Combining the estimate $\eqref{f:3.6}$ and \cite[Theorem 1.2]{SZW4} leads to the desired estimate $\eqref{pri:3.7}$. The proof is done.
 \qed

\end{pf}

\begin{thm}\label{thm:3.1}
Let $\Psi_{\varepsilon,0}$ be given in $\eqref{pde:1.2}$. Suppose that $A\in \Lambda(\mu,\kappa,\tau)$, and $V$ satisfies $\eqref{a:2}$ and $\eqref{a:4}$. Then we obtain
\begin{equation}\label{pri:3.1}
\|\nabla\Psi_{\varepsilon,0}\|_{L^\infty(\Omega)} \leq C,
\qquad \|\Psi_{\varepsilon,0}-I\|_{L^\infty(\Omega)} \leq C\varepsilon \ln(r_0/\varepsilon+2),
\end{equation}
where $C$ depends only on $\mu,\tau,\kappa,d,m$ and $\Omega$.
\end{thm}

\begin{lemma}\label{lemma:3.1}
Suppose that $A\in\Lambda(\mu,\tau,\kappa)$. Let $u_\varepsilon$ be the solution to $L_\varepsilon(u_\varepsilon) = 0$ in $\Omega$ and
$\partial u_\varepsilon/\partial\nu_\varepsilon = g$  on $\partial\Omega$, where
\begin{equation*}
  g = \sum_{ij}\Big(n_i\frac{\partial}{\partial x_j} - n_j\frac{\partial}{\partial x_i}\Big) g_{ij},
\end{equation*}
and $g_{ij}\in C^1(\partial\Omega; \mathbb{R}^m)$. Then
\begin{equation}\label{pri:3.2}
 |\nabla u_\varepsilon(x)| \leq \frac{C}{\delta(x)} \sum_{ij}\|g_{ij}\|_{L^\infty(\partial\Omega)},
\end{equation}
for any $x\in\Omega$, where $\delta(x) = \emph{dist}(x,\partial\Omega)$, and $C$ depends only on $\mu,\kappa,\tau,d,m$ and $\Omega$.
\end{lemma}

\begin{pf}
The proof of this lemma can be found in \cite[pp.920]{SZW4}.
Here we can offer an easier proof based on a better estimate obtained in $\eqref{pri:6.4}$
due to the remarkable progress achieved in \cite{SZ,SZW12}. Proceeding as in the proof of \cite[Lemma 6.2]{SZW4},
we only need to show that
\begin{equation*}
 |u_\varepsilon(x) - u_\varepsilon(y)| \leq C \sum_{ij}\|g_{ij}\|_{L^\infty(\partial\Omega)}~,
\end{equation*}
where $|x-y| < cr$ and $r=\delta(x)$, since
the desired result straightforward comes from the interior Lipschitz estimates (see \cite[Lemma 16]{MAFHL}) for $L_\varepsilon$
\begin{equation*}
 |\nabla u_\varepsilon(x)| \leq  \frac{C}{\delta(x)}\Big(\dashint_{B(x,r)}|u_\varepsilon(y) - u_\varepsilon(x)|^2 dy\Big)^{1/2}.
\end{equation*}
According to the solution given by the related Neumann matrices, we have
\begin{eqnarray*}
 u_\varepsilon(x) - u_\varepsilon(y) & = & \int_{\partial\Omega} \big\{N_\varepsilon(x,z) - N_\varepsilon(y,z)\big\} g(z) dS(z) \\
 & = & -\sum_{ij}\int_{\partial\Omega} \big(n_i\frac{\partial}{\partial z_j} - n_j\frac{\partial}{\partial z_i}\big)
 \big\{N_\varepsilon(x,z)-N_\varepsilon(y,z)\big\}g_{ij}(z)dS(z),
\end{eqnarray*}
where we use the fact
that $n_i\frac{\partial}{\partial z_j} - n_j\frac{\partial}{\partial z_i}$ is a tangential derivative in the second equality
(see \cite[pp.921]{SZW4}), and this leads to
\begin{equation}\label{f:4.5}
\begin{aligned}
|u_\varepsilon(x) - u_\varepsilon(y)|
& \leq  2\int_{\partial\Omega}|\nabla_z\big\{N_\varepsilon(x,z) - N_\varepsilon(y,z)\big\}|dS(z)\sum_{ij}\|g_{ij}\|_{L^\infty(\partial\Omega)} \\
&\leq  2\sum_{ij}\|g_{ij}\|_{L^\infty(\partial\Omega)}\Big\{\int_{|z-Q|\leq cr}+\sum_{k=0}^\infty\int_{\Sigma_k}\Big\}\frac{|x-y|}{|x-z|^d} dS(z)
 =:2\sum_{ij}\|g_{ij}\|_{L^\infty(\partial\Omega)}\Big\{ I_1 + I_2\Big\}.
\end{aligned}
\end{equation}
where $Q\in\partial\Omega$ such that $|x-Q| = r$, and $\Sigma_k = \{z\in\partial\Omega~:~2^kcr\leq |z-Q|\leq 2^{k+1}cr\}$.
Note that we use the third estimate of $\eqref{pri:6.4}$ in the second inequality.

For $I_1$, it is clear to see that $|x-z|\geq r$ since $|z-Q|\leq cr$ and $|x-Q|=\delta(x)=r$. Thus
\begin{equation}\label{f:4.7}
 I_1 = \int_{|z-Q|\leq cr}\frac{|x-y|}{|x-z|^d} dS(z)
 \leq (cr)\cdot\frac{(cr)^{d-1}}{r^{d}}
 \leq C.
\end{equation}

For $I_2$, observing that $|x-z|\approx|z-Q|$ (in other words, they are comparable by triangle inequality)
whenever $z\in\Sigma_k$. Then we have
\begin{equation}\label{f:4.8}
I_2 = \sum_{k=0}^\infty\int_{\Sigma_k}\frac{|x-y|}{|x-z|^d} dS(z)
\leq \sum_{k=0}^\infty (cr)\cdot\frac{(2^{k+1}cr)^{d-1}}{(2^kcr)^{d}}
\leq C.
\end{equation}

Plugging $\eqref{f:4.7}$ and $\eqref{f:4.8}$ back into $\eqref{f:4.5}$,
the estimate $\eqref{pri:3.2}$ thus follows, and we complete the proof.
\qed
\end{pf}

\begin{lemma}\label{lemma:3.2}
Suppose that $A\in\Lambda(\mu,\kappa,\tau)$, and $V$ satisfies $\eqref{a:2}$ and $\eqref{a:4}$.
Then the weak solution $\Psi_{\varepsilon,0}$ to $\eqref{pde:1.2}$ satisfies the estimate
\begin{equation}\label{pri:3.3}
 \big|\nabla \Psi_{\varepsilon,0}(x)\big| \leq C\Big[1+\frac{\varepsilon}{\delta(x)}\Big]
\end{equation}
for any $x\in\Omega$, where $C$ depends only on $\mu, \tau, \kappa, d,m$ and $\Omega$.
\end{lemma}

\begin{pf}
Let
\begin{equation*}
  H_{\varepsilon,0}^{\alpha\beta}(x) = \Psi_{\varepsilon,0}^{\alpha\beta}(x) - \delta^{\alpha\beta} - \varepsilon\chi_{0}^{\alpha\beta}(x/\varepsilon).
\end{equation*}
By definition of $\Psi_{\varepsilon,0}$ and $\chi_{0}$ (see $\eqref{pde:1.2}$ and $\eqref{pde:6.1}$), we have $L_\varepsilon(H_{\varepsilon,0}) = 0$ in $\Omega$.
If we can verify that there exists $g_{ij}^{\alpha\beta}$ such that
\begin{equation}\label{f:3.1}
 \frac{\partial H_{\varepsilon,0}^{\alpha\beta}}{\partial \nu_\varepsilon}
 = \sum_{ij}\Big(n_i\frac{\partial}{\partial x_j} - n_j\frac{\partial}{\partial x_i}\Big)g_{ij}^{\alpha\beta}
 \quad \text{on}~~\partial\Omega, \qquad \text{and} \quad \|g_{ij}^{\alpha\beta}\|_{L^\infty(\partial\Omega)} \leq C\varepsilon,
\end{equation}
then the estimate $\eqref{pri:3.3}$ will follow from Lemma $\ref{lemma:3.1}$ directly.

Hence the main task now is to show $\eqref{f:3.1}$ is true. We first have
\begin{eqnarray*}
 \frac{\partial H_{\varepsilon,0}^{\alpha\gamma}}{\partial \nu_\varepsilon}
 = n_i\Big[\hat{V}_i^{\alpha\gamma} - V_{i}^{\alpha\gamma}(y)-a_{ij}^{\alpha\beta}(y)\frac{\partial\chi_0^{\beta\gamma}}{\partial y_j}\Big]
 = n_i b_{i0}^{\alpha\gamma}(y),
\end{eqnarray*}
where $y = x/\varepsilon$. Then from Lemma $\ref{lemma:6.4}$, it follows that there exist $E_{ji0}^{\alpha\gamma}$ such that
\begin{equation*}
  b_{i0}^{\alpha\gamma} = \frac{\partial}{\partial y_j}\big\{E_{ji0}^{\alpha\gamma}\big\}
  \qquad \text{and} \qquad
  E_{ji0}^{\alpha\gamma} =-E_{ij0}^{\alpha\gamma}.
\end{equation*}
We thus arrive at
\begin{equation}\label{f:3.2}
 \frac{\partial H_{\varepsilon,0}^{\alpha\gamma}}{\partial \nu_\varepsilon}
 = \frac{1}{2}\big(\sum_i n_i b_{i0}^{\alpha\gamma} + \sum_j n_j b_{j0}^{\alpha\gamma}\big)
 = \frac{\varepsilon}{2}\sum_{ij}\Big(n_i\frac{\partial}{\partial x_j} - n_j\frac{\partial}{\partial x_i}\Big)
 \big\{ E_{ji0}^{\alpha\gamma}(x/\varepsilon)\big\}.
\end{equation}

It remains to prove the inequality of $\eqref{f:3.1}$. Since $V$ satisfies $\eqref{a:2}$ and $\eqref{a:4}$, it is not hard to see that
the related corrector $\chi_0$ is a H\"older continuous function.
Therefore it follows from Lemma $\ref{lemma:6.4}$ that
$E_{ji0}^{\alpha\gamma}\in L^\infty(Y)$.

Finally, by setting $g_{ij}^{\alpha\gamma}(x)=\varepsilon E_{ji0}^{\alpha\gamma}(x/\varepsilon)$, we complete the proof.
\qed
\end{pf}

\begin{flushleft}
\textbf{Proof of Theorem \ref{thm:3.1}}\textbf{.}\quad We proceed to prove the first estimate of $\eqref{pri:3.1}$.
From the estimate $\eqref{pri:3.3}$, it immediately follows that
$\|\nabla\Psi_{\varepsilon,0}\|_{L^\infty(\Omega)}\leq C$, whenever $\delta(x) > 2\varepsilon$.

\quad~ If $\delta(x) \leq 2\varepsilon$ with $0<\varepsilon<1$, we will employ the blow-up method to approach this case.
Let $w^{\alpha\gamma}(x) = H_{\varepsilon,0}^{\alpha\gamma}(\varepsilon x)$. After a standard manipulation,
we derive $L_1(w) = 0$ in $\Omega_\varepsilon$ and
\end{flushleft}
\begin{equation*}
\frac{\partial w^{\alpha\gamma}}{\partial \nu_1} = n_i(\varepsilon x)a_{ij}^{\alpha\beta}(x)\frac{\partial w^{\beta\gamma}}{\partial x_j}
= \frac{\varepsilon}{2}\sum_{ij}\Big(n_i(\varepsilon x)\frac{\partial}{\partial x_j} - n_j(\varepsilon x)\frac{\partial}{\partial x_i}\Big)
\big\{E_{ji0}^{\alpha\gamma}(x)\big\} =: h(x)
\quad \text{on} ~ \partial\Omega_\varepsilon,
\end{equation*}
where $\Omega_\varepsilon = \{x\in\mathbb{R}^d:\varepsilon x\in\Omega\}$. Since $\partial\Omega\in C^{1,\eta}$ and $\eta\in[\tau,1)$,
the outward unit normal vector is H\"older continuous. This together with
$\|\nabla E_{ij0}^{\alpha\gamma}\|_{C^{0,\tau}(Y)}\leq C$, implies
$\|h\|_{C^{0,\tau}(\partial\Omega_\varepsilon\cap B(0,2))}\leq C\varepsilon$.
Hence from the Lipschitz estimate $\eqref{pri:2.6}$, it follows that
(for any $B(P,1)\cap\Omega_\varepsilon$ with $P\in\partial\Omega_\varepsilon$, we may assume $P = 0$ by translation)
\begin{eqnarray*}
\|\nabla w\|_{L^\infty(B(0,1)\cap\Omega_\varepsilon)} &\leq & C\Big\{\|h\|_{C^{0,\tau}(\partial\Omega_\varepsilon\cap B(0,2))}
+ \Big(\dashint_{B(0,2)\cap\Omega_\varepsilon}|\nabla w|^2 dx\Big)^{1/2}\Big\} \\
&\leq & C_p\Big\{\varepsilon+ \Big(\dashint_{B(0,2)\cap\Omega_\varepsilon}|\nabla w|^p dx\Big)^{1/p} \Big\},
\end{eqnarray*}
where $0<p<1$, and we use the convex property (see \cite[pp.184]{MGLM}) in the last inequality.
Set $H_{\varepsilon,0}(x) = \big(H_{\varepsilon,0}^{\alpha\gamma}(x)\big)$. By change of variable, we acquire
\begin{eqnarray}\label{f:4.34}
\|\nabla H_{\varepsilon,0}\|_{L^\infty(B(0,\varepsilon)\cap\Omega)}
\leq C_p \Big\{1+ \Big(\dashint_{B(0,2\varepsilon)\cap\Omega}|\nabla H_{\varepsilon,0}|^p dx\Big)^{1/p} \Big\}.
\end{eqnarray}
It remains to estimate the second term in the right hand side above. By substituting the estimate
\begin{equation}\label{f:4.9}
 |\nabla H_{\varepsilon,0}(x)| \leq \frac{C\varepsilon}{\delta(x)}
\end{equation}
into $\eqref{f:4.34}$, we obtain
\begin{eqnarray*}
  \Big(\dashint_{B(0,2\varepsilon)\cap\Omega}|\nabla H_{\varepsilon,0}|^p dx\Big)^{1/p}
  \leq C\varepsilon  \Big(\dashint_{B(0,2\varepsilon)\cap\Omega}\frac{dx}{\delta^p(x)}\Big)^{1/p}
  \leq C,
\end{eqnarray*}
where $C$ depends only on $\mu,\kappa,\tau,d,m,\eta$ and $\Omega$. Up to now, we have accomplished the first estimate in $\eqref{pri:3.1}$.

We now prove the remainder of $\eqref{pri:3.1}$. Due to \cite[Lemma 7.16]{DGNS}, it is not hard to see
\begin{equation}\label{f:4.11}
|H_{\varepsilon,0}(x) - H_{\varepsilon,0}(y)|\leq C\Big\{\int_\Omega \frac{|\nabla H_{\varepsilon,0}(z)|}{|x-z|^{d-1}}dz
+\int_\Omega\frac{|\nabla H_{\varepsilon,0}(z)|}{|y-z|^{d-1}}dz\Big\}
\end{equation}
for any $x,y\in\Omega$. Let $r_0$ be the diameter of $\Omega$. Clearly, we only need to estimate the first integral in the right hand-side,
and the other follows by the same way. In view of $\eqref{f:4.9}$, we have
\begin{equation}\label{f:4.10}
\begin{aligned}
\int_\Omega \frac{|\nabla H_{\varepsilon,0}(z)|}{|x-z|^{d-1}}dz
&\leq \int_{B(x,r_0)} \frac{|\nabla H_{\varepsilon,0}(z)|}{|x-z|^{d-1}}dz \\
&\leq C\varepsilon\int_{B(x,r_0-\varepsilon)}\frac{dz}{|x-z|^{d-1}\delta(z)}
+ C\int_{B(x,r_0)\setminus B(x,r_0-\varepsilon)} \frac{dz}{|x-z|^{d-1}}\\
&\leq C\Big\{\varepsilon\int_0^{r_0-\varepsilon}\int_{\partial B(x,r)}\frac{dS_rdr}{r^{d-1}(r_0-r)}
+ \int_{r_0-\varepsilon}^{r_0}\int_{\partial B(x,r)}\frac{dS_rdr}{r^{d-1}}\Big\}\leq C\big\{\varepsilon\ln{(r_0/\varepsilon)}+\varepsilon\big\},
\end{aligned}
\end{equation}
where we use the observation that $\Omega\subset B(x,r_0)$ and $\delta(z)$ actually becomes $r_0-|x-z|$,
(which do not weaken the singularities of the integrand.) We plug $\eqref{f:4.10}$ back into $\eqref{f:4.11}$, and have
\begin{equation}\label{f:4.12}
 |H_{\varepsilon,0}(x) - H_{\varepsilon,0}(y)|\leq C\varepsilon\ln{(r_0/\varepsilon+2)}.
\end{equation}
Note that $\Psi_{\varepsilon,0}$ defined in $\eqref{pde:1.2}$ is unique up to a constant.
If we assume that $\Psi_{\varepsilon,0}(x_0) = I$ for a fixed point $x_0\in\Omega$, then
$H_{\varepsilon,0}(x_0) = -\varepsilon\chi_0(x_0/\varepsilon)$.
So by setting $y = x_0$ in $\eqref{f:4.12}$,
it is clear to see $|H_{\varepsilon,0}(x)|\leq C\varepsilon\ln{(r_0/\varepsilon+2)}$ for any $x\in\Omega$, and the desired result follows.

As a comment, the result for $\varepsilon >1$ is in fact the case of the standard regularity theory (see $\eqref{pri:2.9}$).
The proof of the theorem is now complete.
\qed

\begin{cor}\label{cor:3.2}
Assume the same conditions as in Theorem $\ref{thm:3.1}$, then we have
\begin{equation}\label{pri:3.8}
\big[\Psi_{\varepsilon,0}-I\big]_{C^{0,\sigma}(\Omega)} \leq C[\varepsilon \ln(r_0/\varepsilon+2)]^{1-\sigma}
\end{equation}
for any $\sigma\in(0,1)$, where $C$ depends only on $\mu,\tau,\kappa,m,d,\sigma$ and $\Omega$.
\end{cor}
\begin{pf}
 In virtue of the interpolation inequality, we have
 \begin{equation*}
 \big[\Psi_{\varepsilon,0}-I\big]_{C^{0,\sigma}(\Omega)}\leq
 4r_0\|\Psi_{\varepsilon,0}-I\|_{L^\infty(\Omega)}^{1-\sigma}
 \|\nabla\Psi_{\varepsilon,0}\|_{L^\infty(\Omega)}^\sigma\leq C[\varepsilon \ln(r_0/\varepsilon+2)]^{1-\sigma},
 \end{equation*}
 where we use $\eqref{pri:3.1}$ in the last inequality.
 \qed
\end{pf}

\begin{cor}\label{cor:3.1}
Assume the same conditions as in Theorem $\ref{thm:3.1}$. then we have
\begin{equation}\label{pri:3.5}
\big\|\nabla \Psi_{\varepsilon,0}\big\|_{C^{0,\tau}(\Omega)} \leq C\max\{\varepsilon^{-\tau},1\}.
\end{equation}
Furthermore, $\Psi_{\varepsilon,0}^{-1}$ exists and satisfies the following estimates:
\begin{equation}\label{pri:3.6}
  2/3 \leq \big\|\Psi_{\varepsilon,0}^{-1}\big\|_{L^\infty(\Omega)} \leq 2, \qquad\qquad
  \big\|\nabla \Psi_{\varepsilon,0}^{-1}\big\|_{L^\infty(\Omega)} \leq C,
\end{equation}
whenever $\varepsilon\leq\varepsilon_0$, where $\varepsilon_0=\varepsilon_0(\mu,\tau,\kappa,d,m,\Omega)$ is sufficiently small.
\end{cor}

\begin{remark}
\emph{The proof of the estimates $\eqref{pri:3.5}$ and \eqref{pri:3.6} is quite similar to that given for \cite[Lemma 4.9]{QXS},
and is thus omitted here.}
\end{remark}

\begin{lemma}\label{lemma:5.7}
$(\text{A nonuniform estimates})$\textbf{.} Suppose $A\in\Lambda(\mu,\tau,\kappa)$,
$V$ satisfies $\eqref{a:2}$, $\eqref{a:4}$, and $B,c$
satisfy $\eqref{a:3}$.
Let $\sigma\in (0,\tau]$ and $p>d$. Assume $f\in C^{0,\sigma}(\Omega;\mathbb{R}^{md})$, $F\in L^p(\Omega;\mathbb{R}^m)$ and
$g\in C^{0,\sigma}(\partial\Omega;\mathbb{R}^m)$. Then
the weak solution to $\eqref{pde:1.1}$ satisfies the estimate
\begin{equation} \label{pri:5.10}
 \|\nabla u_\varepsilon\|_{C^{0,\frac{\sigma}{2}}(\Omega)}
 \leq C\max\{\varepsilon^{\frac{\sigma}{2}-1},1\}\big\{ \|f\|_{C^{0,\sigma}(\Omega)} + \|F\|_{L^p(\Omega)}
 + \|g\|_{C^{0,\sigma}(\partial\Omega)}\big\},
\end{equation}
where $C$ depends only on $\mu,\tau,\kappa,\lambda,p,d,m$ and $\Omega$.
\end{lemma}

\begin{pf}
The idea of this proof can be found in \cite[Lemma 4.10]{QXS}, and we provide with a proof here for the sake of completeness.
We focus on the case of $0<\varepsilon<1$,
since $\eqref{pri:5.10}$ follows directly from the Schauder estimate $\eqref{pri:2.12}$ when $\varepsilon\geq 1$. We first establish
\begin{equation}\label{f:2.6}
\|\nabla u_\varepsilon\|_{L^\infty(\Omega)}
\leq C\varepsilon^{\sigma-1}\big\{\|f\|_{C^{0,\sigma}(\Omega)}+\|F\|_{L^p(\Omega)}+\|g\|_{C^{0,\sigma}(\partial\Omega)}\big\}.
\end{equation}

To do so, let $v_\varepsilon^\beta = u_\varepsilon^\beta - u_\varepsilon^\beta(P)
- \varepsilon\chi_0^{\beta\gamma}(x/\varepsilon)u_\varepsilon^\gamma(P)$ for any $P\in \partial\Omega$.
By translation we may assume $P=0$, and have the new systems:
\begin{equation*}
\left\{\begin{aligned}
 \mathcal{L}_\varepsilon(v_\varepsilon) &= \text{div}(\tilde{f})
 + \tilde{F} &\qquad& \text{in}~~\Omega,\\
 \mathcal{B}_\varepsilon(v_\varepsilon) &= \tilde{g} - n\cdot\tilde{f} &\qquad& \text{on}~\partial\Omega,
\end{aligned}
\right.
\end{equation*}
where
\begin{equation*}
\begin{aligned}
&\tilde{f} = f+\varepsilon V_\varepsilon\chi_{0,\varepsilon}u_\varepsilon(0),\\
&\tilde{F} = F - \big[B_\varepsilon\nabla_y\chi_0+(c_\varepsilon +\lambda I)(I+\varepsilon\chi_0(x/\varepsilon))\big]u_\varepsilon(0),\\
&\tilde{g} = g - n\cdot V_\varepsilon u_\varepsilon(0) - n A_\varepsilon\nabla_y\chi_0 u_\varepsilon(0).
\end{aligned}
\end{equation*}
A routine computation gives rise
to
\begin{equation}\label{f:4.35}
\begin{aligned}
\mathcal{R}(\varepsilon;2D;2\Delta;\tilde{f};\tilde{F};\tilde{g})
& \leq C\mathcal{R}(1;D;\Delta;f;F;g)+C\|u_\varepsilon\|_{L^\infty(D)} \\
& \leq C\|u_\varepsilon\|_{L^\infty(\Omega)}
+ C\big\{\|f\|_{C^{0,\sigma}(\Omega)} + \|F\|_{L^p(\Omega)} + \|g\|_{C^{0,\sigma}(\partial\Omega)}\big\},
\end{aligned}
\end{equation}
where $\mathcal{R}$ is defined in $\eqref{f:2.12}$, and we note that
$[\tilde{g}]_{C^{0,\sigma}(\partial\Omega)}
\leq [g]_{C^{0,\sigma}(\partial\Omega)} + C(\mu,\tau,\kappa,m,d)\varepsilon^{-\sigma}|u_\varepsilon(0)|$.

Then in virtue of the Lipschitz estimate $\eqref{pri:2.8}$ on $\varepsilon$ scale, we acquire
\begin{equation*}
\begin{aligned}
\|\nabla v_\varepsilon\|_{L^\infty(D(0,\varepsilon))}
&\leq C\bigg\{\frac{1}{\varepsilon}\Big(\dashint_{2D}|v_\varepsilon|^2\Big)^{\frac{1}{2}}
+\mathcal{R}(\varepsilon;2D;2\Delta;\tilde{F};\tilde{f};\tilde{g})\bigg\} \\
& \leq C\varepsilon^{\sigma-1}[u_\varepsilon]_{C^{0,\sigma}(\Omega)} + C\|u_\varepsilon\|_{L^\infty(\Omega)}
+ C\big\{\|f\|_{C^{0,\sigma}(\Omega)} + \|F\|_{L^p(\Omega)} + \|g\|_{C^{0,\sigma}(\partial\Omega)}\big\} \\
& \leq C\varepsilon^{\sigma-1}\big\{\|f\|_{C^{0,\sigma}(\Omega)} + \|F\|_{L^p(\Omega)} + \|g\|_{C^{0,\sigma}(\partial\Omega)}\big\},
\end{aligned}
\end{equation*}
where we use $\eqref{f:4.35}$ in the second inequality, and $\eqref{pri:2.3}$ in the last one.
The covering argument (the interior estimate is established in \cite[Lemma 4.10]{QXS}) together with
\begin{equation*}
 \nabla u_\varepsilon = \nabla v_\varepsilon + \nabla_y\chi_0 u_\varepsilon(0)
\end{equation*}
leads to the estimate $\eqref{f:2.6}$.

Proceeding as in the proof above, we only consider boundary estimate, and it is not hard to derive the following estimate from $\eqref{pri:2.7}$:
\begin{equation*}
\begin{aligned}
 \big[\nabla v_\varepsilon\big]_{C^{0,\sigma}(D(\varepsilon))}
&\leq C\varepsilon^{-\sigma}\Big\{\frac{1}{\varepsilon}\Big(\dashint_{D(2\varepsilon)}|v_\varepsilon|^2\Big)^{\frac{1}{2}}
+ \mathcal{R}(\varepsilon;2D;2\Delta;\tilde{F};\tilde{f};\tilde{g}) \Big\} \\
& \leq C\varepsilon^{-1}[u_\varepsilon]_{C^{0,\sigma}(\Omega)} + C\varepsilon^{-\sigma}\big\{\|u_\varepsilon\|_{L^\infty(\Omega)}
+\|f\|_{C^{0,\sigma}(\Omega)}  + \|g\|_{C^{0,\sigma}(\partial\Omega)}\big\} + C\|F\|_{L^p(\Omega)}\\
& \leq C\varepsilon^{-1}\big\{\|f\|_{C^{0,\sigma}(\Omega)} + \|F\|_{L^p(\Omega)} + \|g\|_{C^{0,\sigma}(\partial\Omega)}\big\},
\end{aligned}
\end{equation*}
(where we also use $\eqref{f:4.35}$ in the second inequality and $\eqref{pri:2.3}$ in the last one) and this together with the corresponding interior estimate
(see \cite[Lemma 4.10]{QXS}) implies
\begin{equation}\label{f:2.7}
[\nabla u_\varepsilon]_{C^{0,\sigma}(\Omega)}
\leq C\varepsilon^{-1}\big\{\|f\|_{C^{0,\sigma}(\Omega)} + \|F\|_{L^p(\Omega)} + \|g\|_{C^{0,\sigma}(\partial\Omega)}\big\}.
\end{equation}

Incorporated with $\eqref{f:2.6}$ and $\eqref{f:2.7}$, a simple interpolation inequality indicates
\begin{equation*}
 [\nabla u_\varepsilon]_{C^{0,\frac{\sigma}{2}}(\Omega)}
 \leq 2\|\nabla u_\varepsilon\|_{L^\infty(\Omega)}^{\frac{1}{2}}[\nabla u_\varepsilon]_{C^{0,\sigma}(\Omega)}^{\frac{1}{2}}
 \leq C\varepsilon^{\frac{\sigma}{2}-1}\big\{\|f\|_{C^{0,\sigma}(\Omega)} + \|F\|_{L^p(\Omega)} + \|g\|_{C^{0,\sigma}(\partial\Omega)}\big\},
\end{equation*}
and then the estimate $\eqref{pri:5.10}$ follows. We now complete the proof.
\qed
\end{pf}

\begin{flushleft}
\textbf{Proof of Theorem \ref{thm:1.2}}\textbf{.}\quad
As we mentioned before, the estimate $\eqref{pri:1.2}$ immediately follows from the
standard Lipschitz estimate (see $\eqref{pri:2.9}$) for $\varepsilon\geq\varepsilon_*$, and we only need to consider the following transformation
\begin{equation}\label{T:1}
 u_\varepsilon^\beta(x) = \Psi_{\varepsilon,0}^{\beta\gamma}(x)v_\varepsilon^\gamma(x)
\end{equation}
for $\varepsilon<\varepsilon_*$, where $\varepsilon_* =\min\{\varepsilon_0,\varepsilon_1,\varepsilon_2\}$, and $\varepsilon_0$ is given in
Corollary $\ref{cor:3.1}$ and $\varepsilon_1,\varepsilon_2$ will be fixed later.
\end{flushleft}

On account of $\eqref{T:1}$, the Neumann problem $\eqref{pde:1.1}$ can be transformed into
\begin{equation}\label{f:3.7}
\left\{\begin{aligned}
L_\varepsilon(v_\varepsilon) & = \text{div}(\tilde{f}) + \tilde{F} &\qquad& \text{in}~~\Omega, \\
\frac{\partial v_\varepsilon}{\partial\nu_\varepsilon} &= \tilde{g} - n\cdot\tilde{f} &\qquad& \text{on}~\partial\Omega,
\end{aligned}\right.
\end{equation}
where
\begin{equation*}
\begin{aligned}
& \tilde{f}^\alpha = f^\alpha + A_\varepsilon^{\alpha\beta}(\Psi_{\varepsilon,0}^{\beta\gamma}-\delta^{\beta\gamma})\nabla v_\varepsilon^\gamma
+V_\varepsilon^{\alpha\beta}(\Psi_{\varepsilon,0}^{\beta\gamma}-\delta^{\beta\gamma})v_\varepsilon^\gamma, \\
& \tilde{F}^\alpha = F^\alpha + A_\varepsilon^{\alpha\beta}\nabla\Psi_{\varepsilon,0}^{\beta\gamma}\nabla v_\varepsilon^\gamma
+V_\varepsilon^{\alpha\gamma}\nabla v_\varepsilon^\gamma - B_\varepsilon^{\alpha\beta}\nabla u_\varepsilon^\beta
- (c_\varepsilon^{\alpha\beta}+\lambda\delta^{\alpha\beta})u_\varepsilon^\beta, \\
&\tilde{g}^\alpha = g^\alpha -n\cdot\widehat{V}^{\alpha\gamma}v_\varepsilon^\gamma.
\end{aligned}
\end{equation*}
In view of $\eqref{pde:1.2}$ and $\eqref{C:1}$, it is not hard to verify
\begin{equation*}
 \int_\Omega \tilde{F}^\alpha dx + \int_{\partial\Omega} \tilde{g}^\alpha dS(x) = 0.
\end{equation*}
Hence, applying Lemma $\ref{lemma:3.3}$ to $\eqref{f:3.7}$, we have
\begin{equation}\label{f:3.11}
 \|\nabla v_\varepsilon\|_{L^\infty(\Omega)}
 \leq C\big\{\|\tilde{f}\|_{C^{0,\nu}(\Omega)} + \|\tilde{F}\|_{L^p(\Omega)} + \|\tilde{g}\|_{C^{0,\nu}(\partial\Omega)}\big\},
\end{equation}
where $\nu$ will be given later.

Apparently, our task now is to estimate the right hand side of $\eqref{f:3.11}$.
Because $\tilde{f}$, $\tilde{F}$ and $\tilde{g}$ much involves $v_\varepsilon, u_\varepsilon$ and their derivatives,
based on the previous results obtained in the paper, we need to provide some quantity estimates related to them as follows.
It follows from Theorem $\ref{thm:1.1}$ and Corollary $\ref{cor:2.1}$ that
\begin{equation}\label{f:3.8}
\max\{\|u_\varepsilon\|_{W^{1,p}(\Omega)}, \|u_\varepsilon\|_{C^{0,\sigma}(\Omega)}\}
\leq C\big\{\|f\|_{L^\infty(\Omega)} + \|F\|_{L^p(\Omega)} + \|g\|_{L^\infty(\partial\Omega)}\big\}.
\end{equation}
This, together with $v_\varepsilon = \Psi_{\varepsilon,0}^{-1}u_\varepsilon$ and
$\nabla v_\varepsilon = \nabla(\Psi_{\varepsilon,0}^{-1})u_\varepsilon+\Psi_{\varepsilon,0}^{-1}\nabla u_\varepsilon$, leads to
\begin{equation}\label{f:3.9}
\max\{\|v_\varepsilon\|_{W^{1,p}(\Omega)}, \|v_\varepsilon\|_{C^{0,\sigma}(\Omega)}\}
\leq C\big\{\|f\|_{L^\infty(\Omega)} + \|F\|_{L^p(\Omega)} + \|g\|_{L^\infty(\partial\Omega)}\big\},
\end{equation}
where we use $\eqref{pri:3.6}$ in the above inequality.

We now in the position to approach the estimates associated with $\tilde{F}$, $\tilde{f}$ and $\tilde{g}$.

In view of
$\eqref{pri:3.1}$, $\eqref{f:3.8}$ and $\eqref{f:3.9}$, we first have
\begin{equation}\label{f:3.12}
\|\tilde{F}\|_{L^p(\Omega)}\leq C\big\{\|f\|_{L^\infty(\Omega)} + \|F\|_{L^p(\Omega)} + \|g\|_{L^\infty(\partial\Omega)}\big\}.
\end{equation}

For the estimate of $\|\tilde{f}\|_{C^{0,\nu}(\Omega)}$, we rewrite $\nabla v_\varepsilon$ as
\begin{equation*}
 \nabla v_\varepsilon = -\nabla\Psi_{\varepsilon,0}v_\varepsilon + (I-\Psi_{\varepsilon,0})\nabla v_\varepsilon + \nabla u_\varepsilon.
\end{equation*}
Then set $\nu =\min\{\tau,\sigma\}/2$, and it follows from $\eqref{pri:3.1}$, $\eqref{pri:3.8}$, $\eqref{pri:3.5}$,
$\eqref{pri:5.10}$ and $\eqref{f:3.9}$ that
\begin{equation*}
\begin{aligned}
\bigl[\nabla v_\varepsilon\bigr]_{C^{0,\nu}(\Omega)}
&\leq [\nabla \Psi_{\varepsilon,0}]_{C^{0,\nu}(\Omega)} \| v_\varepsilon\|_{L^\infty(\Omega)}
+\|\nabla\Psi_{\varepsilon,0}\|_{L^\infty(\Omega)}[v_\varepsilon]_{C^{0,\nu}(\Omega)} \\
\quad &+ [I-\Psi_{\varepsilon,0}]_{C^{0,\nu}(\Omega)} \|\nabla v_\varepsilon\|_{L^{\infty}(\Omega)}
+\|I-\Psi_{\varepsilon,0}\|_{L^\infty(\Omega)} [\nabla v_\varepsilon]_{C^{0,\nu}(\Omega)} + [\nabla u_\varepsilon]_{C^{0,\nu}(\Omega)} \\
\quad&\leq  C(\varepsilon^{-\tau}+\varepsilon^{\frac{\sigma}{2}-1})\big\{\|f\|_{C^{0,\sigma}(\Omega)} + \|F\|_{L^p(\Omega)}
+\|g\|_{C^{0,\sigma}(\partial\Omega)}\big\} + C\|\nabla v_\varepsilon\|_{L^\infty(\Omega)}
+  C\varepsilon^{\frac{1}{2}}[\nabla v_\varepsilon]_{C^{0,\nu}(\Omega)},
\end{aligned}
\end{equation*}
and this further indicates
\begin{equation}\label{f:3.10}
\big[\nabla v_\varepsilon\big]_{C^{0,\nu}(\Omega)}
\leq  C\varepsilon^{-\nu^\prime}\big\{\|f\|_{C^{0,\sigma}(\Omega)}+\|F\|_{L^p(\Omega)}+\|g\|_{C^{0,\sigma}(\partial\Omega)}\big\}
+C\|\nabla v_\varepsilon\|_{L^\infty(\Omega)},
\end{equation}
whenever $\varepsilon<\varepsilon_1$, where $\varepsilon_1=\min\{1/(2C)^2,1\}$ and $\nu^\prime =\max\{\tau,1-\sigma/2\}$. We thus have
\begin{equation}\label{f:3.13}
\begin{aligned}
\|\tilde{f}\|_{C^{0,\nu}(\Omega)} 
&\leq  \|f\|_{C^{0,\sigma}(\Omega)} + C\varepsilon\ln(\varepsilon/r_0+2)\|\nabla v_\varepsilon\|_{L^\infty(\Omega)}\\
& + \big[A_\varepsilon(\Psi_{\varepsilon,0}-I)\nabla v_\varepsilon\big]_{C^{0,\nu}(\Omega)}
+ \|V_\varepsilon(\Psi_{\varepsilon,0}-I)v_\varepsilon\|_{C^{0,\nu}(\Omega)} \\
&\leq  C\varepsilon^{1-\nu^\prime-\iota_0}\|\nabla v_\varepsilon\|_{L^\infty(\Omega)}
 + C\big\{\|f\|_{C^{0,\sigma}(\Omega)}+\|F\|_{L^p(\Omega)}+\|g\|_{C^{0,\sigma}(\partial\Omega)}\big\},
\end{aligned}
\end{equation}
where we use $\eqref{pri:3.1}$, $\eqref{pri:3.8}$, $\eqref{f:3.9}$ and $\eqref{f:3.10}$ in last inequality,
and choose a small positive number $\iota_0$ such that $\varepsilon^{\iota_0}\ln(r_0/\varepsilon+2)\leq 1$.

We now proceed to show the quantity of $\|\tilde{g}\|_{C^{0,\sigma}(\partial\Omega)}$, and it is not hard to
see that
\begin{equation}\label{f:3.14}
\begin{aligned}
\|\tilde{g}\|_{C^{0,\sigma}(\partial\Omega)}
&\leq \|g\|_{C^{0,\sigma}(\partial\Omega)} + C\|v_\varepsilon\|_{C^{0,\sigma}(\Omega)}\\
&\leq  C\big\{\|f\|_{C^{0,\sigma}(\Omega)}+\|F\|_{L^p(\Omega)}+\|g\|_{C^{0,\sigma}(\partial\Omega)}\big\}.
\end{aligned}
\end{equation}

Collecting $\eqref{f:3.11}$, $\eqref{f:3.12}$, $\eqref{f:3.13}$ and $\eqref{f:3.14}$ yields
\begin{equation*}
\|\nabla v_\varepsilon\|_{L^\infty(\Omega)} \leq C\varepsilon^{1-\nu^\prime-\iota_0}\|\nabla v_\varepsilon\|_{L^\infty(\Omega)}
+ C\big\{\|f\|_{C^{0,\sigma}(\Omega)}+\|F\|_{L^p(\Omega)}+\|g\|_{C^{0,\sigma}(\partial\Omega)}\big\},
\end{equation*}
and this implies
\begin{equation}\label{f:3.15}
\|\nabla v_\varepsilon\|_{L^\infty(\Omega)}\leq
C\big\{\|f\|_{C^{0,\sigma}(\Omega)}+\|F\|_{L^p(\Omega)}+\|g\|_{C^{0,\sigma}(\partial\Omega)}\big\},
\end{equation}
whenever $\varepsilon<\varepsilon_2$, where $\varepsilon_2 = \min\{1/(2C)^{\frac{1}{1-\nu^\prime-\iota_0}},1\}$.
Recalling the transformation $\eqref{T:1}$, it is straightforward to derive the estimate $\eqref{pri:1.2}$ from $\eqref{f:3.15}$.
Up to now we have proved this theorem.
\qed

\section{Convergence rates}

\begin{lemma}\label{lemma:4.1}
Suppose that $u_\varepsilon,u_0\in H^{1}(\Omega;\mathbb{R}^m)$ satisfy
$\myl{L}{\varepsilon}(u_\varepsilon)=\myl{L}{0}(u_0)$ in $\Omega$, and $\myl{B}{\varepsilon}(u_\varepsilon) = \myl{B}{0}(u_0)$ on
$\partial\Omega$. Let
\begin{equation*}
w_\varepsilon^\beta = u_\varepsilon^\beta - u_0^\beta
-\varepsilon\sum_{k=0}^d\chi_k^{\beta\gamma}(x/\varepsilon)\varphi_k^\gamma,
\end{equation*}
where $\varphi=(\varphi_k^\gamma)\in H^1(\mathbb{R}^d;\mathbb{R}^{md})$.  Then we have
\begin{equation}\label{pde:4.3}
\begin{aligned}
\big[\myl{L}{\varepsilon}(\mys{w})\big]^\alpha
& = -\frac{\partial}{\partial x_i}\left\{\mathcal{K}_i^\alpha
+\big[\hat{a}_{ij}^{\alpha\beta}-a_{ij,\varepsilon}^{\alpha\beta}\big]
\Big[\frac{\partial u_0^\beta}{\partial x_j}-\varphi_j^\beta\Big]
+\big[\hat{V}_i^{\alpha\beta}-V_{i,\varepsilon}^{\alpha\beta}\big]\big[u_0^\beta
-\varphi_0^\beta\big] - \varepsilon\big(\mathcal{I}_{i}^\alpha+\mathcal{J}_i^\alpha\big) \right\}\\
& +[\hat{B}_i^{\alpha\beta}-B_{i,\varepsilon}^{\alpha\beta}]
\big[\frac{\partial u_0^\beta}{\partial x_i}-\varphi^\beta_i\big]
+\big[\hat{c}^{\alpha\beta}-c_{\varepsilon}^{\alpha\beta}\big]\big[u_0^\beta-\varphi^\beta_0\big]
-\varepsilon\big(\mathcal{M}^\alpha + \mathcal{N}^\alpha\big)
\end{aligned}
\end{equation}
and
\begin{equation}\label{pde:4.4}
\begin{aligned}
\big[\myl{B}{\varepsilon}(w_\varepsilon)\big]^{\alpha}
= n_i\Big\{\mathcal{K}_i^\alpha
+\big[\hat{a}_{ij}^{\alpha\beta}-a_{ij,\varepsilon}^{\alpha\beta}\big]
\Big[\frac{\partial u_0^\beta}{\partial x_j}-\varphi^\beta_j\Big]
+\big[\hat{V}_i^{\alpha\beta}-V_{i,\varepsilon}^{\alpha\beta}\big]\big[u_0^\beta
-\varphi_0^\beta\big] - \varepsilon\mathcal{I}_{i}^\alpha\Big\},
\end{aligned}
\end{equation}
where $n_i$ denotes the i'th component of the unit normal vector to $\partial\Omega$, and
\begin{equation}\label{eq:1}
\begin{aligned}
& \mathcal{I}_i^\alpha
 = a_{ij,\varepsilon}^{\alpha\beta}\sum_{k=0}^d\chi_{k,\varepsilon}^{\beta\gamma}
 \frac{\partial}{\partial x_j}\big\{\varphi^\gamma_k\big\}
+V_{i,\varepsilon}^{\alpha\beta}\sum_{k=0}^d\chi_{k,\varepsilon}^{\beta\gamma}
\varphi^\gamma_k,
\qquad
\mathcal{J}_{i}^\alpha = \sum_{k=0}^d\frac{\partial \vartheta_k^{\alpha\gamma}}{\partial y_i}\varphi^\gamma_k,
\qquad
\mathcal{K}_{i}^\alpha =  \sum_{j=0}^db_{ij,\varepsilon}^{\alpha\gamma}\varphi^\gamma_j, \\
& \mathcal{M}^\alpha =
\sum_{k=0}^d\Big[\frac{\partial \vartheta_k^{\alpha\gamma}}{\partial y_i}
+ B_{i,\varepsilon}^{\alpha\beta}\chi_{k,\varepsilon}^{\beta\gamma}\Big]
\frac{\partial}{\partial x_i}\big\{\varphi_k^\gamma\big\},
\qquad
\mathcal{N}^\alpha =
\big[c_{\varepsilon}^{\alpha\beta}+\lambda\delta^{\alpha\beta}\big]
 \sum_{k=0}^d\chi_{k,\varepsilon}^{\beta\gamma}\varphi^\gamma_k,
 \qquad y =x/\varepsilon.
\end{aligned}
\end{equation}
\end{lemma}

\begin{pf}
In view of $\myl{L}{\varepsilon}(u_\varepsilon) = \myl{L}{0}(u_0)$, we have
\begin{equation}\label{f:5.20}
\myl{L}{\varepsilon}(w_\varepsilon) = \myl{L}{0}(u_0) -\myl{L}{\varepsilon}(u_0)
- \myl{L}{\varepsilon}\big(\varepsilon\sum_{k=0}^d\chi_{k,\varepsilon}\varphi_k\big)
\end{equation}
where
\begin{equation}\label{f:5.21}
\begin{aligned}
\big[\myl{L}{0}(u_0)\big]^\alpha &= -\frac{\partial}{\partial x_i}\Big\{\hat{a}_{ij}^{\alpha\beta}\frac{\partial u_0^\beta}{\partial x_j}
+ \hat{V}_i^{\alpha\beta}u_0^\beta\Big\} + \hat{B}_i^{\alpha\beta}\frac{\partial u_0^\beta}{\partial x_i}
+ [\hat{c}^{\alpha\beta}+\lambda\delta^{\alpha\beta}]u_0^\beta,\\
-\big[\myl{L}{\varepsilon}(u_0)\big]^\alpha & = \frac{\partial}{\partial x_i}\Big\{a_{ij,\varepsilon}^{\alpha\beta}\frac{\partial u_0^\beta}{\partial x_j}
+ V_{i,\varepsilon}^{\alpha\beta}u_0^\beta\Big\} - B_{i,\varepsilon}^{\alpha\beta}\frac{\partial u_0^\beta}{\partial x_i}
- [c_\varepsilon^{\alpha\beta}+\lambda\delta^{\alpha\beta}]u_0^\beta,
\end{aligned}
\end{equation}
and
\begin{equation}\label{f:5.22}
\begin{aligned}
- \big[\myl{L}{\varepsilon}\big(\varepsilon\chi_{k,\varepsilon}\varphi_k\big)\big]^\alpha
& = \frac{\partial}{\partial x_i}\Big\{a_{ij,\varepsilon}^{\alpha\beta}\frac{\partial\chi_k^{\beta\gamma}}{\partial y_j}
\varphi^\gamma_k+\varepsilon a_{ij,\varepsilon}^{\alpha\beta}\chi_{k,\varepsilon}^{\beta\gamma}
\frac{\partial}{\partial x_j}\big\{\varphi^\gamma_k\big\}
+\varepsilon V_{i,\varepsilon}^{\alpha\beta}\chi_{k,\varepsilon}^{\beta\gamma}\varphi^\gamma_k\Big\}\\
&- B_{i,\varepsilon}^{\alpha\beta}\frac{\partial\chi_k^{\beta\gamma}}{\partial y_i}\varphi^\gamma_k
-\varepsilon B_{i,\varepsilon}^{\alpha\beta}\chi_{k,\varepsilon}^{\beta\gamma}\frac{\partial}{\partial x_i}\big\{\varphi^\gamma_k\big\}
-\varepsilon[c^{\alpha\beta}_\varepsilon+\lambda\delta^{\alpha\beta}]\chi_{k,\varepsilon}^{\beta\gamma}\varphi^\gamma_k
\end{aligned}
\end{equation}
with summation convention applied to $k$ from $0$ to $d$.
By substituting $\eqref{f:5.21}$ and $\eqref{f:5.22}$ into $\eqref{f:5.20}$, we obtain
\begin{equation}\label{f:4.16}
\begin{aligned}
\big[\myl{L}{\varepsilon}(w_\varepsilon)\big]^\alpha
&= -\frac{\partial}{\partial x_i}\Big\{\hat{a}_{ij}^{\alpha\beta}\frac{\partial u_0^\beta}{\partial x_j}
-a_{ij,\varepsilon}^{\alpha\beta}\frac{\partial u_0^\beta}{\partial x_j}
- a_{ij,\varepsilon}^{\alpha\beta}\frac{\partial\chi_l^{\beta\gamma}}{\partial y_j}\varphi_l^\gamma
+\hat{V}_i^{\alpha\beta}u_0^\beta - V_{i,\varepsilon}^{\alpha\beta}u_0^\beta
-a_{ij,\varepsilon}^{\alpha\beta}\frac{\partial\chi_0^{\beta\gamma}}{\partial y_j}\varphi_0^\gamma\\
& - \varepsilon\myl{I}{i}^\alpha\Big\}
+ \hat{B}_i^{\alpha\beta}\frac{\partial u_0^\beta}{\partial x_i} - B_{i,\varepsilon}^{\alpha\beta}\frac{\partial u_0^\beta}{\partial x_i}
- B_{i,\varepsilon}^{\alpha\beta}\frac{\partial\chi_l^{\beta\gamma}}{\partial y_i}\varphi_l^\gamma
+\hat{c}^{\alpha\beta}u_0^\beta - c_\varepsilon^{\alpha\beta}u_0^\beta
- B_{i,\varepsilon}^{\alpha\beta}\frac{\partial\chi_0^{\beta\gamma}}{\partial y_i}\varphi_0^\gamma\\
& -\varepsilon B_{i,\varepsilon}^{\alpha\beta}\sum_{k=0}^d\chi_k^{\beta\gamma}\frac{\partial}{\partial x_i}\big\{\varphi_k^\gamma\big\}
-\varepsilon[c_\varepsilon^{\alpha\beta} + \lambda\delta^{\alpha\beta}]
\sum_{k=0}^d\chi_{k,\varepsilon}^{\beta\gamma}\varphi^\gamma_k,
\end{aligned}
\end{equation}
where $l=1,\ldots,d$. By definition of $\eqref{def:6.1}$ and $\eqref{def:6.2}$, we have
\begin{equation}\label{f:4.17}
\begin{aligned}
&\hat{a}_{ij}^{\alpha\beta}\frac{\partial u_0^\beta}{\partial x_j}
-a_{ij,\varepsilon}^{\alpha\beta}\frac{\partial u_0^\beta}{\partial x_j}
- a_{ij,\varepsilon}^{\alpha\beta}\frac{\partial\chi_l^{\beta\gamma}}{\partial y_j}\varphi_l^\gamma
 = b_{ij,\varepsilon}^{\alpha\gamma}\varphi_j^\gamma
+ [\hat{a}_{ij}^{\alpha\beta}-a_{ij,\varepsilon}^{\alpha\beta}]\Big[\frac{\partial u_0^\beta}{\partial x_j}-\varphi_j^\beta\Big]\\
&\hat{V}_i^{\alpha\beta}u_0^\beta - V_{i,\varepsilon}^{\alpha\beta}u_0^\beta
-a_{ij,\varepsilon}^{\alpha\beta}\frac{\partial\chi_0^{\beta\gamma}}{\partial y_j}\varphi_0^\gamma
= b_{i0,\varepsilon}^{\alpha\gamma}\varphi_0^\gamma
+ [\hat{V}_i^{\alpha\beta} -V_{i,\varepsilon}^{\alpha\beta}][u_0^\beta -\varphi_0^\beta]
\end{aligned}
\end{equation} and
\begin{equation}\label{f:4.18}
\begin{aligned}
&\hat{B}_i^{\alpha\beta}\frac{\partial u_0^\beta}{\partial x_i} - B_{i,\varepsilon}^{\alpha\beta}\frac{\partial u_0^\beta}{\partial x_i}
- B_{i,\varepsilon}^{\alpha\beta}\frac{\partial\chi_l^{\beta\gamma}}{\partial y_i}\varphi_l^\gamma
= \mathrm{W}_{l,\varepsilon}^{\alpha\gamma}\varphi_l^\gamma
+ [\hat{B}_i^{\alpha\beta}-B_{i,\varepsilon}^{\alpha\beta}]\Big[\frac{\partial u_0^\beta}{\partial x_i} - \varphi_i^\beta\Big]\\
&\hat{c}^{\alpha\beta}u_0^\beta - c_\varepsilon^{\alpha\beta}u_0^\beta
- B_{i,\varepsilon}^{\alpha\beta}\frac{\partial\chi_0^{\beta\gamma}}{\partial y_i}\varphi_0^\gamma
=\mathrm{W}_{0,\varepsilon}^{\alpha\gamma}\varphi_0^\gamma
+ [\hat{c}^{\alpha\beta}-c_\varepsilon^{\alpha\beta}][u_0^\beta-\varphi_0^\beta].
\end{aligned}
\end{equation}
Then collecting $\eqref{f:4.16}$,$\eqref{f:4.17}$ and $\eqref{f:4.18}$ gives
\begin{equation*}
\begin{aligned}
\big[\myl{L}{\varepsilon}(w_\varepsilon)\big]^\alpha
&= -\frac{\partial}{\partial x_i}\Big\{\sum_{k=0}^d b_{ik,\varepsilon}^{\alpha\gamma}\varphi_k^\gamma
+ [\hat{a}_{ij}^{\alpha\beta}-a_{ij,\varepsilon}^{\alpha\beta}]\Big[\frac{\partial u_0^\beta}{\partial x_j}-\varphi_j^\beta\Big]
+ [\hat{V}_i^{\alpha\beta} -V_{i,\varepsilon}^{\alpha\beta}][u_0^\beta -\varphi_0^\beta]
- \varepsilon\mathcal{I}_i^\alpha\Big\}\\
&+ [\hat{B}_i^{\alpha\beta}-B_{i,\varepsilon}^{\alpha\beta}]\Big[\frac{\partial u_0^\beta}{\partial x_i} - \varphi_i^\beta\Big]
-\varepsilon B_{i,\varepsilon}^{\alpha\beta}\sum_{k=0}^d\chi_k^{\beta\gamma}\frac{\partial}{\partial x_i}\big\{\varphi_k^\gamma\big\}
+\sum_{k=0}^d\mathrm{W}_{k,\varepsilon}^{\alpha\gamma}\varphi_k^\gamma
+ [\hat{c}^{\alpha\beta}-c_\varepsilon^{\alpha\beta}][u_0^\beta-\varphi_0^\beta]\\
&-\varepsilon[c_\varepsilon^{\alpha\beta} + \lambda\delta^{\alpha\beta}]
\sum_{k=0}^d\chi_{k,\varepsilon}^{\beta\gamma}\varphi^\gamma_k.
\end{aligned}
\end{equation*}
This together with
\begin{eqnarray*}
\sum_{k=0}^d\mathrm{W}_{k,\varepsilon}^{\alpha\gamma}\varphi_k^\gamma
 = \varepsilon\sum_{k=0}^d\frac{\partial}{\partial x_i}\Big\{\frac{\partial \vartheta_k^{\alpha\gamma}}{\partial y_i}\varphi_k^\gamma\Big\}
-\varepsilon\sum_{k=0}^d\frac{\partial \vartheta_k^{\alpha\gamma}}{\partial y_i}\frac{\partial}{\partial x_i}\big\{\varphi_k^\gamma\big\}
\end{eqnarray*}
leads to the expression $\eqref{pde:4.3}$.

We now turn to prove $\eqref{pde:4.4}$. Note that $\myl{B}{\varepsilon}(u_\varepsilon) = \myl{B}{0}(u_0)$ on $\partial\Omega$, and we have
\begin{equation}\label{f:4.19}
\begin{aligned}
\myl{B}{\varepsilon}(w_\varepsilon) = \myl{B}{0}(u_0) - \myl{B}{\varepsilon}(u_0)
- \myl{B}{\varepsilon}\big(\varepsilon\sum_{k=0}^d\chi_{k,\varepsilon}\varphi_k\big)
\end{aligned}
\end{equation}
where
\begin{equation}\label{f:5.23}
\big[\myl{B}{0}(u_0)\big]^\alpha = n_i\hat{V}_i^{\alpha\beta}u_0^\beta + n_i\hat{a}_{ij}^{\alpha\beta}\frac{\partial u_0^\beta}{\partial x_j},
\qquad
\big[\myl{B}{\varepsilon}(u_0)\big]^\alpha = n_iV_{i,\varepsilon}^{\alpha\beta}u_0^\beta + n_ia_{ij,\varepsilon}^{\alpha\beta}
\frac{\partial u_0^\beta}{\partial x_j},
\end{equation}
and
\begin{equation}\label{f:5.24}
\begin{aligned}
&\big[\myl{B}{\varepsilon}\big(\varepsilon\chi_{k,\varepsilon}\varphi_k\big)\big]^\alpha
=\varepsilon n_iV_{i,\varepsilon}^{\alpha\beta}\chi_{k,\varepsilon}^{\beta\gamma}\varphi_k^\gamma
+n_ia_{ij,\varepsilon}^{\alpha\beta}\frac{\partial\chi_k^{\beta\gamma}}{\partial y_j}\varphi_k^\gamma
+\varepsilon n_i a_{ij,\varepsilon}^{\alpha\beta}\chi_{k,\varepsilon}^{\beta\gamma}
\frac{\partial}{\partial x_j}\big\{\varphi_k^\gamma\big\}.
\end{aligned}
\end{equation}
for $k=0,1,\ldots,d$.
Then inserting the expression $\eqref{f:5.23}$ and $\eqref{f:5.24}$ into $\eqref{f:4.19}$ indicates
\begin{equation*}
\big[\myl{B}{\varepsilon}(w_\varepsilon)\big]^\alpha
= n_i\Big\{\hat{V}_i^{\alpha\beta}u_0^\beta- V_{i,\varepsilon}^{\alpha\beta}u_0^\beta
- a_{ij,\varepsilon}^{\alpha\beta}\frac{\partial\chi_0^{\beta\gamma}}{\partial y_j}\varphi_0^\gamma
+\hat{a}_{ij}^{\alpha\beta}\frac{\partial u_0^\beta}{\partial x_j}
-a_{ij,\varepsilon}^{\alpha\beta}\frac{\partial u_0^\beta}{\partial x_j}
-a_{ij,\varepsilon}^{\alpha\beta}\frac{\partial\chi_l^{\beta\gamma}}{\partial y_j}\varphi_l^\gamma
-\varepsilon\mathcal{I}_i^\alpha\Big\},
\end{equation*}
where summation convention is used to $l$ from $1$ to $d$.
This together with $\eqref{f:4.17}$ finally gives the expression $\eqref{pde:4.4}$, and we complete the proof.
\qed
\end{pf}

\begin{lemma}\label{lemma:4.4}
Let $u_\varepsilon$ and $u_0$ be given in Lemma $\ref{lemma:4.1}$, and $w_\varepsilon$ is defined in $\eqref{f:4.13}$ and satisfies
$\eqref{pde:4.3}$ and $\eqref{pde:4.4}$. Then we have
\begin{equation}\label{pri:4.5}
\begin{aligned}
\|w_{\varepsilon}\|_{H^1(\Omega)}
&\leq C\Big\{\|h_{\varepsilon}\vec{\phi}\|_{L^2(\Omega\setminus\Sigma_{2\varepsilon})}
+ \|\nabla u_0 - \vec{\varphi}\|_{L^2(\Omega)}
+\|u_0 - \varphi_0\|_{L^2(\Omega)} \\
&+\varepsilon\|h_{\varepsilon}\nabla\vec{\phi}\|_{L^2(\Omega)}
+ \varepsilon \|h_{\varepsilon}\vec{\phi}\|_{L^2(\Omega)}\Big\}
+C\varepsilon\|h_\varepsilon\vec{\phi}\|_{L^2(\partial\Omega)},
\end{aligned}
\end{equation}
where $\vec{\varphi}=(\varphi_1,\cdots,\varphi_d)$, $\vec{\phi} = (\varphi_0,\vec{\varphi})$,
and $h$ represents the periodic function depending on some of the periodic functions shown before
such as the coefficients of $\myl{L}{\varepsilon}$,
the correctors $\{\chi_{k}\}_{k=0}^{d}$, and auxiliary functions
$\{b_{ik}, E_{jik},\nabla\vartheta_k\}_{k=0}^d$.
\end{lemma}

\begin{pf}
Let $w_\varepsilon=w_{\varepsilon,1}+w_{\varepsilon,2}$, where $w_{\varepsilon,1}$ and $w_{\varepsilon,2}$
satisfy
\begin{equation}\label{pde:4.1}
\left\{\begin{aligned}
\mys{\mathcal{L}}(w_{\varepsilon,1})
& = \mys{\mathcal{L}}(\mys{w})
+ \text{div}(\mathcal{K}) &\quad \text{in} ~~~\Omega,\\
\mys{\mathcal{B}}(w_{\varepsilon,1})
& = \mys{\mathcal{B}}(\mys{w})- n\cdot(\mathcal{K}) &\quad \text{on} ~\partial\Omega,
\end{aligned}\right.
\end{equation}
and
\begin{equation}\label{pde:4.2}
\left\{\begin{aligned}
\mys{\mathcal{L}}(w_{\varepsilon,2})
& = - \text{div}(\mathcal{K}) &\quad \text{in} ~~~\Omega,\\
\mys{\mathcal{B}}(w_{\varepsilon,2})
& = n\cdot\mathcal{K} &\quad \text{on} ~\partial\Omega,
\end{aligned}\right.
\end{equation}
respectively, where $\mathcal{K}=(\mathcal{K}_{i}^\alpha)$ is defined in $\eqref{eq:1}$.

For the first equation in $\eqref{pde:4.1}$, it follows from $\eqref{pri:2.4}$ that
\begin{equation}\label{f:4.14}
\|w_{\varepsilon,1}\|_{H^1(\Omega)}
\leq C\big\{\|\nabla u_0 - \vec{\varphi}\|_{L^2(\Omega)}
+ \|u_0 - \varphi_0\|_{L^2(\Omega)}+\varepsilon\|h_{\varepsilon}\nabla\vec{\phi}\|_{L^2(\Omega)}
+ \varepsilon \|h_\varepsilon\vec{\phi}\|_{L^2(\Omega)}
+\varepsilon\|h_\varepsilon\vec{\phi}\|_{L^2(\partial\Omega)}\big\},
\end{equation}
where $h$ depends on the coefficients of $\myl{L}{\varepsilon}$,
the correctors $\{\chi_{k}\}_{k=0}^{d}$ and auxiliary functions $\{\nabla\vartheta_k\}_{k=0}^d$.


We now focus on $\eqref{pde:4.2}$. In view of $\eqref{pde:1.4}$, we have
\begin{equation}\label{f:4.1}
\mathrm{B}_\varepsilon[w_{\varepsilon,2},v] = \int_\Omega \mathcal{K}\cdot\nabla v dx =: R(v).
\end{equation}
According to Lemma $\ref{lemma:6.4}$, $R(v)$ in $\eqref{f:4.1}$ satisfies
\begin{eqnarray*}
R(v)& = & \varepsilon\int_\Omega\sum_{k=0}^d\Big\{\frac{\partial}{\partial x_j}\big[{E}_{jik,\varepsilon}^{\alpha\gamma}\big]
\varphi^\gamma_k\Big\}\frac{\partial v^\alpha}{\partial x_i}dx\\
&=&\varepsilon\int_\Omega
\frac{\partial}{\partial x_j}\Big\{\sum_{k=0}^d{E}_{jik,\varepsilon}^{\alpha\gamma}\varphi^\gamma_k\Big\}
\frac{\partial v^\alpha}{\partial x_i} dx
-\varepsilon\int_\Omega
\sum_{k=0}^d\Big\{{E}_{jik,\varepsilon}^{\alpha\gamma}
\frac{\partial}{\partial x_j}\big[\varphi^\gamma_k\big]\Big\}
\frac{\partial v^\alpha}{\partial x_i} dx\\
&=:& R_1(v) - R_2(v).
\end{eqnarray*}
Note that due to the antisymmetry of $E_{jik}$ with respect to $i,j$, we obtain
\begin{eqnarray*}
R_1(v) &=& \varepsilon\int_\Omega\frac{\partial}{\partial x_j}
\Big\{\big[\psi_\varepsilon+(1-\psi_\varepsilon)\big]{E}_{jik,\varepsilon}^{\alpha\gamma}\varphi^\gamma_k
\Big\}\frac{\partial v^\alpha}{\partial x_i} dx\\
&=&\varepsilon\int_\Omega\frac{\partial}{\partial x_j}
\Big\{(1-\psi_\varepsilon){E}_{jik,\varepsilon}^{\alpha\gamma}\varphi^\gamma_k\Big\}\frac{\partial v^\alpha}{\partial x_i} dx
-~\varepsilon\int_\Omega\psi_\varepsilon{E}_{jik,\varepsilon}^{\alpha\gamma}\varphi^\gamma_k
\frac{\partial^2 v^\alpha}{\partial x_i\partial x_j} dx \\
&=&\varepsilon\int_\Omega\frac{\partial}{\partial x_j}
\Big\{(1-\psi_\varepsilon){E}_{jik,\varepsilon}^{\alpha\gamma}\varphi^\gamma_k
\Big\}\frac{\partial v^\alpha}{\partial x_i} dx,
\end{eqnarray*}
where $\psi_\varepsilon\in C_0^\infty(\Omega)$ satisfies $\eqref{f:4.20}$, and $k=0,1,\ldots,d$. Moreover, we have
\begin{equation*}
\begin{aligned}
R_1(v) & = - \int_\Omega
\frac{\partial\psi_\varepsilon}{\partial x_j}{E}_{jik,\varepsilon}^{\alpha\gamma}\varphi^\gamma_k
\frac{\partial v^\alpha}{\partial x_i} dx
+\int_\Omega (1-\psi_\varepsilon)b_{ik,\varepsilon}^{\alpha\gamma}\varphi^\gamma_k
\frac{\partial v^\alpha}{\partial x_i} dx
 + \varepsilon \int_\Omega(1-\psi_\varepsilon){E}_{jik,\varepsilon}^{\alpha\gamma}
\frac{\partial\varphi^\gamma_k}{\partial x_j}\frac{\partial v^\alpha}{\partial x_i} dx,
\end{aligned}
\end{equation*}
and this indicates
\begin{equation}\label{f:4.2}
\begin{aligned}
\big|R_1(v)\big| & \leq \big\{\|{E}_{\varepsilon}\vec{\phi}\|_{L^2(\Omega\setminus\Sigma_{2\varepsilon})}
+\|b_{\varepsilon}\vec{\phi}\|_{L^2(\Omega\setminus\Sigma_{2\varepsilon})}
+ \varepsilon\|{E}_{\varepsilon}\nabla\vec{\phi}\|_{L^2(\Omega\setminus\Sigma_{2\varepsilon})}\big\}
\|\nabla v\|_{L^2(\Omega\setminus\Sigma_{2\varepsilon})}.
\end{aligned}
\end{equation}
Meanwhile we arrive at
\begin{equation}\label{f:4.3}
 \big|R_2(v)\big| \leq \varepsilon\|{E}_\varepsilon\nabla\vec{\phi}\|_{L^2(\Omega)}\|\nabla v\|_{L^2(\Omega)}.
\end{equation}

Let $v = w_{\varepsilon,2}$ in the expression $\eqref{f:4.1}$. Combining $\eqref{f:4.1}$, $\eqref{f:4.2}$, $\eqref{f:4.3}$ and $\eqref{f:2.3}$,
we acquire
\begin{equation}\label{f:4.6}
\begin{aligned}
\|w_{\varepsilon,2}\|_{H^1(\Omega)}
\leq C\big\{\|{E}_{\varepsilon}\vec{\phi}\|_{L^2(\Omega\setminus\Sigma_{2\varepsilon})}
+\|b_{\varepsilon}\vec{\phi}\|_{L^2(\Omega\setminus\Sigma_{2\varepsilon})}
+\varepsilon\|{E}_{\varepsilon}\nabla\vec{\phi}\|_{L^2(\Omega)}\big\}.
\end{aligned}
\end{equation}
Then the conclusion of this lemma immediately comes from $\eqref{f:4.14}$ and $\eqref{f:4.6}$, and we complete the proof.
\qed
\end{pf}

As a reminder, $S_\varepsilon, \bar{S}_\varepsilon,\psi_r$ is defined in $\eqref{f:6.5}$, $\eqref{def:6.3}$ and $\eqref{f:4.20}$, respectively.

\begin{lemma}\label{lemma:4.5}
Let $\Omega$ be a bounded $C^{1,1}$ domain. Suppose that the coefficients of $\myl{L}{\varepsilon}$ satisfy $\eqref{a:1}$,
$\eqref{a:2}$ and $\eqref{a:3}$. Assume that $u_\varepsilon,u_0$ are the weak solutions to $\eqref{pde:1.5}$
with $F\in L^2(\Omega;\mathbb{R}^m)$ and $g\in B^{1/2,2}(\partial\Omega;\mathbb{R}^m)$. Then we have
\begin{equation}\label{pri:4.12}
\big\|u_\varepsilon - u_0
-\varepsilon\chi_{0,\varepsilon}S_\varepsilon(\psi_{4\varepsilon}u_0)
-\varepsilon\chi_{k,\varepsilon}S_\varepsilon(\psi_{4\varepsilon}\nabla_k u_0)\big\|_{H^1(\Omega)}
\leq C\varepsilon^{\frac{1}{2}}\big\{\|F\|_{L^2(\Omega)}+\|g\|_{B^{1/2,2}(\partial\Omega)}\big\},
\end{equation}
where $C$ depends only on $\mu,\kappa,m,d$ and $\Omega$.
\end{lemma}

\begin{pf}
By choosing $\varphi_0 = \mys{S}(\psi_{4\varepsilon}u_0)$ and $\varphi_k = \mys{S}(\psi_{4\varepsilon}\nabla_k u_0)$ in $\eqref{f:4.13}$, we let
\begin{equation*}
 w_\varepsilon = u_\varepsilon - u_0
-\varepsilon\chi_{0,\varepsilon}S_\varepsilon(\psi_{4\varepsilon}u_0)
-\varepsilon\chi_{k,\varepsilon}S_\varepsilon(\psi_{4\varepsilon}\nabla_k u_0).
\end{equation*}
Then it follows from Lemma $\ref{lemma:4.4}$ that
\begin{equation}\label{f:6.10}
\begin{aligned}
 \|w_\varepsilon\|_{H^1(\Omega)}
 &\leq C \bigg\{ \| u_0 - S_\varepsilon(\psi_{4\varepsilon} u_0)\|_{L^2(\Omega)}
 + \|\nabla u_0 - S_\varepsilon(\psi_{4\varepsilon}\nabla u_0)\|_{L^2(\Omega)}
 + \varepsilon\|h_\varepsilon S_\varepsilon(\psi_{4\varepsilon}u_0)\|_{L^2(\Omega)}  \\
 & + \varepsilon\|h_\varepsilon S_\varepsilon(\psi_{4\varepsilon}\nabla u_0)\|_{L^2(\Omega)}
   + \varepsilon\|h_\varepsilon \nabla S_\varepsilon(\psi_{4\varepsilon}u_0)\|_{L^2(\Omega)}
   + \varepsilon\|h_\varepsilon \nabla S_\varepsilon(\psi_{4\varepsilon}\nabla u_0)\|_{L^2(\Omega)}
 \bigg\}.
\end{aligned}
\end{equation}
We note that $S_\varepsilon(\psi_{4\varepsilon} u_0)$ and $S_\varepsilon(\psi_{4\varepsilon} \nabla u_0)$
is supported in $\Sigma_{3\varepsilon}$. To complete the proof, we need the following estimates.

Due to $\eqref{pri:6.2}$, we have
\begin{equation}\label{f:6.6}
\begin{aligned}
 \| u_0 - S_\varepsilon(\psi_{4\varepsilon} u_0)\|_{L^2(\Omega)}
 &\leq \| (1-\psi_{4\varepsilon})u_0\|_{L^2(\Omega)}
  + \| \psi_{4\varepsilon}u_0 - S_\varepsilon(\psi_{4\varepsilon} u_0)\|_{L^2(\Omega)} \\
 & \leq \|u_0\|_{L^2(\Omega\setminus\Sigma_{8\varepsilon})}
 + C\varepsilon\|\nabla(\psi_{4\varepsilon}u_0)\|_{L^2(\Omega)} \\
 & \leq C\|u_0\|_{L^2(\Omega\setminus\Sigma_{8\varepsilon})} + C\varepsilon\|\nabla u_0\|_{L^2(\Sigma_{4\varepsilon})}
\end{aligned}
\end{equation}
and
\begin{equation}\label{f:6.7}
 \| \nabla u_0 - S_\varepsilon(\psi_{4\varepsilon} \nabla u_0)\|_{L^2(\Omega)}
  \leq C\|\nabla u_0\|_{L^2(\Omega\setminus\Sigma_{8\varepsilon})} + C\varepsilon\|\nabla^2 u_0\|_{L^2(\Sigma_{4\varepsilon})}.
\end{equation}
From Lemma $\ref{lemma:6.1}$, it follows that
\begin{equation}\label{f:6.8}
\begin{aligned}
\|h_\varepsilon S_\varepsilon(\psi_{4\varepsilon}u_0)\|_{L^2(\Omega)} +
\|h_\varepsilon \nabla S_\varepsilon(\psi_{4\varepsilon}u_0)\|_{L^2(\Omega)}
& \leq C\|\psi_{4\varepsilon}u_0\|_{L^2(\Omega)} + C\|\nabla(\psi_{4\varepsilon}u_0)\|_{L^2(\Omega)}  \\
& \leq C\|u_0\|_{L^2(\Sigma_{4\varepsilon})} + C\varepsilon^{-1}\|u_0\|_{L^2(\Omega\setminus\Sigma_{8\varepsilon})}
+ C\|\nabla u_0\|_{L^2(\Sigma_{4\varepsilon})}
\end{aligned}
\end{equation}
and
\begin{equation}\label{f:6.9}
\|h_\varepsilon S_\varepsilon(\psi_{4\varepsilon}\nabla u_0)\|_{L^2(\Omega)} +
\|h_\varepsilon \nabla S_\varepsilon(\psi_{4\varepsilon}\nabla u_0)\|_{L^2(\Omega)}
\leq C\|\nabla u_0\|_{L^2(\Sigma_{4\varepsilon})} + C\varepsilon^{-1}\|\nabla u_0\|_{L^2(\Omega\setminus\Sigma_{8\varepsilon})}
+ C\|\nabla^2 u_0\|_{L^2(\Sigma_{4\varepsilon})}.
\end{equation}
Plugging the estimates $\eqref{f:6.6}$, $\eqref{f:6.7}$, $\eqref{f:6.8}$ and $\eqref{f:6.9}$ into $\eqref{f:6.10}$, we obtain
\begin{equation}\label{f:6.11}
\|w_\varepsilon\|_{H^1(\Omega)}
\leq C\|u_0\|_{H^1(\Omega\setminus\Sigma_{8\varepsilon})} + C\varepsilon\|u_0\|_{H^2(\Sigma_{4\varepsilon})}.
\end{equation}

Set $\varrho = (\varrho_1,\cdots,\varrho_d)\in C_0^1(\mathbb{R}^d;\mathbb{R}^d)$ be a
vector field such that $\big<\varrho,n\big> \geq c >0$ on $\partial\Omega$, where $n$ denotes the outward unit normal vector to $\partial\Omega$.
Since the divergence theorem, we have
\begin{equation}\label{f:6.25}
\begin{aligned}
c\int_{\partial\Omega} \big(|\nabla u_0|^2 + |u_0|^2\big) dS
&\leq \int_{\partial\Omega}<\varrho,n> (|\nabla u_0|^2 + |u_0|^2)dS\\
&= \int_{\Omega} \text{div}(\varrho)\big(|\nabla u_0|^2+|u_0|^2\big) dx
+ 2\int_{\Omega} \varrho\big(\nabla^2 u_0 \nabla u_0 + \nabla u_0 u_0\big)dx \\
&\leq C\big\{\|u_0\|_{L^2(\Omega)}^2 + \|\nabla u_0\|_{L^2(\Omega)}^2 + \|\nabla^2 u_0\|_{L^2(\Omega)}\big\}
\leq C\|u_0\|_{H^2(\Omega)}^2
\end{aligned}
\end{equation}
By the same argument, it is not hard to see
\begin{equation*}
 \Big(\int_{S_t} \big(|\nabla u_0|^2 + |u|^2 \big)dS_t\Big)\leq C\|u_0\|_{H^2(\Sigma_{t})}^2\leq C\|u_0\|_{H^2(\Omega)}^2
\end{equation*}
holds for any $t\in(0,8\varepsilon)$, where $S_t$ and $\Sigma_t$ are defined in Remark $\ref{re:2.2}$, and $C$ does not depend on $t$.
Then by $\eqref{f:2.11}$, we obtain
\begin{equation}\label{f:6.24}
\begin{aligned}
 \|u_0\|_{H^1(\Omega\setminus\Sigma_{8\varepsilon})}
 & = \Big(\int_{\Omega\setminus\Sigma_{8\varepsilon}}\big(|\nabla u_0|^2 + |u|^2 \big)dx\Big)^{1/2}\\
 & = \Big(\int_0^{8\varepsilon}\int_{S_t} \big(|\nabla u_0|^2 + |u|^2 \big)dS_t dt\Big)^{1/2}
 \leq C\varepsilon^{\frac{1}{2}}\|u_0\|_{H^2(\Omega)}.
\end{aligned}
\end{equation}
This together with $\eqref{f:6.11}$ gives
\begin{equation*}
\|w_\varepsilon\|_{H^1(\Omega)}\leq C\varepsilon^{\frac{1}{2}}\|u_0\|_{H^2(\Omega)}
\leq C\varepsilon^{\frac{1}{2}}\big\{\|F\|_{L^2(\Omega)} + \|g\|_{B^{1/2,2}(\partial\Omega)}\big\}.
\end{equation*}
We note that the assumption of $\partial\Omega\in C^{1,1}$ is only used in the last inequality, and we are done.
\qed

\end{pf}

\begin{remark}\label{rm:5.1}
\emph{We replace $w_\varepsilon$ in the proof of Lemma $\ref{lemma:4.5}$ into
\begin{equation}\label{f:5.32}
w_\varepsilon = u_\varepsilon - u_0 -\varepsilon\chi_{0,\varepsilon}\bar{S}_\varepsilon(\myu{u})
-\varepsilon\chi_{k,\varepsilon}\bar{S}_\varepsilon(\nabla_k \myu{u}),
\end{equation}
where $\myu{u}$ is the extension of $u_0$ such that $\myu{u} = u_0$ on $\Omega$, and $\myu{u}\in H^2(\mathbb{R}^d;\mathbb{R}^m)$
satisfying $\|\myu{u}\|_{H^2(\mathbb{R}^d)}\leq C\|u_0\|_{H^2(\Omega)}$.
Also, $\bar{S}_\varepsilon$ is given in the sense of $\eqref{def:6.3}$.
It follows from $\eqref{pri:4.5}$ that
\begin{equation*}
\begin{aligned}
\|w_{\varepsilon}\|_{H^1(\Omega)}
&\leq C\Big\{\|{h}_{\varepsilon}\mys{\bar{S}}(\myu{u})\|_{L^2(\Omega\setminus\Sigma_{2\varepsilon})}
+\|{h}_{\varepsilon}\mys{\bar{S}}(\nabla \myu{u})\|_{L^2(\Omega\setminus\Sigma_{2\varepsilon})}
+\|\myu{u} - \mys{\bar{S}}(\myu{u})\|_{L^2(\Omega)}\\
&+ \|\nabla \myu{u} - \mys{\bar{S}}(\nabla \myu{u})\|_{L^2(\Omega)}
+\varepsilon\|h_{\varepsilon}\mys{\bar{S}}(\nabla^2\myu{u})\|_{L^2(\Omega)}
+ \varepsilon \|h_{\varepsilon}\mys{\bar{S}}(\nabla \myu{u})\|_{L^2(\Omega)} \\
& + \varepsilon \|h_{\varepsilon}\mys{\bar{S}}(\myu{u})\|_{L^2(\Omega)}
+\varepsilon\|h_\varepsilon\mys{\bar{S}}(\myu{u})\|_{L^2(\partial\Omega)}
+\varepsilon\|h_\varepsilon\mys{\bar{S}}(\nabla \myu{u})\|_{L^2(\partial\Omega)}\Big\},
\end{aligned}
\end{equation*}
where we need the counterparts of $\eqref{pri:6.1}$ and $\eqref{pri:6.2}$ for $\bar{S}_\varepsilon$(or see Remark $\ref{re:4.3}$).
The different task from the proof of Lemma $\ref{lemma:4.5}$ is to estimate
the terms of $h_\varepsilon \bar{S}_\varepsilon(\tilde{u})$ and $h_\varepsilon \bar{S}_\varepsilon(\nabla\tilde{u})$ in
$L^2(\Omega\setminus\Sigma_{2\varepsilon})$ and $L^{2}(\partial\Omega)$ by Lemma $\ref{lemma:6.3}$.
Then we can also arrive at
\begin{equation}\label{f:4.45}
\big\|w_\varepsilon\big\|_{H^1(\Omega)}
\leq C\varepsilon^{\frac{1}{2}}\big\{\|F\|_{L^2(\Omega)}+\|g\|_{B^{1/2,2}(\partial\Omega)}\big\}
\end{equation}
without any difficulty, which actually gives the counterpart of \cite[Theorem 4.2]{TS2} in our case. }
\end{remark}

\begin{lemma}\label{lemma:6.5}
Let $\Omega$ be a bounded $C^{1,1}$ domain. Suppose that the coefficients of $\myl{L}{\varepsilon}$ satisfy $\eqref{a:1}$,
$\eqref{a:2}$ and $\eqref{a:3}$. Assume that $u_\varepsilon,u_0$ are the weak solutions to $\eqref{pde:1.5}$
with $F\in L^2(\Omega;\mathbb{R}^m)$ and $g\in B^{1/2,2}(\partial\Omega;\mathbb{R}^m)$. Then we have
\begin{equation}\label{pri:6.7}
\big\|u_\varepsilon - u_0
-\varepsilon\chi_{0,\varepsilon}S_\varepsilon(\psi_{4\varepsilon}\tilde{u}_0)
-\varepsilon\chi_{k,\varepsilon}S_\varepsilon(\psi_{4\varepsilon}\nabla_k \tilde{u}_0)\big\|_{L^2(\Omega)}
\leq C\varepsilon\big\{\|F\|_{L^2(\Omega)}+\|g\|_{B^{1/2,2}(\partial\Omega)}\big\},
\end{equation}
where $\tilde{u}_0\in H^2(\mathbb{R}^d;\mathbb{R}^m)$ is the extension of $u_0$, and $C$ depends only on $\mu,\kappa,m,d$ and $\Omega$.
\end{lemma}

\begin{pf}
We will employ duality argument.
For any $\Phi\in L^2(\Omega;\mathbb{R}^m)$, let $\phi_\varepsilon$ and $\phi_0$ satisfy
\begin{equation}\label{pde:4.8}
\left\{\begin{aligned}
\myl{L}{\varepsilon}^*(\phi_\varepsilon) &= \Phi &\quad&\text{in}~~\Omega,\\
\myl{B}{\varepsilon}^*(\phi_\varepsilon) &= 0 &\quad&\text{on}~\partial\Omega,
\end{aligned}\right.
\qquad\quad
\left\{\begin{aligned}
\myl{L}{0}^*(\phi_0) &= \Phi &\quad&\text{in}~~\Omega,\\
\myl{B}{0}^*(\phi_0) &= 0 &\quad&\text{on}~\partial\Omega,
\end{aligned}\right.
\end{equation}
respectively. Due to $\eqref{pri:4.12}$, we can have
\begin{equation}\label{f:6.23}
\big\|\phi_\varepsilon - \phi_0 -\varepsilon\chi_0^*(x/\varepsilon)S_\varepsilon(\psi_{10\varepsilon}\phi_0)
-\varepsilon\chi_k^*(x/\varepsilon)S_\varepsilon(\psi_{10\varepsilon}\nabla_k\phi_0)\big\|_{H^{1}(\Omega)}
\leq C\varepsilon^{\frac{1}{2}}\|\Phi\|_{L^2(\Omega)},
\end{equation}
where $\chi_k^*$ are correctors associated with $\mathcal{L}_\varepsilon^*$,
and satisfy the same estimate as $\chi_k$ do, where $k=0,\cdots,d$.

Let
\begin{equation}\label{f:6.31}
 w_\varepsilon = u_\varepsilon - u_0
-\varepsilon\chi_{0,\varepsilon}S_\varepsilon(\psi_{4\varepsilon}\tilde{u}_0)
-\varepsilon\chi_{k,\varepsilon}S_\varepsilon(\psi_{4\varepsilon}\nabla_k \tilde{u}_0),
\end{equation}
where $\tilde{u}_0\in H^2(\mathbb{R}^d)$ is the extended function of $u_0$ such that $\tilde{u}_0 = u_0$ in $\Omega$
and $\|\tilde{u}_0\|_{H^2(\mathbb{R}^d)}\leq C\|u_0\|_{H^2(\Omega)}$.
According to the formula $\eqref{f:2.10}$ and Lemma $\ref{lemma:4.1}$, we obtain
\begin{equation}\label{f:6.14}
 \int_\Omega w_\varepsilon \Phi dx =
 \int_\Omega f\cdot\nabla \phi_\varepsilon dx + \int_\Omega F \phi_\varepsilon dx,
\end{equation}
where
\begin{equation*}
\begin{aligned}
 & f = \mathcal{K} - \varepsilon(\mathcal{I}+\mathcal{J}) + [\hat{A} - A_\varepsilon]\big[\nabla u_0 - S_\varepsilon(\psi_{4\varepsilon}\nabla \tilde{u}_0)\big]
 + [\hat{V} - V_\varepsilon]\big[u_0-S_\varepsilon(\psi_{4\varepsilon}\tilde{u}_0)\big], \\
 & F = -\varepsilon\big(\mathcal{M} + \mathcal{N}\big) + [\hat{B}-B_{\varepsilon}]
\big[\nabla u_0 - S_\varepsilon(\psi_{4\varepsilon}\nabla \tilde{u}_0)\big]
+\big[\hat{c}-c_{\varepsilon}\big]\big[u_0-S_\varepsilon(\psi_{4\varepsilon}\tilde{u}_0)\big].
\end{aligned}
\end{equation*}
Below we do some calculations in more details. By the fact that $\myu{u}$ is the extension of $u_0$, we have
\begin{equation*}
  \nabla u_0 - S_\varepsilon(\psi_{4\varepsilon}\nabla\myu{u})
= \nabla\myu{u} - S_\varepsilon(\nabla\myu{u})
+ S_\varepsilon\big((1-\psi_{4\varepsilon})\nabla\myu{u}\big)
- (1-\psi_{4\varepsilon})\nabla \myu{u}
+ (1-\psi_{4\varepsilon})\nabla u_0
\quad \text{in}~ \Omega,
\end{equation*}
and then
\begin{equation}\label{f:6.15}
\begin{aligned}
&\bigg|\int_\Omega [\hat{A}  -A_\varepsilon][\nabla u_0 - S_\varepsilon(\psi_{4\varepsilon}\nabla\myu{u})]\cdot\nabla\phi_\varepsilon dx \bigg|
\leq C  \Bigg\{\int_\Omega |(1-\psi_{4\varepsilon})\nabla u_0||\nabla\phi_\varepsilon|dx \\ 
& \qquad + \int_\Omega |\nabla\tilde{u}_0 -S_\varepsilon(\nabla\myu{u})||\nabla\phi_\varepsilon|dx
+ \int_\Omega \big|(1-\psi_{4\varepsilon})\nabla\myu{u}-S_\varepsilon\big((1-\psi_{4\varepsilon})\nabla\myu{u}\big)\big|
|\nabla\phi_\varepsilon|dx \Bigg\} \\
&\leq C\bigg\{\|\nabla u_0\|_{L^2(\Omega\setminus\Sigma_{8\varepsilon})}\|\nabla \phi_\varepsilon\|_{L^2(\Omega\setminus\Sigma_{8\varepsilon})}
+\|\nabla\myu{u}-S_\varepsilon(\nabla\myu{u})\|_{L^2(\mathbb{R}^d)}\|\nabla\phi_\varepsilon\|_{L^2(\Omega)} \\
& \qquad + \|(1-\psi_{4\varepsilon})\nabla\myu{u}-S_\varepsilon\big((1-\psi_{4\varepsilon})\nabla\myu{u}\big)\|_{L^2(\mathbb{R}^d)}
\|\nabla\phi_\varepsilon\|_{L^2(\Omega\setminus\Sigma_{9\varepsilon})}\bigg\} \\
&\leq C\|\nabla u_0\|_{L^2(\Omega\setminus\Sigma_{9\varepsilon})}\|\nabla \phi_\varepsilon\|_{L^2(\Omega\setminus\Sigma_{9\varepsilon})}
+C\varepsilon\|\nabla^2\myu{u}\|_{L^2(\mathbb{R}^d)}\|\nabla\phi_\varepsilon\|_{L^2(\Omega)} \\
& \qquad + C\varepsilon\|\nabla[(1-\psi_{4\varepsilon})\nabla\myu{u}]\|_{L^2(\mathbb{R}^d)}
\|\nabla\phi_\varepsilon\|_{L^2(\Omega\setminus\Sigma_{9\varepsilon})}\\
&\leq C\|\nabla u_0\|_{L^2(\Omega\setminus\Sigma_{9\varepsilon})}\|\nabla\phi_\varepsilon\|_{L^2(\Omega\setminus\Sigma_{9\varepsilon})}
+ C\varepsilon\|u_0\|_{H^2(\Omega)}\|\nabla\phi_\varepsilon\|_{L^2(\Omega)},
\end{aligned}
\end{equation}
where we use Cauchy inequality in the second inequality
and the observation that
$S_\varepsilon\big((1-\psi_{4\varepsilon})\nabla\myu{u}\big)$ restricted to $\Omega$ is supported in $\Omega\setminus\Sigma_{9\varepsilon}$.
In the third one, we employ the estimate $\eqref{pri:6.2}$,
and the last one follows from
\begin{equation*}
\|\nabla[(1-\psi_{4\varepsilon})\nabla\myu{u}]\|_{L^2(\mathbb{R}^d)}
\leq C\varepsilon^{-1}\|\nabla \myu{u}\|_{L^2(\Omega\setminus\Sigma_{9\varepsilon})}
+ C\|\nabla^2\myu{u}\|_{L^2(\mathbb{R}^d)}
\leq C\varepsilon^{-1}\|\nabla u_0\|_{L^2(\Omega\setminus\Sigma_{9\varepsilon})}
+ C\|u_0\|_{H^2(\Omega)}.
\end{equation*}
Similarly, by noting $u_0-S_\varepsilon(\psi_{4\varepsilon}\myu{u})=\myu{u}-S_\varepsilon(\myu{u})
+ S_\varepsilon\big((1-\psi_{4\varepsilon})\myu{u}\big)
-(1-\psi_{4\varepsilon})\myu{u}
+(1-\psi_{4\varepsilon})u_0$ in $\Omega$, we obtain
\begin{equation}\label{f:6.16}
\Big|\int_\Omega  [\hat{V} - V_\varepsilon]\big[u_0-S_\varepsilon(\psi_{4\varepsilon}\tilde{u}_0)\big]\cdot \nabla\phi_\varepsilon dx\Big|
\leq C\|u_0\|_{L^2(\Omega\setminus\Sigma_{9\varepsilon})}\|\nabla\phi_\varepsilon\|_{L^2(\Omega\setminus\Sigma_{9\varepsilon})}
+ C\varepsilon\|u_0\|_{H^2(\Omega)}\|\nabla\phi_\varepsilon\|_{L^2(\Omega)}.
\end{equation}

As demonstrated before, it is not hard to derive
\begin{equation}\label{f:6.17}
\begin{aligned}
&\Big|\int_\Omega [\hat{B}-B_{\varepsilon}]
\big[\nabla u_0 - S_\varepsilon(\psi_{4\varepsilon}\nabla \tilde{u}_0)\big]\phi_\varepsilon dx\Big|
\leq C\|\nabla u_0\|_{L^2(\Omega\setminus\Sigma_{9\varepsilon})}\|\phi_\varepsilon\|_{L^2(\Omega\setminus\Sigma_{9\varepsilon})}
+ C\varepsilon\|u_0\|_{H^2(\Omega)}\|\phi_\varepsilon\|_{L^2(\Omega)}, \\
&\Big|\int_\Omega\big[\hat{c}-c_{\varepsilon}\big]\big[u_0-S_\varepsilon(\psi_{4\varepsilon}\tilde{u}_0)\big]
\phi_\varepsilon dx\Big|
\leq C\|u_0\|_{L^2(\Omega\setminus\Sigma_{9\varepsilon})}\|\phi_\varepsilon\|_{L^2(\Omega\setminus\Sigma_{9\varepsilon})}
+ C\varepsilon\|u_0\|_{H^2(\Omega)}\|\phi_\varepsilon\|_{L^2(\Omega)}.
\end{aligned}
\end{equation}

Now, for simplicity of presentation,
let $h_{ijk,\varepsilon}^{\alpha\gamma}(x) = a_{ij,\varepsilon}^{\alpha\beta}\chi_{k,\varepsilon}^{\beta\gamma}$ with $k=0,\cdots,d$.
Set $\nabla_i = \frac{\partial}{\partial x_i}, \nabla^2_{ij} = \frac{\partial^2}{\partial x_i\partial x_j}$,
and the notation of ``$\nabla_0 \myu{u} $'' means $\myu{u}$ itself.
Note that
\begin{equation}\label{f:6.13}
\begin{aligned}
\int_\Omega h_{ijk,\varepsilon}^{\alpha\gamma}
\nabla_j\big\{S_\varepsilon(\psi_{4\varepsilon}\nabla_k\myu{u}^\gamma)\big\} \nabla_i\phi_\varepsilon^\alpha dx
 = \int_\Omega h_{ijk,\varepsilon}^{\alpha\gamma}
S_\varepsilon\big(\psi_{4\varepsilon}\nabla_{jk}^2\myu{u}^\gamma\big)\nabla_i\phi_\varepsilon^\alpha dx
+\int_\Omega h_{ijk,\varepsilon}^{\alpha\gamma}
S_\varepsilon\big(\nabla_j\psi_{4\varepsilon}\nabla_k\myu{u}^\gamma\big)\nabla_i\phi_\varepsilon^\alpha dx,
\end{aligned}
\end{equation}
and then we have
\begin{equation}\label{f:6.12}
\begin{aligned}
\bigg|\int_\Omega h_{ijk,\varepsilon}^{\alpha\gamma}
&\nabla_j\big\{S_\varepsilon(\psi_{4\varepsilon}\nabla_k\myu{u}^\gamma)\big\} \nabla_i\phi_\varepsilon^\alpha dx\bigg| \\
& \leq  \|h_{ijk,\varepsilon}^{\alpha\gamma}
S_\varepsilon\big(\psi_{4\varepsilon}\nabla_{jk}^2\myu{u}^\gamma\big)\|_{L^2(\mathbb{R}^d)}
\|\nabla_i\phi_\varepsilon^\alpha\|_{L^2(\Omega)}
+ \|h_{ijk,\varepsilon}^{\alpha\gamma}
S_\varepsilon\big(\nabla_j\psi_{4\varepsilon}\nabla_{k}\myu{u}^\gamma\big)\|_{L^2(\mathbb{R}^d)}
\|\nabla_i\phi_\varepsilon^\alpha\|_{L^2(\Omega\setminus\Sigma_{9\varepsilon})}\\
&\leq C\big(\|\nabla^2\myu{u}\|_{L^2(\mathbb{R}^d)}+\|\nabla\myu{u}\|_{L^2(\mathbb{R}^d)}\big)\|\nabla\phi_\varepsilon\|_{L^2(\Omega)}
+ C\varepsilon^{-1}\|\myu{u}\|_{H^1(\Omega\setminus\Sigma_{9\varepsilon})}\|\nabla\phi_\varepsilon\|_{L^2(\Omega\setminus\Sigma_{9\varepsilon})} \\
&\leq C\|u_0\|_{H^2(\Omega)}\|\nabla\phi_\varepsilon\|_{L^2(\Omega)}
+ C\varepsilon^{-1}\|u_0\|_{H^1(\Omega\setminus\Sigma_{9\varepsilon})}\|\nabla\phi_\varepsilon\|_{L^2(\Omega\setminus\Sigma_{9\varepsilon})},
\end{aligned}
\end{equation}
In the first inequality, we mention that
$S_\varepsilon\big(\nabla_j\psi_{4\varepsilon}\nabla_k\myu{u}^\gamma\big)$ restricted to $\Omega$
is actually supported in $\Omega\setminus\Sigma_{9\varepsilon}$.
In the second one, we use Lemma $\ref{lemma:6.1}$. In the last one, we note that $\myu{u}$ is the extension of $u_0$.

If we set $h_{ik,\varepsilon}^{\alpha\gamma} = V_{i,\varepsilon}^{\alpha\beta}\chi_{k,\varepsilon}^{\beta\gamma}$ with $k=0,\cdots,d$,
then the following estimate immediately holds by Lemma $\ref{lemma:6.1}$,
\begin{equation*}
\bigg|\int_\Omega h_{ik,\varepsilon}^{\alpha\gamma}S_\varepsilon(\psi_{4\varepsilon}\nabla_k\myu{u}^\gamma)\nabla_i\phi_\varepsilon^\alpha dx\bigg|
\leq \|h_{ik,\varepsilon}^{\alpha\gamma}S_\varepsilon(\psi_{4\varepsilon}\nabla_k\myu{u}^\gamma)\|_{L^2(\mathbb{R}^d)}
\|\nabla_i\phi_\varepsilon^\alpha\|_{L^2(\Omega)}
\leq C\| u_0\|_{H^1(\Omega)}\|\nabla\phi_\varepsilon\|_{L^2(\Omega)}.
\end{equation*}
This together with the estimate $\eqref{f:6.12}$ gives
\begin{equation}\label{f:6.18}
 \Big|\int_\Omega \mathcal{I}\cdot \nabla \phi_\varepsilon dx\Big|
 \leq C\|u_0\|_{H^2(\Omega)}\|\nabla\phi_\varepsilon\|_{L^2(\Omega)}
 + C\varepsilon^{-1}\|u_0\|_{H^1(\Omega\setminus\Sigma_{9\varepsilon})}\|\nabla\phi_\varepsilon\|_{L^2(\Omega\setminus\Sigma_{9\varepsilon})}
\end{equation}
by recalling the expression of $\mathcal{I}$ in $\eqref{eq:1}$. By the same token, it is not hard to acquire
\begin{equation}\label{f:6.19}
\begin{aligned}
&\Big|\int_\Omega\mathcal{J}\cdot\nabla\phi_\varepsilon dx\Big|
\leq C\|u_0\|_{H^1(\Omega)}\|\nabla\phi_\varepsilon\|_{L^2(\Omega)},\\
&\Big|\int_\Omega \mathcal{M}\phi_\varepsilon dx\Big|
\leq C\|u_0\|_{H^2(\Omega)}\|\phi_\varepsilon\|_{L^2(\Omega)}
+ C\varepsilon^{-1}\|u_0\|_{H^1(\Omega\setminus\Sigma_{9\varepsilon})}\|\phi_\varepsilon\|_{L^2(\Omega\setminus\Sigma_{9\varepsilon})},\\
&\Big|\int_\Omega \mathcal{N}\phi_\varepsilon dx\Big|
\leq C\|u_0\|_{H^1(\Omega)}\|\phi_\varepsilon\|_{L^2(\Omega)}.
\end{aligned}
\end{equation}

The rest thing is to study the term of $\int_\Omega \mathcal{K}\cdot \nabla \phi_\varepsilon dx$.  In view of $\eqref{f:4.1}$, we have
\begin{equation*}
\int_\Omega \mathcal{K}\cdot \nabla \phi_\varepsilon dx = R(\phi_\varepsilon) = R_1(\phi_\varepsilon) - R_2(\phi_\varepsilon).
\end{equation*}
It follows from the estimate $\eqref{f:4.2}$ that $R_1(\phi_\varepsilon)=0$,
because $S_\varepsilon(\psi_{4\varepsilon}\myu{u})$ and
$S_\varepsilon(\psi_{4\varepsilon}\nabla\myu{u})$ are supported in $\Sigma_{3\varepsilon}$.
The calculations of estimating $R_2(\phi_\varepsilon)$
are quite similar to $\eqref{f:6.13}$ and $\eqref{f:6.12}$, and therefore some explanations are omitted.
\begin{equation*}
\begin{aligned}
R_2(\phi_\varepsilon) &= \varepsilon \int_\Omega E_{jik,\varepsilon}^{\alpha\gamma}\nabla_j\big\{
S_\varepsilon(\psi_{4\varepsilon}\nabla_k\myu{u}^\gamma)\big\}\nabla_i\phi_\varepsilon^\alpha dx\\
&\leq \varepsilon \|E_{jik,\varepsilon}^{\alpha\gamma}S_\varepsilon(\nabla_j\psi_{4\varepsilon}\nabla_k\myu{u}^\gamma)\|_{L^2(\mathbb{R}^d)}\|\nabla_i\phi_\varepsilon^\alpha\|_{L^2(\Omega\setminus\Sigma_{9\varepsilon})}
+\varepsilon\|E_{jik,\varepsilon}^{\alpha\gamma}S_\varepsilon(\psi_{4\varepsilon}\nabla_{jk}^2\myu{u}^\gamma)\|_{L^2(\mathbb{R}^d)}
\|\nabla_i\phi_\varepsilon^\alpha\|_{L^2(\Omega)} \\
&\leq C\|u_0\|_{H^1(\Omega\setminus\Sigma_{9\varepsilon})}\|\nabla\phi_\varepsilon\|_{L^2(\Omega\setminus\Sigma_{9\varepsilon})}
+ C\varepsilon\|u_0\|_{H^2(\Omega)}\|\nabla\phi_\varepsilon\|_{L^2(\Omega)}.
\end{aligned}
\end{equation*}
This implies
\begin{equation}\label{f:6.20}
\Big|\int_\Omega \mathcal{K}\cdot\nabla\phi_\varepsilon dx \Big|
\leq  C\|u_0\|_{H^1(\Omega\setminus\Sigma_{9\varepsilon})}\|\nabla\phi_\varepsilon\|_{L^2(\Omega\setminus\Sigma_{9\varepsilon})}
+ C\varepsilon\|u_0\|_{H^2(\Omega)}\|\nabla\phi_\varepsilon\|_{L^2(\Omega)}.
\end{equation}

Inserting the estimates $\eqref{f:6.15}$, $\eqref{f:6.16}$, $\eqref{f:6.17}$, $\eqref{f:6.18}$, $\eqref{f:6.19}$ and $\eqref{f:6.20}$ into
$\eqref{f:6.14}$, we derive
\begin{equation}\label{f:6.21}
\begin{aligned}
\Big|\int_{\Omega}w_\varepsilon\Phi dx\Big|
\leq C\Big\{\|u_0\|_{H^1(\Omega\setminus\Sigma_{9\varepsilon})}\|\phi_\varepsilon\|_{H^1(\Omega\setminus\Sigma_{9\varepsilon})}
+ \varepsilon\|u_0\|_{H^2(\Omega)}\|\phi_\varepsilon\|_{H^1(\Omega)}\Big\},
\end{aligned}
\end{equation}

To improve the order of the convergence rate,
the next task is to replace $\phi_\varepsilon$ into the first order corrector $\xi_\varepsilon$
in the first term of the right-hand side of $\eqref{f:6.21}$,
where
$\xi_\varepsilon = \phi_\varepsilon - \phi_0 -\varepsilon\chi_{0,\varepsilon}^*S_\varepsilon(\psi_{10\varepsilon}\phi_0)
-\varepsilon\chi_{j,\varepsilon}^*S_\varepsilon(\psi_{10\varepsilon}\nabla_j \phi_0)$.
Observing that $S_\varepsilon(\psi_{10\varepsilon}\phi_0)$ and $S_\varepsilon(\psi_{10\varepsilon}\nabla_j \phi_0)$ are supported in
$\Sigma_{9\varepsilon}$, we arrive at
\begin{equation}\label{f:6.22}
\|\phi_\varepsilon\|_{H^1(\Omega\setminus\Sigma_{9\varepsilon})}
\leq \|\xi_\varepsilon\|_{H^1(\Omega\setminus\Sigma_{9\varepsilon})}
+ \|\phi_0\|_{H^1(\Omega\setminus\Sigma_{9\varepsilon})}
\end{equation}
Thus we plug $\eqref{f:6.22}$ back into $\eqref{f:6.21}$ and obtain
\begin{equation}
\begin{aligned}
\Big|\int_{\Omega}w_\varepsilon\Phi dx\Big|
&\leq C\Big\{\|u_0\|_{H^1(\Omega\setminus\Sigma_{9\varepsilon})}\|\phi_0\|_{H^1(\Omega\setminus\Sigma_{9\varepsilon})}
+ \|u_0\|_{H^1(\Omega\setminus\Sigma_{9\varepsilon})}\|\xi_\varepsilon\|_{H^1(\Omega)}
+ \varepsilon\|u_0\|_{H^2(\Omega)}\|\phi_\varepsilon\|_{H^1(\Omega)}\Big\} \\
&\leq C\Big\{\|u_0\|_{H^1(\Omega\setminus\Sigma_{9\varepsilon})}\|\phi_0\|_{H^1(\Omega\setminus\Sigma_{9\varepsilon})}
+ \varepsilon^{\frac{1}{2}}\|u_0\|_{H^1(\Omega\setminus\Sigma_{9\varepsilon})}\|\Phi\|_{L^2(\Omega)}
+ \varepsilon\|u_0\|_{H^2(\Omega)}\|\Phi\|_{L^2(\Omega)}\Big\},
\end{aligned}
\end{equation}
where we use the estimate $\eqref{f:6.23}$ and Lemma $\ref{lemma:2.1}$ in the last inequality.

The problem reduces to estimate the layer quantity. On account of the estimate $\eqref{f:6.24}$, in fact we have
\begin{equation*}
\|u_0\|_{H^1(\Omega\setminus\Sigma_{9\varepsilon})}\leq C\varepsilon^{\frac{1}{2}}\|u_0\|_{H^2(\Omega)},
\qquad
\|\phi_0\|_{H^1(\Omega\setminus\Sigma_{9\varepsilon})}\leq C\varepsilon^{\frac{1}{2}}\|\phi_0\|_{H^2(\Omega)}.
\end{equation*}
Thus it follows that
\begin{equation*}
\Big|\int_{\Omega}w_\varepsilon\Phi dx\Big|
\leq C\varepsilon\|u_0\|_{H^2(\Omega)}\|\Phi\|_{L^2(\Omega)},
\end{equation*}
and we consequently obtain the desired estimate $\eqref{pri:6.7}$ by using duality and $H^{2}$ estimates
(where the assumption of $\partial\Omega\in C^{1,1}$ has been used). The proof is complete.
\qed

\end{pf}

\begin{lemma}\label{lemma:4.2}
Let $\Omega$ be a bounded Lipschitz domain.
Suppose that the coefficients of $\myl{L}{\varepsilon}$ satisfy $\eqref{a:1}$, $\eqref{a:2}$ and $\eqref{a:3}$,
and $A$ additionally satisfies $A=A^*$. Let $u_\varepsilon,u_0$ be the weak solutions to $\eqref{pde:1.5}$
with $F\in L^p(\Omega;\mathbb{R}^m)$ and $g\in L^2(\partial\Omega;\mathbb{R}^m)$, where $p=\frac{2d}{d+1}$. Then
\begin{equation}\label{pri:4.6}
\big\|u_\varepsilon - u_0 -\varepsilon\chi_{0,\varepsilon}S^2_\varepsilon(\psi_{2\varepsilon}u_0)
-\varepsilon\chi_{k,\varepsilon}S^2_\varepsilon(\psi_{2\varepsilon}\nabla_k u_0)\big\|_{H^1(\Omega)}
\leq C\varepsilon^{\frac{1}{2}}\big\{\|F\|_{L^p(\Omega)}+\|g\|_{L^2(\partial\Omega)}\big\},
\end{equation}
where $C$ depends only on $\mu,\kappa,m,d$ and $\Omega$.
\end{lemma}

\begin{pf}
We note that $\|w_\varepsilon\|_{H^1(\Omega)}$  is exactly the left-hand side of $\eqref{pri:4.6}$
by setting $\varphi_0 = \mys{S}^2(\psi_{2\varepsilon}u_0)$ and $\varphi_k = \mys{S}^2(\psi_{2\varepsilon}\nabla_k u_0)$ in $\eqref{f:4.13}$.
Then it follows from $\eqref{pri:4.5}$ that
\begin{equation}\label{f:5.25}
\begin{aligned}
\|w_{\varepsilon}\|_{H^1(\Omega)}
&\leq C\Big\{
\|\nabla u_0 - \mys{S}^2(\psi_{2\varepsilon}\nabla u_0)\|_{L^2(\Omega)}
+\|u_0 - \mys{S}^2(\psi_{2\varepsilon}u_0)\|_{L^2(\Omega)}
+\varepsilon\|h_{\varepsilon}\nabla\mys{S}^2(\psi_{2\varepsilon} u_0)\|_{L^2(\Omega)}\\
&+ \varepsilon \|h_{\varepsilon}\nabla\mys{S}^2(\psi_{2\varepsilon}\nabla u_0)\|_{L^2(\Omega)}
+ \varepsilon \|h_{\varepsilon}\mys{S}^2(\psi_{2\varepsilon}u_0)\|_{L^2(\Omega)}
+\varepsilon \|h_{\varepsilon}\mys{S}^2(\psi_{2\varepsilon}\nabla u_0)\|_{L^2(\Omega)}\Big\}.
\end{aligned}
\end{equation}
Before proceeding further, let us do some calculations:
\begin{equation}\label{f:5.26}
\begin{aligned}
\|\nabla u_0 - \mys{S}^2(\psi_{2\varepsilon}\nabla u_0)\|_{L^2(\Omega)}
&\leq\|(1-\psi_{2\varepsilon})\nabla u_0\|_{L^2(\Omega)}
+ \|\psi_{2\varepsilon}\nabla u_0 - \mys{S}(\psi_{2\varepsilon}\nabla u_0)\|_{L^2(\Omega)} \\
& +\big\|\mys{S}\big[\psi_{2\varepsilon}\nabla u_0 - \mys{S}(\psi_{2\varepsilon}\nabla u_0)\big]\|_{L^2(\Omega)} \\
&\leq \|\nabla u_0\|_{L^2(\Omega\setminus\Sigma_{4\varepsilon})}
+ C\|\psi_{2\varepsilon}\nabla u_0 - \mys{S}(\psi_{2\varepsilon}\nabla u_0)\|_{L^2(\Omega)},
\end{aligned}
\end{equation}
\begin{equation}\label{f:5.27}
\begin{aligned}
\|h_{\varepsilon}\nabla\mys{S}^2(\psi_{2\varepsilon}\nabla u_0)\|_{L^2(\Omega)}
&\leq C\|\nabla\mys{S}(\psi_{2\varepsilon}\nabla u_0)\|_{L^2(\Omega)} \\
&\leq C\big\{\varepsilon^{-1}\|\nabla u_0\|_{L^2(\Omega\setminus\Sigma_{4\varepsilon})}
+ \|\mys{S}(\psi_{2\varepsilon}\nabla^2 u_0)\|_{L^2(\Omega)}\big\},
\end{aligned}
\end{equation}
and
\begin{equation}\label{f:5.28}
\|h_{\varepsilon}\mys{S}^2(\psi_{2\varepsilon}\nabla u_0)\|_{L^2(\Omega)}
\leq C\|\mys{S}(\psi_{2\varepsilon}\nabla u_0)\|_{L^2(\Omega)}
\end{equation}
where we mainly use the estimate $\eqref{pri:6.1}$ in the second inequality of $\eqref{f:5.26}$, and in the first inequality of $\eqref{f:5.27}$,
as well as in $\eqref{f:5.28}$. After a similar computation, we have
\begin{equation}\label{f:5.29}
\begin{aligned}
&\|h_{\varepsilon}\mys{S}^2(\psi_{2\varepsilon} u_0)\|_{L^2(\Omega)}
\leq C\|\mys{S}(\psi_{2\varepsilon}u_0)\|_{L^2(\Omega)}, \\
& \| u_0 - \mys{S}^2(\psi_{2\varepsilon} u_0)\|_{L^2(\Omega)}
\leq \|u_0\|_{L^2(\Omega\setminus\Sigma_{4\varepsilon})}
+ C\|\psi_{2\varepsilon} u_0 - \mys{S}(\psi_{2\varepsilon} u_0)\|_{L^2(\Omega)}, \\
&\|h_{\varepsilon}\nabla\mys{S}^2(\psi_{2\varepsilon} u_0)\|_{L^2(\Omega)}
\leq C\Big\{\frac{1}{\varepsilon}\|u_0\|_{L^2(\Omega\setminus\Sigma_{4\varepsilon})}
+ \|\mys{S}(\psi_{2\varepsilon}\nabla u_0)\|_{L^2(\Omega)}\Big\}. \\
\end{aligned}
\end{equation}
By substituting $\eqref{f:5.26}$, $\eqref{f:5.27}$, $\eqref{f:5.28}$ and $\eqref{f:5.29}$ into $\eqref{f:5.25}$, we find
\begin{equation}\label{f:4.21}
\begin{aligned}
\|w_\varepsilon\|_{H^1(\Omega)}& \leq C\Big\{
\|u_0\|_{H^1(\Omega\setminus\Sigma_{4\varepsilon})}
+\|\psi_{2\varepsilon}\nabla u_0 - \mys{S}(\psi_{2\varepsilon}\nabla u_0)\|_{L^2(\Omega)}
+\|\psi_{2\varepsilon}u_0 - \mys{S}(\psi_{2\varepsilon}u_0)\|_{L^2(\Omega)}\Big\} \\
& + C\varepsilon\Big\{\|\mys{S}(\psi_{2\varepsilon}\nabla^2 u_0)\|_{L^2(\Omega)}
+ \|\mys{S}(\psi_{2\varepsilon}\nabla u_0)\|_{L^2(\Omega)}
+\|\mys{S}(\psi_{2\varepsilon}u_0)\|_{L^2(\Omega)} \Big\}.
\end{aligned}
\end{equation}

We now handle $\|u_0\|_{H^1(\Omega\setminus\Sigma_{4\varepsilon})}$ in the right hand side of $\eqref{f:4.21}$.
First of all, we rewrite ($\text{H}_0$) of $\eqref{pde:1.5}$ as
\begin{equation}
 L_0(u_0) = F + (\widehat{V} -\widehat{B})\nabla u_0 - (\widehat{c}+\lambda I)u_0
 \quad \text{in}~\Omega,
 \quad\qquad
 \partial u_0/\partial\nu_0 = g - n\cdot\widehat{V}u_0
 \quad \text{on}~\partial\Omega,
\end{equation}
where $L_0 = -\text{div}(\widehat{A}\nabla)$ and $\partial/\partial\nu_0 = n\cdot\widehat{A}\nabla$.
Let $u_0 = v + \rho$,
and $v,\rho$ satisfy
\begin{equation}\label{pde:4.6}
(\text{D}_1)~L_0(v) = \tilde{F} \quad\text{in}~\mathbb{R}^d,
\qquad\qquad (\text{D}_2)~\left\{\begin{aligned}
 L_0(\rho) &= 0 &\quad\text{in}~~~\Omega,\\
 \frac{\partial\rho}{\partial\nu_0} &= g - n\cdot\widehat{V}u_0 - \frac{\partial v}{\partial\nu_0} &\quad\text{on}~\partial\Omega,
\end{aligned}\right.
\end{equation}
respectively. Note that $\tilde{F} = F + (\widehat{V} -\widehat{B})\nabla u_0 - (\widehat{c}+\lambda I)u_0$ in $\Omega$,
and $\tilde{F}=0$ in $\mathbb{R}^d\setminus\Omega$.

For $(\text{D}_1)$, let $\Gamma_0$ denote the fundamental solution of $L_0$, and then we have
$v=\Gamma_0\ast\tilde{F}$ in $\mathbb{R}^d$. Moreover,
we have the estimate
\begin{equation*}
|\nabla v(x)|\leq C\int_{\mathbb{R}^d}\frac{|\tilde{F}(y)|}{|x-y|^{d-1}}dy,
\end{equation*}
and by Hardy-Littlewood-Sobolev inequality on fractional integration (see \cite[Theorem 7.25]{MGLM}),
\begin{equation}\label{f:6.26}
\|\nabla v\|_{L^2(\mathbb{R}^d)}\leq C\|\tilde{F}\|_{L^{\frac{2d}{d+2}}(\mathbb{R}^d)}
\leq C\big\{\|F\|_{L^{\frac{2d}{d+2}}(\Omega)}+\|u_0\|_{H^1(\Omega)}\big\}
\leq C\big\{\|F\|_{L^{\frac{2d}{d+2}}(\Omega)}+\|g\|_{L^2(\partial\Omega)}\big\}.
\end{equation}
Also, it follows from the Calder\'on-Zygmund theorem (see \cite[Theorem 7.22]{MGLM}) that
\begin{equation}\label{f:4.22}
\|\nabla^2 v\|_{L^{\frac{2d}{d+1}}(\mathbb{R}^d)}
\leq C\|\tilde{F}\|_{L^{\frac{2d}{d+1}}(\mathbb{R}^d)}
\leq C\big\{\|F\|_{L^{\frac{2d}{d+1}}(\Omega)}+\|u_0\|_{H^1(\Omega)}\big\}
\leq C\big\{\|F\|_{L^{\frac{2d}{d+1}}(\Omega)}+\|g\|_{L^2(\partial\Omega)}\big\}.
\end{equation}
Note that we use $\eqref{pri:2.4}$ in the last inequalities of $\eqref{f:6.26}$ and $\eqref{f:4.22}$ with the fact of
$L^2(\partial\Omega;\mathbb{R}^m)\subset B^{-1/2,2}(\partial\Omega;\mathbb{R}^m)$.
Then due to Sobolev inequality, we have
\begin{equation}\label{f:6.27}
\begin{aligned}
 &\|v\|_{L^{\frac{2d}{d-2}}(\mathbb{R}^d)}\leq C\|\nabla v\|_{L^2(\mathbb{R}^d)}
 \leq C \big\{\|F\|_{L^{\frac{2d}{d+2}}(\Omega)} + \|g\|_{L^2(\partial\Omega)} \big\}, \\
 &\|\nabla v\|_{L^{\frac{2d}{d-1}}(\mathbb{R}^d)}
 \leq C \|\nabla^2 v\|_{L^{\frac{2d}{d+1}}(\mathbb{R}^d)}
 \leq C \big\{\|F\|_{L^{\frac{2d}{d+1}}(\Omega)}+\|g\|_{L^2(\partial\Omega)}\big\}.
\end{aligned}
\end{equation}

Hence, by suitable modification to the proof of the estimate $\eqref{f:6.25}$, we acquire
\begin{equation}\label{f:4.23}
\begin{aligned}
c\int_{\partial\Omega} \big(|\nabla v|^2 + |v|^2\big) dS
&\leq \int_{\partial\Omega}<\varrho,n> (|\nabla v|^2 + |v|^2)dS\\
&= \int_\Omega \text{div}(\varrho)\big(|\nabla v|^2+|v|^2\big) dx
+ 2\int_\Omega \varrho\big(\nabla^2 v \nabla v + \nabla v v\big)dx \\
&\leq C\big\{\|v\|_{L^{\frac{2d}{d-2}}(\Omega)}^2 +\|\nabla v\|_{L^2(\Omega)}^2
+ \|\nabla v\|_{L^{\frac{2d}{d-1}}(\Omega)}\|\nabla^2 v\|_{L^{\frac{2d}{d+1}}(\Omega)}\big\} \\
&\leq C\big\{\|F\|_{L^{\frac{2d}{d+1}}(\Omega)}^2+\|g\|_{L^2(\partial\Omega)}^2\big\}.
\end{aligned}
\end{equation}
where we use H\"older's inequality and Young's inequality in the second inequality,
and the estimates $\eqref{f:6.26}$, $\eqref{f:4.22}$ and $\eqref{f:6.27}$ in the last one.
As explained for the estimate $\eqref{f:6.24}$, it is not hard to see
\begin{equation}\label{f:4.24}
\Big\{\frac{1}{\varepsilon}\int_{\Omega\setminus\Sigma_{4\varepsilon}}\big(|\nabla v|^2 + |v|^2\big) dx\Big\}^{\frac{1}{2}}
\leq C\big\{\|F\|_{L^{\frac{2d}{d+1}}(\Omega)}+\|g\|_{L^2(\partial\Omega)}\big\}.
\end{equation}

We now turn to study $(\text{D}_2)$. From the Rellich identity, the estimates $\eqref{pri:2.4}$ and $\eqref{f:4.23}$
(see Remark $\ref{re:4.2}$), it follows that
\begin{equation}\label{f:4.43}
\begin{aligned}
 \|(\nabla\rho)^*\|_{L^2(\partial\Omega)}
 &\leq C\big\{\|g\|_{L^2(\partial\Omega)}+ \|F\|_{L^{\frac{2d}{d+1}}(\Omega)}\big\}.
\end{aligned}
\end{equation}
According to Lemma $\ref{lemma:2.6}$, we have
\begin{equation*}
\|\mathcal{M}(\rho)\|_{L^2(\partial\Omega)} \leq C\|\rho\|_{H^1(\Omega\setminus\Sigma_{c_0})}
\leq C\|\rho\|_{H^1(\Omega)}
\leq C\big\{\|g\|_{L^2(\partial\Omega)}+ \|F\|_{L^{\frac{2d}{d+2}}(\Omega)}\big\},
\end{equation*}
where we use $\eqref{pri:2.4}$ in the last inequality, and $C$ depends on $d,m,c_0$ and the character of $\Omega$. This together with
$\eqref{f:4.43}$ leads to
\begin{equation}\label{f:4.25}
\begin{aligned}
\Big\{\frac{1}{\varepsilon}\int_{\Omega\setminus\Sigma_{4\varepsilon}}\big(|\nabla \rho|^2 + |\rho|^2\big) dx\Big\}^{\frac{1}{2}}
& \leq C\big\{\|(\nabla\rho)^*\|_{L^2(\partial\Omega)}+\|\mathcal{M}(\rho)\|_{L^2(\partial\Omega)}\big\} \\
&\leq C\big\{\|g\|_{L^2(\partial\Omega)}+\|F\|_{L^{\frac{2d}{d+1}}(\Omega)}\big\}.
\end{aligned}
\end{equation}
where we use the estimate $\eqref{pri:6.5}$ (for $p=2$ and $r=4\varepsilon$) in the first inequality.

Hence combining $\eqref{f:4.24}$ and $\eqref{f:4.25}$, we have
\begin{equation}\label{f:4.29}
\Big\{\int_{\Omega\setminus\Sigma_{4\varepsilon}}\big(|\nabla u_0|^2 + |u_0|^2\big) dx\Big\}^{\frac{1}{2}}
\leq C\varepsilon^{\frac{1}{2}}\big\{\|g\|_{L^2(\partial\Omega)}+\|F\|_{L^{\frac{2d}{d+1}}(\Omega)}\big\}.
\end{equation}

We now estimate $\|\psi_{2\varepsilon}\nabla u_0 - \mys{S}(\psi_{2\varepsilon}\nabla u_0)\|_{L^2(\Omega)}$ and
$\|\psi_{2\varepsilon}u_0 - \mys{S}(\psi_{2\varepsilon}u_0)\|_{L^2(\Omega)}$ in the right hand side of $\eqref{f:4.21}$.
Since $\psi_{2\varepsilon}\nabla u_0 - \mys{S}(\psi_{2\varepsilon}\nabla u_0)$ is supported in $\Sigma_\varepsilon$, it is
equivalent to estimating $\|\psi_{2\varepsilon}\nabla u_0 - \mys{S}(\psi_{2\varepsilon}\nabla u_0)\|_{L^2(\Sigma_\varepsilon)}$.
It immediately follows from $\eqref{pri:6.2}$, $\eqref{pri:6.3}$ and $\eqref{f:4.24}$ that
\begin{equation}\label{f:4.27}
\begin{aligned}
&\|\psi_{2\varepsilon}\nabla u_0 - \mys{S}(\psi_{2\varepsilon}\nabla u_0)\|_{L^2(\Sigma_\varepsilon)} \\
&\leq \|\nabla v - \mys{S}(\nabla v)\|_{L^2(\Sigma_\varepsilon)}
+\|(\psi_{2\varepsilon}-1)\nabla v\|_{L^2(\Sigma_\varepsilon)}
 +\|S_\varepsilon\big((\psi_{2\varepsilon}-1)\nabla v\big)\|_{L^2(\Sigma_\varepsilon)}
+ \|\psi_{2\varepsilon}\nabla\rho - \mys{S}(\psi_{2\varepsilon}\nabla\rho)\|_{L^2(\Sigma_\varepsilon)} \\
&\leq \|\nabla v - \mys{S}(\nabla v)\|_{L^2(\mathbb{R}^d)}
+\|(\psi_{2\varepsilon}-1)\nabla v\|_{L^2(\Omega)}
 +\|S_\varepsilon\big((\psi_{2\varepsilon}-1)\nabla v\big)\|_{L^2(\Sigma_\varepsilon)}
+ \|\psi_{2\varepsilon}\nabla\rho - \mys{S}(\psi_{2\varepsilon}\nabla\rho)\|_{L^2(\Omega)} \\
& \leq  C\varepsilon^{\frac{1}{2}}\|\nabla^2 v\|_{L^{\frac{2d}{d+1}}(\mathbb{R}^d)}
+ C\|\nabla v\|_{L^2(\Omega\setminus\Sigma_{4\varepsilon})}
+ C\varepsilon\|\nabla(\nabla\rho\psi_{2\varepsilon})\|_{L^2(\Omega)} \\
& \leq C\varepsilon^{\frac{1}{2}}\big\{\|F\|_{L^{\frac{2d}{d+1}}(\Omega)} + \|g\|_{L^2(\partial\Omega)} \big\}
+ C\big\{\varepsilon\|\nabla^2\rho\|_{L^2(\Sigma_{2\varepsilon})}
+ \|\nabla\rho\|_{L^2(\Omega\setminus\Sigma_{4\varepsilon})}\big\}.
\end{aligned}
\end{equation}
Note that $\rho$ satisfies $(\text{D}_2)$, from the interior estimate for $L_{0}$ (see Remark $\ref{re:5.3}$), we have
\begin{equation*}
 |\nabla^2\rho(x)|\leq \frac{C}{\delta(x)}\Big(\dashint_{B(x,\delta(x)/8)}|\nabla\rho(y)|^2 dy\Big)^{\frac{1}{2}}.
\end{equation*}
This gives
\begin{equation}\label{f:4.42}
 \int_{\Sigma_{2\varepsilon}}|\nabla^2\rho|^2 dx
 \leq C\int_{\Sigma_{2\varepsilon}}\dashint_{B(x,\delta(x)/8)}\frac{|\nabla\rho(y)|^2}{[\delta(x)]^2}dydx
 \leq C\int_\varepsilon^\infty\int_{S_t}\frac{|(\nabla\rho)^*(x^\prime)|^2}{t^2}dS_tdt.
\end{equation}
Note that $\delta(x)\approx t$, and $x^\prime\in S_{\delta(x)/4}$ such that $|\nabla\rho(y)|\leq(\nabla\rho)^*(x^\prime)$
for any $y\in B(x,\delta(x)/8)$. By using the observation that
$\|(\nabla\rho)^*\|_{L^2(S_t)}\leq\|(\nabla\rho)^*\|_{L^2(\partial\Omega)}$ holds for all $t\in[0,\infty)$,
it follows from $\eqref{f:4.43}$ and $\eqref{f:4.42}$ that
\begin{equation}\label{f:4.28}
 \|\nabla^2\rho\|_{L^2(\Sigma_{2\varepsilon})}
 \leq C\varepsilon^{-\frac{1}{2}}\|(\nabla\rho)^*\|_{L^2(\partial\Omega)}
 \leq C\varepsilon^{-\frac{1}{2}}\big\{\|g\|_{L^2(\partial\Omega)} + \|F\|_{L^{\frac{2d}{d+1}}(\Omega)}\big\}.
\end{equation}

Now, collecting $\eqref{f:4.25},\eqref{f:4.27}$ and $\eqref{f:4.28}$ leads to
\begin{equation}\label{f:4.30}
\|\psi_{2\varepsilon}\nabla u_0 - \mys{S}(\psi_{2\varepsilon}\nabla u_0)\|_{L^2(\Omega)}
\leq C\varepsilon^{\frac{1}{2}}\big\{\|g\|_{L^2(\partial\Omega)}+\|F\|_{L^{\frac{2d}{d+1}}(\Omega)}\big\}.
\end{equation}

By the same procedure as above, we have
\begin{equation}\label{f:4.31}
\begin{aligned}
\|\psi_{2\varepsilon} u_0 - \mys{S}(\psi_{2\varepsilon} u_0)\|_{L^2(\Omega)}
& \leq C\varepsilon\|\nabla(\psi_{2\varepsilon}u_0)\|_{L^2(\Omega)}
\leq C\|u_0\|_{L^2(\Omega\setminus\Sigma_{4\varepsilon})} + C\varepsilon\|\nabla u_0\|_{L^2(\Omega)}\\
&\leq  C\varepsilon^{\frac{1}{2}}\big\{\|g\|_{L^2(\partial\Omega)}+\|F\|_{L^{\frac{2d}{d+1}}(\Omega)}\big\},
\end{aligned}
\end{equation}
where we use $\eqref{pri:6.2}$ in the first inequality, and $\eqref{pri:2.4}$ and $\eqref{f:4.29}$ in the last one.

To accomplish the proof, we still need the following estimates:
\begin{equation}\label{f:4.32}
\begin{aligned}
\|S_\varepsilon(\psi_{2\varepsilon}\nabla^2u_0)\|_{L^2(\Omega)}
&\leq \|S_\varepsilon(\psi_{2\varepsilon}\nabla^2v)\|_{L^2(\Omega)}
+ \|S_\varepsilon(\psi_{2\varepsilon}\nabla^2\rho)\|_{L^2(\Omega)} \\
& \leq C\varepsilon^{-\frac{1}{2}}\|\psi_{2\varepsilon}\nabla^2 v\|_{L^{\frac{2d}{d+1}}(\mathbb{R}^d)} + C\|\nabla^2\rho\|_{L^2(\Sigma_{2\varepsilon})}\\
&\leq C\varepsilon^{-\frac{1}{2}}\big\{\|g\|_{L^2(\partial\Omega)}+\|F\|_{L^{\frac{2d}{d+1}}(\Omega)}\big\},
\end{aligned}
\end{equation}
where we use $\eqref{pri:6.3}$ in the second inequality, and $\eqref{f:4.22}$, $\eqref{f:4.28}$ in the last one. Also,
\begin{equation}\label{f:4.33}
\begin{aligned}
\|S_\varepsilon(\psi_{2\varepsilon}\nabla u_0)\|_{L^2(\Omega)}
+ \|S_\varepsilon(\psi_{2\varepsilon}u_0)\|_{L^2(\Omega)}
\leq C\|u_0\|_{H^1(\Omega)}\leq C\big\{\|g\|_{L^2(\partial\Omega)}+\|F\|_{L^{\frac{2d}{d+2}}(\Omega)}\big\}.
\end{aligned}
\end{equation}

Finally, combining $\eqref{f:4.21}$, $\eqref{f:4.29}$, $\eqref{f:4.30}$, $\eqref{f:4.31}$, $\eqref{f:4.32}$ and $\eqref{f:4.33}$,
we obtain the estimate $\eqref{pri:4.6}$ and the proof is done.
\qed
\end{pf}

\begin{remark}\label{re:4.2}
\emph{To see the estimate $\eqref{f:4.43}$, we first show the Rellich identity for the equation: $L_0(u) = \tilde{F}$ in $\Omega$
and $\partial u/\partial\nu_0 = \tilde{g}$ on $\partial\Omega$,
where $\tilde{F}\in L^2(\Omega;\mathbb{R}^m)$ and $\tilde{g}\in L^2(\partial\Omega;\mathbb{R}^m)$
with $\int_\Omega \tilde{F}^\alpha dx + \int_{\partial\Omega} \tilde{g}^{\alpha} dS = 0$.
\begin{equation}
\begin{aligned}
\int_{\partial\Omega} \varrho_k n_k \hat{a}_{ij}^{\alpha\beta}\frac{\partial u^\beta}{\partial x_j}\frac{\partial u^\alpha}{\partial x_i} dS
&= \int_{\Omega} \frac{\partial \varrho_k}{\partial x_k} \hat{a}_{ij}^{\alpha\beta}\frac{\partial u^\beta}{\partial x_j}\frac{\partial u^\alpha}{\partial x_i}
+ 2\int_\Omega\frac{\partial u^\alpha}{\partial x_k}\tilde{F}^\alpha \varrho_k dx \\
&- 2\int_\Omega\frac{\partial u^\alpha}{\partial x_k} \hat{a}_{ij}^{\alpha\beta}\frac{\partial u^\beta}{\partial x_j}\frac{\partial \varrho_k}{\partial x_i} dx
+ 2\int_{\partial\Omega} \tilde{g}^{\alpha}\frac{\partial u^\alpha}{\partial x_k} \varrho_k dS,
\end{aligned}
\end{equation}
where $\varrho$ is a $C^1_0(\mathbb{R}^d;\mathbb{R}^d)$ vector field similarly defined as in $\eqref{f:6.25}$,
and we use the assumption of $A=A^*$.
Coupled with the $H^1$ estimate (see $\eqref{pri:2.10}$ for $p=2$ and $f=0$), it is not hard to see that
\begin{equation}\label{pri:4.10}
 \|\nabla u\|_{L^2(\partial\Omega)}
 \leq C\big\{\|\tilde{F}\|_{L^2(\Omega)}+ \|\tilde{g}\|_{L^2(\partial\Omega)}\big\}.
\end{equation}
Then, we apply the estimate $\eqref{pri:4.10}$ to the solution $\rho$ to $(\text{D}_2)$ in $\eqref{pde:4.6}$
(where $\tilde{F} = 0$ and $\tilde{g} = g -n\cdot\widehat{V}u_0 - \partial v/\partial v$), and acquire
\begin{equation}\label{pri:5.9}
\begin{aligned}
\|\nabla \rho\|_{L^2(\partial\Omega)}
& \leq C\big\{\|g\|_{L^2(\partial\Omega)}+\|u_0\|_{L^2(\partial\Omega)}
+ \|\nabla v\|_{L^2(\partial\Omega)}\big\} \\
& \leq C\big\{\|g\|_{L^2(\partial\Omega)}
+ \|u_0\|_{H^1(\Omega)}+ \|\nabla v\|_{L^2(\partial\Omega)} \big\}
&\leq C\big\{\|g\|_{L^2(\partial\Omega)} + \|F\|_{\frac{2d}{d+1}(\Omega)}\big\},
\end{aligned}
\end{equation}
where we use the trace theorem (see \cite[Lemma 5.17]{GML}) in second inequality,
and the estimate $\eqref{pri:2.4}$ for $u_0$ 
and the estimate $\eqref{f:4.23}$ are employed in the last one.}

\emph{Moreover, according to the solvability of $L^2$-Dirichlet problems (see \cite[Theorem 1.3]{SZW11}), it is not hard to see
$\|(\nabla\rho)^*\|_{L^2(\partial\Omega)} \leq C\|\nabla\rho\|_{L^2(\partial\Omega)}$.
This together with $\eqref{pri:5.9}$ leads to
\begin{equation*}
\|(\nabla\rho)^*\|_{L^2(\partial\Omega)}
\leq C\big\{\|g\|_{L^2(\partial\Omega)} + \|F\|_{\frac{2d}{d+1}(\Omega)}\big\}.
\end{equation*}}
\end{remark}

\begin{remark}\label{re:5.3}
\emph{Let $L_0 = -\text{div}(\widehat{A}\nabla)$,
and $u\in H^1_{loc}(\Omega;\mathbb{R}^m)$ be a weak solution to $L_0(u)=0$ in $\Omega$. For any $B(P,R)\subset 4B\subset\Omega$,
we may assume $P=0$ and $R=1$ from the translation and rescaling arguments.
Then due to the interior $H^k$ regularity theory (see \cite[Theorem 4.11]{MGLM}), we have
$\|u\|_{H^k(B)}\leq C(k,m,d,\mu)\|u\|_{L^2(2B)}$, where
$H^k(\Omega;\mathbb{R}^m)= W^{k,2}(\Omega;\mathbb{R}^m)$.
Moreover, it follows from the Sobolev imbedding theorem (for $2k>d$) that
$|u(0)|\leq \|u\|_{L^\infty(B)}\leq C\|u\|_{H^k(B)}\leq C\|u\|_{L^2(2B)}$. Note that $v=\nabla^2_{ij}u$ also satisfy
$L_0(v) = 0$ in $\Omega$. In fact, we thus have $|\nabla^2 u(0)|\leq C\|\nabla^2 u\|_{L^2(2B)}$. By Cacciopolli's inequality $\eqref{pri:2.13}$
(or see \cite[Theorem 4.1]{MGLM}), we acquire $|\nabla^2 u(0)|\leq C\|\nabla u\|_{L^2(4B)}$. Let $\tilde{u}(x) = u(Rx)$ be the scaled function,
a routine computation gives
\begin{equation*}
\big|\nabla^2 u(0)\big|\leq \frac{C}{R}\Big(\dashint_{B(0,4R)}|\nabla u|^2dy\Big)^{1/2}.
\end{equation*}}
\end{remark}

\begin{remark}
\emph{Let $u_0 = v+\rho$ and $v,\rho$ be given in $\eqref{pde:4.6}$.
As shown in the proof of Lemma $\ref{lemma:4.2}$, we can also prove the following results:
\begin{equation}\label{f:6.28}
\begin{aligned}
 \|\nabla^2 u_0\|_{L^2(\Sigma_{N\varepsilon})}
 &\leq \|\nabla^2 v\|_{L^2(\mathbb{R}^d)} + \|\nabla^2\rho\|_{L^2(\Sigma_{N\varepsilon})} \\
 &\leq C\big\{\|F\|_{L^2(\Omega)}+\|g\|_{L^2(\partial\Omega)}\big\} + C\varepsilon^{-\frac{1}{2}}\|(\nabla\rho)^*\|_{L^2(\partial\Omega)}
 \leq C\varepsilon^{-\frac{1}{2}}\big\{\|F\|_{L^2(\Omega)}+\|g\|_{L^2(\partial\Omega)}\big\},
\end{aligned}
\end{equation}
where $N$ is a positive integer.
Note that we use the Calder\'on-Zygmund theorem (see \cite[Theorem 7.22]{MGLM}) for $v$
and the estimate $\eqref{f:4.28}$ in the second inequality, and H\"older's inequality in the last one.
Also,
\begin{equation}\label{f:6.29}
\|u_0\|_{H^1(\Omega\setminus\Sigma_{N\varepsilon})}
\leq C\varepsilon^{\frac{1}{2}}\big\{\|F\|_{L^{\frac{2d}{d+1}}(\Omega)}+\|g\|_{L^2(\partial\Omega)}\big\},
\end{equation}
which is the counterpart of the estimate $\eqref{f:4.29}$.
Note that if we combine $\eqref{f:6.11}$, $\eqref{f:6.28}$ and $\eqref{f:6.29}$ by setting $N=4$,
then we have
\begin{equation}\label{f:6.30}
\big\|u_\varepsilon - u_0
-\varepsilon\chi_{0,\varepsilon}S_\varepsilon(\psi_{4\varepsilon}u_0)
-\varepsilon\chi_{k,\varepsilon}S_\varepsilon(\psi_{4\varepsilon}\nabla_k u_0)\big\|_{H^1(\Omega)}
\leq C\varepsilon^{\frac{1}{2}}\big\{\|F\|_{L^2(\Omega)}+\|g\|_{L^{2}(\partial\Omega)}\big\}
\end{equation}
under the assumption that $\Omega$ is just a bounded Lipshitz domain.
The estimate $\eqref{f:6.30}$ is similar to $\eqref{pri:4.6}$ in Lemma $\ref{lemma:4.2}$, but
the estimate $\eqref{pri:4.6}$ is sharp in the sense of the integrability of the given data $F$, and it is also
the reason why we mollify $\psi_{2\varepsilon}(\nabla u_0)$ and $\psi_{2\varepsilon}u_0)$ twice in Lemma $\ref{lemma:4.2}$. }
\end{remark}

\begin{thm}\label{thm:4.1}
Let $\Omega$ be a bounded Lipschitz domain. Suppose that the coefficients of $\myl{L}{\varepsilon}$ satisfy
$\eqref{a:1}-\eqref{a:3}$, and additional condition $A=A^*$. Assume $u_\varepsilon\in H^1(\Omega;\mathbb{R}^m)$ is the weak solution to
the Neumann problem $\myl{L}{\varepsilon}(u_\varepsilon) = F$ in $\Omega$, and $\myl{B}{\varepsilon}(u_\varepsilon) = g$ on $\partial\Omega$,
where $F\in L^p(\Omega;\mathbb{R}^m)$ for $p=\frac{2d}{d+1}$ and $g\in L^2(\partial\Omega;\mathbb{R}^m)$. Then
\begin{equation}\label{pri:4.7}
 \Big\{\frac{1}{\varepsilon}\int_{\Omega\setminus\Sigma_{\varepsilon}}\big(|u_\varepsilon|^2+|\nabla u_\varepsilon|^2\big) dx \Big\}^{\frac{1}{2}} \leq C\big\{\|F\|_{L^p(\Omega)} + \|g\|_{L^2(\partial\Omega)}\big\},
\end{equation}
where $C$ depends only on $\mu,\kappa,m,d,p$ and $\Omega$.
\end{thm}

\begin{pf} Without loss of generality we may assume that $\|F\|_{L^p(\Omega)} + \|g\|_{L^2(\partial\Omega)} = 1$.
Let
\begin{equation*}
w_\varepsilon = u_\varepsilon - u_0 -\varepsilon\chi_{0,\varepsilon}S^2_\varepsilon(\psi_{4\varepsilon}u_0)
-\varepsilon\chi_{k,\varepsilon}S^2_\varepsilon(\psi_{4\varepsilon}\nabla_k u_0)
\end{equation*}
Hence, it follows from $\eqref{pri:4.6}$ and $\eqref{f:4.29}$ that
\begin{equation*}
 \|u_\varepsilon\|_{H^1(\Omega\setminus\Sigma_{\varepsilon})}
 \leq \|u_0\|_{H^1(\Omega\setminus\Sigma_{\varepsilon})}
 + \|w_\varepsilon\|_{H^1(\Omega)}\leq C\varepsilon^{\frac{1}{2}}.
\end{equation*}
Note that $S^2_\varepsilon(\psi_{4\varepsilon}u_0)$ and $S^2_\varepsilon(\psi_{4\varepsilon}\nabla_k u_0)$
is supported in $\Sigma_\varepsilon$, therefore their $H^1$-norms vanish on $\Omega\setminus\Sigma_\varepsilon$. The proof is completed.
\qed
\end{pf}

\begin{remark}\label{re:5.2}
\emph{ If we additionally assume $A\in C^1(\mathbb{R}^d)$ in Theorem $\ref{thm:4.1}$, and $F =0$.
 Then, the estimate $\eqref{pri:4.7}$ in fact leads to (the Rellich estimate)
 $\|\nabla u_\varepsilon\|_{L^2(\partial\Omega)} \leq C\|g\|_{L^2(\partial\Omega)}$,
 where $C$ is independent of $\varepsilon$.
 The proof will be given in another place. Also, by referring to \cite[Remark 3.1]{SZW12}, the reader can prove it
 without real difficulties.}
\end{remark}

\begin{flushleft}
\textbf{Proof of Theorem \ref{thm:1.3}}\textbf{.}\quad
(i) With the help of the preceding Lemma $\ref{lemma:6.5}$ we can now prove the first part of Theorem $\ref{thm:1.3}$.
Note that
\begin{equation*}
 \begin{aligned}
\|u_\varepsilon - u_0\|_{L^2(\Omega)}
&\leq \|u_\varepsilon-u_0-\varepsilon\chi_{0,\varepsilon}S_\varepsilon(\psi_{4\varepsilon}\myu{u})
-\varepsilon\chi_{j,\varepsilon}S_\varepsilon(\psi_{4\varepsilon}\nabla_j\myu{u})\|_{L^2(\Omega)} \\
&+\varepsilon\|\chi_{0,\varepsilon}S_\varepsilon(\psi_{4\varepsilon}\myu{u})\|_{L^2(\Omega)}
+\varepsilon\|\chi_{j,\varepsilon}S_\varepsilon(\psi_{4\varepsilon}\nabla_j\myu{u})\|_{L^2(\Omega)}\\
&\leq C\varepsilon\big\{\|F\|_{L^2(\Omega)}+\|g\|_{B^{1/2,2}(\partial\Omega)}\big\} + C\varepsilon\|u_0\|_{H^1(\Omega)}\\
&\leq C\varepsilon\big\{\|F\|_{L^2(\Omega)}+\|g\|_{B^{1/2,2}(\partial\Omega)}\big\},
 \end{aligned}
\end{equation*}
\end{flushleft}
where we use the estimates $\eqref{pri:6.1}$ and $\eqref{pri:6.7}$ in the second inequality. In the last one, we mention that
$B^{1/2,2}(\partial\Omega)\subset B^{-1/2,2}(\partial\Omega)$.

(ii) The second part of Theorem $\ref{thm:1.3}$ follows from Lemma $\ref{lemma:4.2}$. As shown before, we have
\begin{equation*}
 \begin{aligned}
\|u_\varepsilon - u_0\|_{L^2(\Omega)}
&\leq \|u_\varepsilon-u_0-\varepsilon\chi_{0,\varepsilon}S_\varepsilon^2(\psi_{2\varepsilon}u_0)
-\varepsilon\chi_{j,\varepsilon}S_\varepsilon^2(\psi_{2\varepsilon}\nabla_ju_0)\|_{L^2(\Omega)} \\
&+\varepsilon\|\chi_{0,\varepsilon}S_\varepsilon^2(\psi_{2\varepsilon}u_0)\|_{L^2(\Omega)}
+\varepsilon\|\chi_{j,\varepsilon}S_\varepsilon^2(\psi_{2\varepsilon}\nabla_ju_0)\|_{L^2(\Omega)}\\
&\leq C\varepsilon^{\frac{1}{2}}\big\{\|F\|_{L^{\frac{2d}{d+1}}(\Omega)}+\|g\|_{L^{2}(\partial\Omega)}\big\}
+ C\varepsilon\big\{\|S_\varepsilon(\psi_{2\varepsilon}u_0)\|_{L^2(\mathbb{R}^d)}
+ \|S_\varepsilon(\psi_{2\varepsilon}\nabla u_0)\|_{L^2(\mathbb{R}^d)}\big\}\\
&\leq C\varepsilon^{\frac{1}{2}}\left\{\|F\|_{L^{\frac{2d}{d+1}}(\Omega)}+\|g\|_{L^{2}(\partial\Omega)}
+\|u_0\|_{L^{\frac{2d}{d+1}}(\Omega)}+\|\nabla u_0\|_{L^{\frac{2d}{d+1}}(\Omega)}\right\}\\
&\leq C\varepsilon^{\frac{1}{2}}\left\{\|F\|_{L^{\frac{2d}{d+1}}(\Omega)}+\|g\|_{L^{2}(\partial\Omega)}\right\}.
 \end{aligned}
\end{equation*}
In the second inequality, we employ the estimates $\eqref{pri:6.1}$ and $\eqref{pri:4.6}$.
In the third one, we use the estimate $\eqref{pri:6.3}$. For the last one, we note that
$\|u_0\|_{W^{1,p}(\Omega)}\leq C\|u_0\|_{H^1(\Omega)}\leq C\{\|F\|_{L^p(\Omega)}+\|g\|_{L^2(\partial\Omega)}\}$, where
$p=\frac{2d}{d+1}$.

By duality argument, we will prove the following estimate
\begin{equation}\label{pri:6.8}
\|u_\varepsilon - u_0\|_{L^\frac{2d}{d-1}(\Omega)}\leq C\varepsilon\|u_0\|_{H^2(\Omega)},
\end{equation}

Let $w_\varepsilon$ be defined in $\eqref{f:6.31}$. For any $\Phi\in L^{\frac{2d}{d+1}}(\Omega;\mathbb{R}^m)$,
let $\phi_\varepsilon,\phi_0$ be the solutions of $\eqref{pde:4.8}$ (the existence is given by Lemma $\ref{lemma:2.1}$).
Then it follows from Lemma $\ref{lemma:4.2}$ that
\begin{equation}\label{f:6.32}
\big\|\phi_\varepsilon - \phi_0 -\varepsilon\chi_{0,\varepsilon}^*S^2_\varepsilon(\psi_{12\varepsilon}\phi_0)
-\varepsilon\chi_{k,\varepsilon}^*S^2_\varepsilon(\psi_{12\varepsilon}\nabla_k \phi_0)\big\|_{H^1(\Omega)}
\leq C\varepsilon^{\frac{1}{2}}\|\Phi\|_{L^{\frac{2d}{d+1}}(\Omega)},
\end{equation}
where $\chi_k^*$ are correctors associated with $\mathcal{L}_\varepsilon^*$ for $k=0,\cdots,d$.

Proceeding as in the proof of Lemma $\ref{lemma:6.5}$, we arrive at $\eqref{f:6.21}$, which is
\begin{equation}
\begin{aligned}\label{f:6.33}
\Big|\int_{\Omega}w_\varepsilon\Phi dx\Big|
\leq C\Big\{\|u_0\|_{H^1(\Omega\setminus\Sigma_{9\varepsilon})}\|\phi_\varepsilon\|_{H^1(\Omega\setminus\Sigma_{9\varepsilon})}
+ \varepsilon\|u_0\|_{H^2(\Omega)}\|\phi_\varepsilon\|_{H^1(\Omega)}\Big\}.
\end{aligned}
\end{equation}
Due to Lemma $\ref{lemma:2.1}$ and H\"older inequality,
it is easy to see
\begin{equation}\label{f:6.35}
\|\phi_\varepsilon\|_{H^1(\Omega)}\leq C\|\Phi\|_{L^{\frac{2d}{d+2}}(\Omega)}\leq C\|\Phi\|_{L^{\frac{2d}{d+1}}(\Omega)}.
\end{equation}

Set $\xi_\varepsilon = \phi_\varepsilon - \phi_0 -\varepsilon\chi_{0,\varepsilon}^*S^2_\varepsilon(\psi_{12\varepsilon}\phi_0)
-\varepsilon\chi_{k,\varepsilon}^*S^2_\varepsilon(\psi_{12\varepsilon}\nabla_k \phi_0)$, and then we have
\begin{equation}\label{f:6.34}
\begin{aligned}
\|\phi_\varepsilon\|_{H^1(\Omega\setminus\Sigma_{9\varepsilon})}
\leq \|\xi_\varepsilon\|_{H^1(\Omega\setminus\Sigma_{9\varepsilon})}
+ \|\phi_0\|_{H^1(\Omega\setminus\Sigma_{9\varepsilon})}
\leq C\varepsilon^{\frac{1}{2}}\|\Phi\|_{L^{\frac{2d}{d+1}}(\Omega)}
\end{aligned}
\end{equation}
where we use the estimates $\eqref{f:6.29}$ and $\eqref{f:6.32}$ in the second inequality.
Collecting $\eqref{f:6.33}$, $\eqref{f:6.35}$, $\eqref{f:6.34}$, we obtain
\begin{equation*}
\begin{aligned}
\Big|\int_{\Omega}w_\varepsilon\Phi dx\Big|
\leq C\varepsilon^{\frac{1}{2}}\|u_0\|_{H^1(\Omega\setminus\Sigma_{9\varepsilon})}\|\Phi\|_{L^{\frac{2d}{d+1}}(\Omega)}
+ C\varepsilon\|u_0\|_{H^2(\Omega)}\|\Phi\|_{L^{\frac{2d}{d+1}}(\Omega)}
\leq C\varepsilon\|u_0\|_{H^2(\Omega)}\|\Phi\|_{L^{\frac{2d}{d+1}}(\Omega)},
\end{aligned}
\end{equation*}
where we use the estimate $\eqref{f:6.24}$ in the last inequality. This implies
\begin{equation*}
 \|w_\varepsilon\|_{L^{\frac{2d}{d-1}}(\Omega)} \leq C\varepsilon\|u_0\|_{H^2(\Omega)}.
\end{equation*}
and by the same procedure as we did in (i), it is not hard to derive the estimate $\eqref{pri:6.8}$,
where we need to employ the estimate $\eqref{pri:6.1}$ for $p=\frac{2d}{d-1}$,
and $\|\chi_{k}\|_{L^p(Y)}
\leq C\|\chi_k\|_{H^1(Y)}$
with $k=0,\cdots,d$, (due to Sobolev embedding theorem). The details are left to readers, and we complete the proof.
\qed

\begin{center}
\textbf{Acknowledgements}
\end{center}

The author wishes to express his sincere appreciation to Professor Zhongwei Shen for his constant guidance and encouragement. Without his
illuminating instruction, this paper could not have reached its present form.
The author is greatly indebted to the referee for many useful comments and helpful suggestions.
The author also wants to express his heartfelt gratitude to Professor Peihao Zhao, who led him into the world of mathematics.
This work was supported by the National Natural Science Foundation of China (Grant NO.11471147), and in part by
the Chinese Scholar Council (File No. 201306180043).




\clearpage
\end{document}